\documentclass[12pt,a4paper]{article}
\title
{\bf Finite order automorphisms and real forms 
of affine Kac-Moody algebras in the smooth and algebraic category
}
\author{Ernst Heintze and Christian Gro{\ss}} 
\usepackage{amsmath}
\usepackage{amsthm}
\usepackage{amsfonts}
\usepackage{amssymb}
\usepackage[T1]{fontenc}
\usepackage[ansinew]{inputenc}
\usepackage{latexsym}
\usepackage{array}
\newtheorem{theorem}{Theorem}[section]
\newtheorem{lemma}[theorem]{Lemma}
\newtheorem{definition}[theorem]{Definition}
\newtheorem{remark}[theorem]{Remark}
\numberwithin{equation}{section}
\newtheorem{corollary}[theorem]{Corollary}
\newtheorem{proposition}[theorem]{Proposition}

\def\D{\mathbb{D}}
\def\N{\mathbb{N}}
\def\R{\mathbb{R}}
\def\Q{\mathbb{Q}}
\def\C{\mathbb{C}}

\def\F{\mathbb{F}}
\def\Z{\mathbb{Z}}

\newcommand\J{\mathfrak J}

\newcommand\U{\mathcal U}

\newcommand\h{\mathfrak h}
\newcommand\g{\mathfrak g}
\renewcommand\u{\mathfrak u}
\newcommand\p{\mathfrak p}

\renewcommand\k{\mathfrak k}
\renewcommand\a{\mathfrak a}
\newcommand\s{\mathfrak s}
\renewcommand\l{\mathfrak l}
\newcommand\e{\mathfrak e}
\newcommand\f{\mathfrak f}
\renewcommand\b{\mathfrak b}
\renewcommand\c{\mathfrak c}
\renewcommand\d{\mathfrak d}
\newcommand\su{{\mathfrak s} {\mathfrak u}}
\newcommand\so{{\mathfrak s} {\mathfrak o}}
\renewcommand\sp{{\mathfrak s} {\mathfrak p}}

\newcommand\KC{\mathcal K}
\newcommand\PC{\mathcal P}
\newcommand\GC{\mathcal G}

\def\Hom{\mathrm{Hom}}
\def\Aut{\mathrm{Aut}}
\def\End{\mathrm{End}}
\def\Int{\mathrm{Int}}

\def\alg{\mathrm{alg}}
\def\ord{\mathrm{ord}}
\def\supp{\mathrm{supp}}
\def\rank{\mathrm{rank}}
\def\per{\mathrm{per}}
\def\ad{\mathrm{ad}}
\def\Ad{\mathrm{Ad}}
\def\id{\mathrm{id}}

\date{}

\begin{document}
\vspace {-1 cm}
\maketitle
%\markboth{\today}{\today}

\begin{abstract}\setlength{\baselineskip}{1.5 em}

Let $\g$ be a real or complex (finite dimensional) simple Lie algebra and $\sigma\in\Aut\g$. We study automorphisms of the twisted loop algebra $L(\g,\sigma)$ of smooth $\sigma$-periodic maps from $\R$ to $\g$ as well as of the ''smooth`` affine Kac-Moody algebra $\hat L(\g,\sigma)$, which is a $2$-dimensional extension of $L(\g,\sigma)$. It turns out that these automorphisms which either preserve or reverse  the orientation of loops, and are correspondingly called to be of first and second kind, can be described essentially by curves of automorphisms of $\g$. If the order of the automorphisms is finite then the corresponding curves in $\Aut\g$ allow to define certain invariants and these turn out to parametrize the conjugacy classes of the automorphisms. If their order is 2 (and $\g$ is either compact or complex) we carry this out in detail and deduce a complete classification of involutions and real forms (which correspond to conjugate linear involutions) of smooth affine Kac-Moody algebras. The resulting classification can be seen as an extension of Cartan's classification of symmetric spaces, i.e.~of involutions on $\g$. For example conjugacy classes of involutions of the second kind on $\hat L(\g,\sigma)$ are classified by equivalence classes of pairs $(\varrho_+,\varrho_-)$ where $\varrho_\pm\in\Aut\g$ are involutions or the identity, and $\varrho_-\varrho_+$ is conjugate to $\sigma$ in $\Aut\g/\Int\g$. If $\g$ is compact then conjugate linear extensions of involutions (and the identity) from $\hat L(\g,\sigma)$ to conjugate linear involutions on $\hat L(\g_\C,\sigma_\C)$ yield a bijection between their conjugacy classes and this gives existence and uniqueness of Cartan decompositions of real forms of complex smooth affine Kac-Moody algebras.

The affine Kac-Moody algebras introduced by Kac and Moody are isomorphic to a $2$-dimensional extension of the algebra of twisted loops in $\g$ whose Fourier expansion is finite (assuming $\sigma$ to be of finite order). We show that our methods work equally well also in this case when combined with a basic result of Levstein and lead essentially to the same results.\\

\noindent key words: Kac-Moody algebras, loop algebras, automorphisms of finite order, real forms, Cartan decompositions
\end{abstract}

\section{Introduction}
We study in this paper automorphisms of finite order and real forms of ``smooth'' affine Kac-Moody algebras, that is of certain extensions of the algebra of smooth twisted loops (for a precise definition see below). These objects have been considered already extensively in the algebraic category, where the loops are assumed to have {\sl finite} Fourier expansion (\cite{BP}, \cite{Kob}, \cite{Lev},  \cite{Bau}, \cite{BR}, \cite{Rou1}, \cite{Rou2}, \cite{Rou3}, \cite{Cor1}, \cite{Cor2}, \cite{Cor3}, \cite{And}, \cite{B$_3$R}, \cite{KW}, \cite{Bat}, \cite{JZ}, \cite{BMR}, \cite{BMR'}). In particular involutions and real forms have finally been classified in the algebraic case in \cite{B$_3$R} and \cite{BMR}.

Our approach is very different, much more elementary and direct. It does not use the structure theory of Kac-Moody algebras but rather reduces the problems as fast as possible to the finite dimensional case. Interesting enough, it also works in the algebraic setting and seems to give even there more complete answers and new insights. For example it turns out that involutions and real forms of affine Kac-Moody algebras are either in close connection with hyperpolar actions on compact Lie groups or else with the group $\pi_0((\Aut\g)^\varrho)$ of connected components of the centralizer of an involution~$\varrho$ in the group of automorphism of a simple Lie algebra~$\g$ (cf.~Chapter 6).

To describe our approach and results in more detail, let $\g$ be a finite dimensional simple Lie algebra over $\F:=\R$ or $\C$ and $\sigma\in\Aut\g$ be an arbitrary automorphism, not necessarily of finite order. We then call
\[
L(\g,\sigma):=\{u:\R\to\g\mid u(t+2\pi)=\sigma u(t), u\in C^\infty\}
\]
a twisted loop algebra and $\hat L(\g,\sigma):=L(\g,\sigma)+\F c+\F d$\, a (smooth) affine Kac-Moody algebra. Here $c$ lies in the center, $d$ acts on the loops as derivation and the bracket between two loops is the pointwise bracket plus a certain multiple of $c$ (cf.~Chapter 3).

An isomorphism $\hat\varphi:\hat L(\g,\sigma)\to\hat L(\tilde\g,\tilde\sigma)$ between two such algebras induces an isomorphism $\varphi:L(\g,\sigma)\to L(\tilde\g,\tilde\sigma)$ between the loop algebras. The isomorphisms $\varphi$ or $\hat\varphi$ are called {\sl standard} if $\varphi u(t)=\varphi_t(u(\lambda(t)))$ where $\lambda:\R\to\R$ is a diffeomorphism and $\varphi_t:\g\to\tilde\g$ is a smooth curve of isomorphisms. Our first main result (Theorem~\ref{2.10}) says that isomorphisms between loop algebras are always standard. The essential point of the proof consists in showing that for each fixed $t_0\in\R$ there exists an $s_0\in\R$ such that $\varphi(fu)(t_0)=f(s_0)\cdot\varphi u(t_0)$ for all $u\in L(\g,\sigma)$ and all $2\pi$-periodic smooth functions~$f$. This follows by means of a classical theorem of Burnside. In order that $\varphi u(t)$ is $\tilde\sigma$-periodic for all $u$ one necessarily has $\lambda(t+2\pi)=\lambda(t)+\epsilon 2\pi$ for some $\epsilon\in\{\pm 1\}$ and $\varphi_{t+2\pi}=\tilde\sigma\varphi_t\sigma^{-\epsilon}$. The isomorphism $\varphi$ (as well as $\hat\varphi$ if it induces $\varphi$) is called of first kind if $\epsilon=1$ and of second kind if $\epsilon=-1$, i.e.~if $\lambda$ is orientation preserving, resp. reversing. Conversely, given $\lambda$ and $\varphi_t$ which satisfy the above conditions, the mapping $\varphi:L(\g,\sigma)\to L(\tilde\g,\tilde\sigma)$ with $\varphi u(t)=\varphi_t(u(\lambda(t))$ defines an isomorphism. It extends to an isomorphism between the Kac-Moody algebras precisely if $\lambda'(t)$ is constant, i.e.~$\lambda(t)=\epsilon t+t_0$ for some $t_0\in\R$. Moreover the extension is almost unique and this implies that conjugacy classes of automorphisms of finite order on $\hat L(\g,\sigma)$ and $L(\g,\sigma)$ are in bijective correspondence. Thus the study of isomorphisms between affine Kac-Moody and loop algebras is essentially reduced to the study of curves of isomorphisms between the corresponding finite dimensional simple Lie algebras. This result offers the possibility to classify automorphisms of finite order up to conjugation in an elementary way, since one easily sees how $\varphi_t$ and $\lambda(t)$ change if an automorphism $\varphi$ of the form $\varphi u(t)=\varphi_t(u(\lambda(t)))$ is conjugated by another automorphism $\psi$ with $\psi u(t)=\psi_t(u(\mu(t)))$ (Lemma~\ref{3.10}). The problem that has to be solved then, is to extract an invariant out of $\varphi_t$ and $\lambda$ which does not change under these modifications and which determines the conjugacy class. This is done in Chapter~4 for automorphisms of the first kind and in Chapter~5 for automorphisms of the second kind. To this end we define for each $i\in\{1,2\}$ and each $q\in\N$ ($q$ even if $i=2$) so-called sets of invariants $\J_i^q(\g,\sigma)$ as follows:
\begin{eqnarray*}
&\J_1^q(\g,\sigma):=\{(p,\varrho,[\beta])\mid p\in\{0,1,\dots,q-1\},\ \varrho\in{\cal A}^r,\ \beta\in(\Aut\g)^\varrho,\ \hat L(\g,\varrho^l\beta^{q'})\cong\hat L(\g,\sigma)\}\ ,\\
&\J^{2q}_2(\g,\sigma):=\{[\varphi_+,\varphi_-]\mid\varphi_\pm\in\Aut\g,\ \varphi_+^2=\varphi_-^2,\ \ord(\varphi_\pm^2)=q,\ \hat L(\g,\varphi_-^{-1}\varphi_+)\cong\hat L(\g,\sigma)\}
\end{eqnarray*}
where $r,l$ and $q'$ are certain integers depending only on $p$ and $q$ (cf.~\ref{4.3}), ${\cal A}^r$ is a set of representatives of conjugacy classes of automorphisms of $\g$ of order $r$ and $[\beta]$ and $[\varphi_+,\varphi_-]$ denote equivalence classes with respect to some equivalence relation. Note that $\hat L(\g,\sigma)$ and $\hat L(\g,\tilde\sigma)$ are isomorphic if and only if $\sigma$ and $\tilde\sigma$ are conjugate in $\Aut\g/\Int\g$~(\ref{3.5}). We then define a mapping
$\Aut^q_i\hat L(\g,\sigma)\to \J_i^q(\g,\sigma)$ from the set of automorphisms of order $q$ and kind~$i$ to the set of invariants which is constant along conjugacy classes. Actually this mapping gives the same value if $\hat\varphi\in\Aut_i^q\hat L(\g,\sigma)$ is conjugated more generally by an {\sl isomorphism} $\hat\psi:\hat L(\g,\sigma)\to\hat L(\g,\tilde\sigma)$, in which case we call $\hat\psi\hat\varphi\hat\psi^{-1}$ to be {\sl quasiconjugate} to $\hat\varphi$. This remark is useful in proving surjectivity of the above mapping since it is fairly easy to find for each invariant a $\tilde\sigma\in\Aut\g$ with $\hat L(\g,\tilde\sigma)$ isomorphic to $\hat L(\g,\sigma)$ and a $\hat\varphi\in\Aut_i^q\hat L(\g,\tilde\sigma)$ with the given invariant. Moreover $\tilde\sigma$ and $\hat\varphi$ can be chosen in such a way that $\varphi u(t)=\varphi_0 u(\epsilon t+t_0)$ with $\varphi_0\in\Aut\g$ (constant). But the main point is to show that two automorphisms with the same invariant are conjugate (Theorems~\ref{4.11} and~\ref{5.7}) and that hence the sets of invariants parametrize the conjugacy classes. From the remark above it follows then that any automorphism of finite order of $\hat L(\g,\sigma)$ is quasiconjugate to one with $\varphi u(t)=\varphi_0(u(\epsilon t+t_0))$ where $\varphi_0$ is constant (and $\hat\varphi c=\epsilon c,\ \hat\varphi d=\epsilon d$). In \cite{HPTT} these special automorphisms had been studied and it had been asked whether any automorphism of finite order is {\sl conjugate} to such a special one. The answer is ``no'' in general as examples show (cf.~Remarks~\ref{4.14} and~\ref{5.11}) but ``yes'' if we allow to change $\sigma$, i.e.~to apply a quasiconjugation.

In Chapter~6 we specialize the above results to involutions i.e.~to automorphisms of order two. We classify these explicitly up to conjugation if $\g$ is compact or complex. The classification is in both cases the same and amounts to determine explicitly $\J^2_1(\g,\sigma)$ and $\J_2^2(\g,\sigma)$, which in turn follows from a refinement of E. Cartan's classification of involutions of $\g$. While one has to determine for each involution $\varrho\in\Aut\g$ representatives of the conjugacy classes of $\pi_0((\Aut\g)^\varrho)$ in the first case one has to determine pairs $(\varrho_+,\varrho_-)$ of $\varrho_\pm\in\Aut\g$ with $\varrho_\pm^2=\id$ up to a certain equivalence relation in the second case. The finite groups $\pi_0((\Aut\g)^\varrho)$ of connected components of the centralizer of $\varrho$ in $\Aut\g$ have already been computed by Cartan and Takeuchi and are listed e.g. in~\cite{Loo}.
An explicit determination of their conjugacy classes is relegated to Appendix A where also a careful discussion of these groups together with a simplified computation of them is given.

If $\g$ is complex then conjugate linear automorphisms of finite order of $\hat L(\g,\sigma)$, i.e.~automorphisms of the realification that anticommute with multiplication by $i=\sqrt{-1}$, can be treated as in Chapter 4. They are also standard and one can associate invariants to them that parametrize their conjugacy classes in a completely analogous way. The conjugate linear involutions are in bijection with real forms of $\hat L(\g,\sigma)$, which are precisely their fixed point sets. We denote by $\overline{\Aut}^2_i\hat L(\g,\sigma)$ the conjugate linear involutions of type $i=1$ or $2$ (according to whether the induced mapping between the loop algebras preserves or reverses orientation, respectively). If $\u$ is a $\sigma$-invariant compact real form of $\g$ then $\hat L(\u,\sigma)$ (and its images under an isomorphism )are called {\sl compact} real forms of $\hat L(\g,\sigma)$. The mappings
\[
\Aut^1_1\hat L(\u,\sigma)\cup\Aut^2_1\hat L(\u,\sigma)\to\overline{\Aut}^2_1\hat L(\g,\sigma)
\mbox{\quad and \quad}
\Aut_2^2\hat L(\u,\sigma)\to\overline{\Aut}^2_2\hat L(\g,\sigma)
\]
which map $\hat\varphi$ to its conjugate linear extension $\hat\varphi_\C\circ\hat\omega$ where $\hat\omega$ denotes complex conjugation with respect to $\hat L(\u,\sigma)$, induce a bijection between conjugacy classes. This fact that follows easily by inspecting the corresponding mappings between the sets of invariants. In particular the (equivalence classes of) non compact real forms are in bijection with (conjugacy classes of) involutions of the compact real form, like in finite dimensions. If $\hat L(\u,\sigma)={\cal K}+{\cal P}$ is an eigenspace decomposition of an involution then ${\cal K}+i{\cal P}$ is a noncompact real form of $\hat L(\g,\sigma)$, and in this way all real forms are obtained. Moreover, each real form has therefore a Cartan decomposition ${\cal K}+\tilde{\cal P}$ with ${\cal K}+i\tilde{\cal P}$ compact and this is unique up to conjugation as the above mappings are injective on the set of conjugacy classes.

In Chapter 8, the last chapter of the paper we carry over our methods from the smooth to the algebraic setting and prove that also in this case automorphisms of finite order and real forms are parametrized by the same invariants as in the $C^\infty$-case. The arguments are similar but need at several points modifications. For example not all isomorphisms are standard in the algebraic case. We let
\[
L_\alg(\g,\sigma)=\{u\in L(\g,\sigma)\mid u(t)={\sum\limits_{\vert n\vert\le N}}u_ne^{int/l},N\in\N, u_n\in\g_\C\}
\]
where $\sigma$ is of finite order, $\sigma^l=\id$, and
\[
\hat L_\alg(\g,\sigma)=L_\alg(\g,\sigma)\oplus\F c\oplus\F d\ .
\]
It then follows that automorphisms of $L_\alg(\g,\sigma)$ are compositions of standard automorphisms with automorphisms $\tau_r$ which map $\sum u_ne^{int/l}$ to $\sum u_nr^{n/l}e^{int/l}$ where $r>0$. But the main difficulty is to show that two automorphisms of finite order with the same invariants (which in spite of the $\tau_r$ can be defined as in the $C^\infty$-case) are conjugate. It is at this point where we use a basic result of F.~Levstein \cite{Lev} which says that automorphisms of finite order of $\hat L_\alg(\g,\sigma)$ leave some Cartan algebra invariant. This implies that any automorphism of finite order is conjugate to a very special one and after a further quasiconjugation in fact to one of the form $\varphi u(t)=\varphi_0(u(\epsilon t+t_0))$ where $\varphi_0\in\Aut\g$ is constant. The conjugacy problem can then be solved by using certain hyperpolar actions on compact Lie groups.

The results of this paper have been announced in \cite{Hei1} and \cite{Hei2}. Most of them had been obtained many years ago, but it took us some time to fill in all details. 

\section{Isomorphisms between smooth loop algebras}

Let $\g$ be a finite dimensional Lie algebra over $\F=\R$ or $\C$ and $\sigma\in\Aut\g$ (not necessarily of finite order). Then the (smooth, twisted) loop algebra
\[
L(\g,\sigma):=\{u:\R\to\g\mid u(t+2\pi)=\sigma u(t)\ \forall\ t,\ u\in C^\infty\}
\]
is a Lie algebra with pointwise bracket
\[
[u,v]_0(t):=[u(t),v(t)]\ .
\]
$L(\g):=L(\g,id)$ is also called the untwisted algebra. 

\begin{remark}
\rm One may weaken the regularity assumption in the definition of $L(\g,\sigma)$ and consider e.g.
\[
L_k(\g,\sigma):=\{u:\R\to\g\mid u(t+2\pi)=\sigma u(t),\ u \hbox{ locally of Sobolev class }H^k\}
\]
for any $k\ge 1$. Although $[L_k,L_k]\subset L_{k-1}$ and thus the bracket is contained in $L_k$ only  after restriction it to dense subspaces like $L_k\times L_{k+1}$, the results of this paper nevertheless go through also for $L_k(\g,\sigma)$ without difficulties.

Another regularity class of interest is the class of algebraic loops
which are given by finite Laurent series, assuming $\sigma$ to be of
finite order. The corresponding algebra
$L_\alg(\g,\sigma)$ and its automorphisms will be studied
in Chapter 8.
\end{remark}

A homomorphism $\varphi:L(\g,\sigma)\to L(\tilde\g,\tilde \sigma)$ between to loop algebras is only supposed to be $\F$-linear and to preserve brackets, no continuity assumptions are made. Simple examples of homomorphisms are mappings $\varphi:L(\g,\sigma)\to L(\tilde\g,\tilde\sigma)$ with
$(\varphi u)(t)=\varphi_t(u(\lambda(t)))$ where $\varphi_t:\g\to\tilde\g$ are homomorphisms and $\lambda:\R\to\R$ is a function such that $t\mapsto \varphi_t$ and $\lambda$ are smooth $(= C^\infty)$.

\begin{definition}
A homomorphism $\varphi:L(\g,\sigma)\to L(\tilde \g,\tilde\sigma)$ is called \textit{standard} if it is of the above form
\[
(\varphi u)(t)=\varphi_t(u(\lambda(t)))\ .
\]
\end{definition}

The main goal of this chapter is to show that all isomorphisms $L(\g,\sigma)\to L(\tilde\g,\tilde\sigma)$ are standard if $\tilde\g$ is simple.

The following result is obvious.

\begin{lemma}\label{2.3}
Let $t_0\in\R$ and $I$ be an open interval around $t_0$. Then there exists for each smooth function $u:I\to\g$ an $\tilde u\in L(\g,\sigma)$ with $\tilde u(t)=u(t)$ in a neighborhood of $t_0$. In particular the evaluation map $L(\g,\sigma)\to\g,\ u\mapsto u(t_0)$, is surjective. \qed
\end{lemma}

The assumption that $\varphi_t$ and $\lambda_t$ are smooth in the definition of a standard homomorphism can be almost deleted.

\begin{lemma}\label{Lemma I}
Let $\varphi:L(\g,\sigma)\to L(\tilde\g,\tilde\sigma)$ be a surjective homomorphism. If $\varphi$ is of the form $\varphi u(t)=\varphi_t(u(\lambda(t)))$ for some homomorphisms $\varphi_t:\g\to\tilde\g$ and some function $\lambda:\R\to\R$ then $\varphi$ is standard, that is there exist $\tilde\varphi_t$ and $\tilde\lambda(t)$ depending smoothly on $t$ with $\varphi u(t)=\tilde\varphi_t(u(\tilde\lambda(t)))$ for all $u\in L(\g,\sigma)$.
\end{lemma}

\begin{proof}
By assumption $\varphi(f\cdot u)=(f\circ\lambda)\cdot\varphi u$ for all $2\pi$-periodic smooth $f:\R\to\F$ and all $u\in L(\g,\sigma)$. Since there exists for each $t_0\in\R$\ \, a\ \, $u\in L(\g,\sigma)$ with $\varphi u(t_0)\ne 0\ ,\ f\circ\lambda$ is smooth. In particular $e^{i\lambda}:\R\to S^1$ is smooth and thus has a smooth lift $\tilde\lambda:\R\to\R$. Hence $e^{i\lambda}=e^{i\tilde\lambda}$ and $\lambda(t)-\tilde\lambda(t)=2\pi k_t$ with $k_t\in\Z$. Let $\tilde\varphi_t:=\varphi_t\sigma^{k_t}$. Then $\tilde\varphi_tu(\tilde\lambda(t))=\varphi_tu(\lambda(t))=\varphi u(t)$. Moreover $t\mapsto\tilde\varphi_t$ is smooth as $u$ can be chosen to be locally constant by \ref{2.3} and $t\mapsto\tilde\varphi_t u(\tilde\lambda(t))$ is smooth.
\end{proof}

\begin{remark}\label{2.5}
\rm The representation of a standard homomorphism as $\varphi u(t)=\varphi_tu(\lambda(t))$ with  $\varphi_t$ and $\lambda$ smooth is still not unique. In fact, $\varphi_t(u(\lambda(t)))=\tilde\varphi_t(u(\tilde\lambda(t)))$ (with $\tilde\varphi_t$ and $\tilde\lambda$ also smooth) if and only if there exists $k\in\Z$ such that $\tilde\varphi_t=\varphi_t\sigma^k$ and $\tilde\lambda(t)=\lambda(t)-2k\pi$.
\end{remark}

Let $C^\infty_{\per}(\R,\F):=\{f:\R\to\F\mid f\in C^\infty,\ f(t+2\pi)=f(t)\ \forall\ t\in\R\}$ be the algebra of $2\pi$-periodic smooth functions.

\begin{lemma}\label{2.6}
A surjective homomorphism $\varphi:L(\g,\sigma)\to L(\tilde\g,\tilde\sigma)$ is standard if and only if there exists a (not necessarily smooth) mapping $\lambda:\R\to\R$ with
\[
\varphi(fu)=(f\circ\lambda)\cdot\varphi(u)
\]
for all $f\in C^\infty_{\rm per}(\R,\F)$ and $u\in L(\g,\sigma)$.
\end{lemma}

\begin{proof}
The condition is clearly necessary. Conversely, if satisfied let $t_0\in\R$ be fixed and $s_0:=\lambda(t_0)$. By \ref{Lemma I} it suffices to show that $\varphi u(t_0)$ depends only on $u(s_0)$ or equivalently that $\varphi u(t_0)=0$ if $u(s_0)=0$.

If $u$ vanishes even in a neighborhood of $s_0$ then there exists $f\in C^\infty_{\per}(\R,\F)$ with $f\cdot u\equiv 0$ and $f(s_0)=1$. Hence $0=\varphi(fu)=(f\circ\lambda)\cdot\varphi(u)$ and thus $\varphi u(t_0)=0$. Therefore $\varphi u(t_0)$ depends only on $u(t)$ in a neighborhood of $s_0$ for any $u\in L(\g,\sigma)$.

Finally let $u(s_0)=0$. By means of \ref{2.3} there exist $u_1,\dots,u_n$ in $L(\g,\sigma)$ such that $u_1(s),\dots, u_n(s)$ are a basis of $\g$ for all $s$ close to $s_0$. Thus $u(s)=\Sigma f_i(s)u(s)$ in a neighborhood of $s_0$ for some $f_i\in C^\infty_{\per}(\R,\F)$ and $f_i(s_0)=0$. Hence $\varphi u(t_0)=(\Sigma(f_i\circ\lambda)\cdot\varphi(u_i))(t_0)=0$.
\end{proof}

The next lemma is needed to extend the main result from $\F=\C$ to $\F=\R$.

\begin{lemma}\label{2.7}
Let $\g=\g_+\oplus\g_-$ be a decomposition of $\g$ into two ideals and $\sigma\in\Aut\g$.
\begin{enumerate}
\item[(i)] If $\sigma$ leaves the ideals invariant then
\[
L(\g_+\oplus\g_-,\sigma)\cong L(\g_+,\sigma_+)\oplus L(\g_-, \sigma_-)
\]
where $\sigma_\pm$ denote the restrictions. The isomorphism is given by $u\mapsto (u_+,u_-)$ if $u(t)=u_+(t)+u_-(t)$.
\item[(ii)] If $\sigma$ interchanges the ideals then
\[
L(\g_+\oplus\g_-,\sigma)\cong L(\g_+,\sigma_+^2)\cong L(\g_-,\sigma_-^2)
\]
The isomorphisms are given by $u_+(t)+u_-(t)\mapsto u_\pm(2t)$.
\end{enumerate}
\end{lemma}

The proofs are almost straightforward. Note that in case (ii) $u_+(t)+u_-(t)\in L(\g_+\oplus\g_-,\sigma)$ implies that $u_\pm(t+4\pi)=\sigma^2 u_\pm(t)$ and $u_\pm(t)=\sigma(u_\mp(t-2\pi))$. In particular  $u_\pm(2t)$ are in $L(\g_\pm,\sigma^2_\pm)$ and determine each other. \qed

\begin{lemma}\label{2.8}
Let $\alpha:C^\infty_{\per}(\R,\C)\to\C$ be a (not necessarily continuous) homomorphism of algebras which does not vanish identically. Then there exists $s_0\in\R$ with $\alpha(f)=f(s_0)$ for all $f\in C^\infty_{\per}(\R,\C)$.
\end{lemma}

\begin{proof}
Since $\alpha$ is not identically zero, $\alpha(\tilde 1)=1$ where $\tilde 1$ denotes the function $f(t)\equiv 1$. We first show that $\alpha$ is continuous with respect to the sup-norm $\Vert.\Vert$, more precisely that $\vert\alpha(f)\vert\le\Vert f\Vert$ for all $f$. In fact, if $\vert\alpha(f)\vert>\Vert f\Vert$ for some $f$ then $g:=\alpha(f)\cdot\tilde1-f$ vanishes nowhere. Hence $1/g\in C^\infty_{\per}(\R,\C)$ and thus $\alpha(g)\ne 0$ because of $\alpha(g)\cdot\alpha(1/g)=1$ in contradiction to $\alpha(g)=\alpha(\alpha(f)\tilde 1-f)=0$.

Let $z_0:=\alpha(e^{it})$. Then $\vert z_0^{\pm 1}\vert=\vert\alpha(e^{\pm it})\vert\le 1$ and hence $z_0=e^{is_0}$ for some $s_0\in\R$. Moreover $\alpha({\sum\limits^N_{-N}}a_ne^{int})={\sum\limits^N_{-N}}a_ne^{ins_0}$, that is $\alpha(f)=f(s_0)$ if $f$ has a finite Fourier expansion. By continuity $\alpha(f)=f(s_0)$ then holds for all $f\in C^\infty_{\per}(\R,\C)$.
\end{proof}

\begin{theorem}\label{2.9}
Let $\varphi:L(\g,\sigma)\to L(\tilde\g,\tilde\sigma)$ be a surjective homomorphism with $\tilde\g$ simple. Then $\varphi$ is standard.
\end{theorem}

\begin{proof}
\begin{enumerate}
\item [(i)] Let $\F=\C$.

By Lemma \ref{2.6} it suffices to find for each $t_0\in\R$ an $s_0\in\R$ with
\[
\varphi(fu)(t_0)=f(s_0)\varphi u(t_0)
\]
for all $u\in L(\g,\sigma)$ and all $f\in C^\infty_{\per}(\R,\C)$. Let $u$ and $f$ be fixed and $a:=\varphi(fu)(t_0),\ b:=\varphi u(t_0)$ and $x_i:=\varphi u_i(t_0)$ where $u_1,\ldots,u_k\in L(\g,\sigma)$ and $k\in\N$ are arbitrary. Then $\ad a\ \ad x_1\cdots \ad x_k\ \ad b=\ad b\ \ad x_1\cdots \ad x_k\ \ad a$ since $[fu,[u_1,\ldots,[u_k,[u,v]_0\cdots]_0=[u,[u_1,\ldots,[u_k,[fu,v]_0\cdots]_0$ for all $v\in L(\g,\sigma)$. The associated subalgebra of $End\tilde\g$ spanned by the products $\ad x_1\cdots \ad x_k$ acts irreducibly on $\tilde\g$ because $\tilde\g$ is simple. By a theorem of Burnside  (cf. \cite{Lan}, XXVII 3.3 and 3.4) it therefore coincides with $\End\tilde\g$. Thus $\ad a\ X\ \ad b= \ad b\ X\ \ad a$ for all $X\in\End\tilde \g$ implying that $\ad a$ and $\ad b$ are linearly dependent. Hence also $a$ and $b$ are linearly dependent as $\tilde \g$ has no center. In particular if $u_0\in L(\g,\sigma)$ satisfies $\varphi u_0(t_0)\ne 0$ then $\varphi (f u_0)(t_0)=\alpha (f)\varphi u_0(t_0)$ for all $f\in C^\infty_{\per}(\R,\C)$, where $\alpha: C^\infty_{\per}(\R,\C)\to\C$ is linear with $\alpha(\tilde 1)=1$. We claim that
\[
\varphi (fu)(t_0)=\alpha(f)\varphi u(t_0)
\]
holds in fact for all $u\in L(\g,\sigma)$, whence $\alpha$ is an algebra homomorphism and the theorem follows from \ref{2.8} and \ref{2.6}.

To prove the claim we consider two cases. If $\varphi u(t_0)$ and $\varphi u_0(t_0)$ are linearly dependent then $\varphi(u-\lambda u_0)(t_0)=0$ for some $\lambda\in\C$ and thus $\varphi(f(u-\lambda u_0))(t_0)=0$ as $[\varphi(fv)(t_0),\ \varphi w(t_0)]=[\varphi v(t_0),\ \varphi (fw)(t_0)]=0$ for all $v,w\in L(\g,\sigma)$ with $\varphi v(t_0)=0$. From this the claim follows.

If $\varphi u(t_0)$ and $\varphi u_0(t_0)$ are linearly independent then $\varphi(f(u+u_0))(t_0)=\varphi(fu)(t_0)+\alpha(f)\varphi(u_0)(t_0)$ on one hand and a multiple of $\varphi(u+u_0)(t_0)$ on the other hand. From this again the claim follows.

\item[(ii)] If $\F=\R$ we consider the complexification $\varphi_\C:L(\g_\C,\sigma_\C)\to L(\tilde \g_\C,\tilde\sigma_\C)$ of $\varphi$. If $\tilde\g_\C$ is simple $\varphi_\C$ and hence also $\varphi$ are standard by (i). If $\tilde\g_\C$ is not simple then $\tilde\g$ has a complex structure $J$ and $\tilde\g_\C$ is the direct sum of the two ideals $\tilde\g_\pm:=\{X\pm iJX\mid X\in\tilde\g\}$ which are simple. $X\mapsto \frac{1}{2} (X+iJX)+\frac{1}{2}(X-iJX)$ defines an isomorphism between $\tilde\g_\C$ and $\tilde\g_+\oplus\tilde\g_-$. Either $\tilde\sigma_\C$ preserves or interchanges these ideals. In the first case Lemma \ref{2.7} yields the homomorphism
\[
L(\g_\C,\sigma_\C)\to L(\tilde\g_+,\tilde\sigma_+),\ u(t)\mapsto \frac{1}{2} (\varphi u(t)+iJ\varphi u(t))\ ,
\]
and in the second case
\[
L(\g_\C,\sigma_\C)\to L(\tilde\g_+,\tilde\sigma_+^2),\ u(t)\mapsto \frac{1}{2} (\varphi u(2t)+iJ\varphi u(2t))\ .
\]
Since these homomorphisms are surjective they are standard by (i) and it follows that also $\varphi$ is standard.
\end{enumerate}
\vspace{-27pt}
\end{proof}

Specializing to the case of isomorphims we can sharpen our results.

\begin{theorem}\label{2.10}
Let $\g,\tilde\g$ be simple.
\begin{enumerate}
\item[(i)] If $\varphi:L(\g,\sigma)\to L(\tilde\g,\tilde\sigma)$ is an isomorphism then there exist $\epsilon\in\{\pm 1\}$, a diffeomorphism $\lambda:\R\to\R$ with $\lambda(t+2\pi)=\lambda(t)+\epsilon 2\pi$, and a smooth curve of automorphisms $\varphi_t:\g\to\tilde\g$ with $\varphi_{t+2\pi}=\tilde\sigma\varphi_t\sigma^{-\epsilon}$ such that
\[
\varphi u(t)=\varphi_t u(\lambda(t))
\]
for all $u\in L(\g,\sigma)$.
\item[(ii)] Conversely, $\epsilon,\lambda$ and $\{\varphi_t\}$ as above define an isomorphism $\varphi:L(\g,\sigma)\to L(\tilde\g,\tilde\sigma)$ by $\varphi u(t)=\varphi_t(u(\lambda(t)))$.
\item[(iii)] $\tilde\epsilon, \tilde\lambda$ and $\{\tilde\varphi_t\}$ define the same isomorphism if and only if $\tilde\epsilon=\epsilon$ and there exists $k\in\Z$ with $\tilde\lambda(t)=\lambda(t)-2k\pi,\ \tilde\varphi_t=\varphi_t\sigma^k$.
\end{enumerate}
\end{theorem}

\begin{proof}
%\begin{enumerate}
%\item[(i)]
$\varphi$ as well as $\varphi^{-1}$ are standard by Theorem~\ref{2.9}, that is $\varphi u(t)=\varphi_tu(\lambda(t))$, $\varphi^{-1}u(t)=\psi_tu(\mu(t))$ for some smooth maps $\lambda,\mu\colon\R\to\R$ and some homomorphisms $\varphi_t\colon\g\to\nobreak\tilde\g$, $\psi_t\colon\tilde\g\to\g$ depending smoothly on $t$. Therefore $\varphi\circ\varphi^{-1}=\id$ is equivalent to $\varphi_t\psi_{\lambda(t)}v (\mu\circ\lambda(t))=v(t)$ for all $v\in L(\tilde\g,\tilde\sigma)$. Hence by \ref{2.5} there exists $m\in\Z$ with $\mu\circ\lambda=\id-2m\pi$ and $\varphi_t\psi_{\lambda(t)}=\tilde\sigma^m$. In particular the $\varphi_t$ are isomorphisms. Now $\varphi u\in L(\tilde\g,\tilde\sigma)$ for all $u\in L(\g,\sigma)$ is equivalent to $\varphi_{t+2\pi}u(\lambda(t+2\pi))=\tilde\sigma\varphi_u u(\lambda(t))$ for all $u$ and hence to $\lambda(t+2\pi)=\lambda(t)+2k\pi$ and $\varphi_{t+2\pi}=\tilde\sigma\varphi_t\sigma^{-k}$ for some $k\in\Z$. Similarly $\mu(t+2\pi)=\mu(t)+2l\pi$ for some $l\in\Z$. From $\mu\circ\lambda=\id-2m\pi$ we get $k\cdot l=1$. Hence $\epsilon:=k=l\in\{\pm 1\}$. This proves (i),  (ii) follows easily, and (iii) is a consequence of \ref{2.5}.
%\end{enumerate}
%\vspace{-12pt}
\end{proof}

\begin{remark}\rm
The $\lambda$ in the representation of $\varphi$ as $\varphi u(t)=\varphi_t(u(\lambda(t)))$ is orientation preserving (reversing) if and only if $\epsilon=1$ (resp.~$\epsilon=-1)$. In particular $\lambda$ determines $\epsilon$.
\end{remark}

\begin{definition}
Let $\g$ and $\tilde\g$ be simple. An isomorphism $\varphi:L(\g,\sigma)\to L(\tilde\g,\tilde\sigma)$ is called to be of first (resp.~second) kind if $\epsilon=1$ (resp.~$\epsilon=-1$).
\end{definition}

For example $\varphi u(t):=u(-t)$ defines an isomorphism of the second kind between $L(\g,\sigma)$ and $L(\g,\sigma^{-1})$. In particular $L(\g,\sigma)$ and $L(\g,\sigma^{-1})$ are isomorphic.

\begin{corollary}\label{2.13}
Let $\g$ and $\tilde\g$ be simple.
\begin{enumerate}
\item[(i)] If $L(\g,\sigma)$ and $L(\tilde\g,\tilde\sigma)$ are isomorphic then $\g$ and $\tilde\g$ are isomorphic.
\item[(ii)] $L(\g,\sigma)$ and $L(\g,\tilde\sigma)$ are isomorphic if and only if $\sigma$ and $\tilde\sigma$ are conjugate in $\Aut\g/\Int\g$, i.e.~$\tilde\sigma=\alpha\beta\sigma\beta^{-1}$ with $\alpha\in\Int\g$ and $\beta\in\Aut\g$.
\end{enumerate}
\end{corollary}

\begin{proof}
\begin{enumerate}
\item[(i)] follows directly from Corollary~\ref{2.10}.
\item[(ii)] If $\varphi:L(\g,\sigma)\to L(\g,\tilde\sigma)$ is an isomorphism then $\varphi u(t)=\varphi_t(u(\lambda(t)))$ for some smooth family $\{\varphi_t\}$ of automorphisms of $\g$ and $\varphi_{t+2\pi}=\tilde\sigma\varphi_t\sigma^{-\epsilon}$ for some $\epsilon\in\{\pm 1\}$. Hence $\tilde\sigma=\varphi_{t+2\pi}\sigma^\epsilon\varphi_t^{-1}$ and $\tilde\sigma$ is conjugate to $\sigma^\epsilon$ in $\Aut\g/\Int\g$. But each element in $\Aut\g/\Int\g$ is conjugate to its inverse as will be noted in the following remark. Conversely if $\sigma$ and $\tilde\sigma$ are conjugate in $\Aut\g/\Int\g$ there exist $\alpha\in \Aut\g$ and $\beta\in\Int\g$ with $\tilde\sigma=\alpha\beta\sigma\alpha^{-1}$. Since $\alpha$ and $\alpha\beta$ lie in the same connected component of  $\Aut\g$ there exists a smooth mapping $[0,2\pi]\to\Aut\g,\ t\mapsto\varphi_t$, with $\varphi_t\equiv \alpha$ near $t=0$ and $\varphi_t\equiv\alpha\beta$ near $t=2\pi$. It extends smoothly to all of $\R$ by $\varphi_{t+2k\pi}:=\tilde\sigma^k\varphi_t\sigma^{-k}$ for $k\in\Z$ and $t\in[0,2\pi]$ satisfying $\varphi_{t+2\pi}=\tilde\sigma\varphi_t\sigma^{-1}$. Thus $\varphi u(t):=\varphi_tu(t)$ yields an isomorphism between $L(\g,\sigma)$ and $L(\g,\tilde\sigma)$.
\end{enumerate}
\vspace{-27pt}
\end{proof}

\begin{remark}\label{2.14}\rm 
If $\g$ is simple and either compact or complex then $\Aut\g/\Int\g$ is isomorphic to $1,\Z_2$ or the symmetric group $S_3$ and any element in these groups is conjugate to its inverse. If $\g$ is real and simple but non compact and $\g=\k+\p$ is a Cartan decomposition with corresponding involution $\varrho$ then $\Aut\g/\Int\g\cong(\Aut\g)^\varrho/((\Aut\g)^\varrho)_0\cong(\Aut\g^{*})^\varrho/((\Aut\g^{*})^\varrho)_0$ where $\g^{*}=\k+i\p^{*}$ is the associated compact algebra with corresponding involution and $G_0$ denotes for any group $G$ the connected component containing the identity (cf. B.2 (i) of the Appendix). It is known that these groups are isomorphic to either $1,\Z_2,\Z_2\times\Z_2$, the dihedral group $\D_4$ or the symmetric group $S_4$, cf. \cite{Loo} or Appendix A. Also in these groups any element is conjugate to its inverse (cf. the discussion in Chapter 6).
\end{remark}

\begin{corollary}\label{2.15}
Let $\g$ be a simple Lie algebra over $\F=\R$ or $\C$ and $\sigma\in\Aut\g$. Then there exists an automorphism $\tilde\sigma$ of $\g$ of finite order with
\[
L(\g,\sigma)\cong L(\g,\tilde\sigma)\ .
\]
\end{corollary}

\begin{proof}
By \ref{2.13}, it suffices to show that any connected component of $\Aut\g$ contains an element of finite order. But this follows from the next lemma using that the compact group $(\Aut\g)^\varrho$ in \ref{2.14} meets every connected component of $\Aut\g$.
\end{proof}

\begin{lemma}\label{2.16}
Let $G$ be a Lie group and $H$ a compact subgroup  that meets every connected component of $G$. Then each connected component of $G$ contains an element of finite order.
\end{lemma}

\begin{proof}
Let $G_1$ be a connected component of $G$ and $h\in G_1\cap H$. Then the closure of $\{h^n\mid n\in\Z\}$ is compact and abelian and thus isomorphic to $T\times F$ where $T$ is a torus and $F$ is finite. After changing $h$ by an appropriate element of $T$ we get an element of finite order and this still lies in $G_1$.
\end{proof}

To study real forms of complex loop algebras we extend the main result to conjugate linear automorphisms $\varphi$, i.e. with $\varphi(iu)=-i\varphi(u)$ for all $u$.

\begin{corollary}\label{2.17}
Let $\g$ be a complex simple Lie algebra and $\sigma\in\Aut\g$. Then the conjugate linear automorphisms of $L(\g,\sigma)$ are precisely the $\varphi$ of the form $\varphi u(t)=\varphi_t(u(\lambda(t)))$ for all $u\in L(\g,\sigma)$ where $\lambda:\R\to\R$ is a diffeomorphism with $\lambda(t+2\pi)=\lambda(t)+\epsilon 2\pi$ for some $\epsilon\in\{\pm 1\}$ and the $\varphi_t$ are conjugate linear automorphisms of $\g$ with $\varphi_{t+2\pi}=\sigma\varphi_t\sigma^{-\epsilon}$ depending smoothly on $t$.
\end{corollary}

\begin{proof}
Since the composition of two conjugate linear automorphisms is $\C$-linear, it is by \ref{2.10} enough   to prove the existence of {\sl one} conjugate linear automorphism of the above form. To this end let $\u$ be a compact real form of $\g$ and $\omega:\g\to\g$ the conjugation with respect to $\u$. Since there exists $\alpha\in\Int\g$ with $\sigma^{-1} \u=\alpha \u,\ \tilde\sigma:=\sigma\alpha$ leaves $\u$ invariant and $L(\g,\sigma)$ and $L(\g,\tilde\sigma)$ are isomorphic. We therefore may assume that $\sigma$ leaves $\u$ invariant. But then $\varphi u(t):=\omega(u(t))$ is a conjugate linear automorphism of the desired form.
\end{proof}

\section {Isomorphisms of smooth affine Kac-Moody algebras}

Let $\g$ be from now on a simple Lie algebra over $\F=\R$ or $\C$ and $\sigma\in\Aut\g$. On $L(\g,\sigma)$ there exists a natural symmetric bilinear form given by $(u,v):=\frac{1}{2\pi}{\int\limits^{2\pi}_0}(u(t),v(t))_0\ dt$ where $(,)_0$ denotes the Killing form on $\g$. It  satisfies $([u,v]_0,w)=(u,[v,w]_0)$ and $(u',v)=-(u,v')$ for all $u,v,w\in L(\g,\sigma)$, where
$'$ denotes differentiation. Let
\[
\hat L(\g,\sigma):=L(\g,\sigma)\oplus \F c\oplus\F d
\]
as vector space with bracket
\[
[u+\alpha c+\beta d,\ v+\gamma c+\delta d]:=[u,v]_0+\beta v' -\delta u'+(u',v)c
\]
for all $u,v\in L(\g,\sigma)$ and $\alpha,\beta,\gamma,\delta\in\F$. Then $\hat L(\g,\sigma)$ is a Lie algebra which we call a (smooth, twisted) affine Kac-Moody algebra. The natural bilinear form on $L(\g,\sigma)$ extends to a natural bilinear form on $\hat L(\g,\sigma)$ by $(u+\alpha c+\beta d,\ v+\gamma c+\delta d):=(u,v)+\alpha\delta+\beta\gamma$ i.e. with $c,d\ \bot\ L(\g,\sigma)$ and $(c,c)=(d,d)=0,\ (c,d)=1$. It is biinvariant in the sense that
\[
([x,y],z)=(x,[y,z])
\]
for all $x,y,z\in\hat L(\g,\sigma)$.

Note that $L(\g,\sigma)$ is only a subspace, not a subalgebra of $\hat L(\g,\sigma)$.

\begin{proposition}\
\begin{enumerate}
\item[(i)] The derived algebra $\hat L'(\g,\sigma)$ of $\hat L(\g,\sigma)$ is equal to $L(\g,\sigma)\oplus\F c$
\item[(ii)] $\F c$ is the center of $\hat L(\g,\sigma)$ and $\hat L'(\g,\sigma)$
\item[(iii)] $L(\g,\sigma)$ is  isomorphic to $\hat L'(\g,\sigma)/\F c$
\item[(iv)] $L(\g,\sigma)$ is equal to its derived algebra.
\end{enumerate}
\end{proposition}

\begin{proof}
(ii) and the implications (iv) $\Rightarrow$ (i) $\Rightarrow$ (iii) are straightforward.

To prove (iv), let $u\in L(\g,\sigma)$. By using the lift of a smooth partition of unity on $S^1$ to $\R$ we may assume that $\supp(u)\cap [0,2\pi]$ is arbitrarily small, in particular that there exist $u_1,\dots,u_n\in L(\g,\sigma)$ by \ref{2.3} which are constant on this set and equal to the elements $x_i,\dots, x_n$ of a basis of $\g$. Expressing $u$ as ${\sum\limits_i}f_iu_i$ with $f_i\in C^\infty_{\per}(\R,\F)$ and the $x_i$ as $\sum\limits_{j,k} a_{ijk}\ [x_j,x_k]$ (which is possible since $\g'=\g)$ gives $u={\sum\limits_{i,j,k}}a_{ijk}\ [f_i u_j,u_k]$.
\end{proof}

The goal of this chapter is to extend the results of the previous one to (linear and conjugate linear) isomorphisms between affine Kac-Moody algebras. In an intermediate step we first consider isomorphisms between the derived algebras.

Let $\tilde\g$ be a second simple Lie algebra over $\F$ and $\tilde\sigma\in\Aut\tilde\g$. A linear or conjugate linear isomorphism $\varphi:L(\g,\sigma)\to L(\tilde\g,\tilde\sigma)$ is by Theorem \ref {2.10} and Corollary \ref{2.17} of the form $\varphi u(t)=\varphi_t(u(\lambda(t)))$ with $\lambda(t+2\pi)=\lambda(t)+\epsilon 2\pi$ for some $\epsilon\in\{\pm 1\}$ and isomorphisms $\varphi_t:\g\to\tilde\g$. Each $\varphi'_t\varphi_t^{-1}:\tilde\g\to\tilde\g$ (where $\varphi'_t={d\over dt}\varphi_t)$ is a derivation (even if $\varphi_t$ is conjugate linear) and hence of the form ad $x(t)$ for a unique $x(t)\in\tilde\g$ depending smoothly on $t$.
Moreover $\varphi_{t+2\pi}=\sigma\varphi_t\sigma^{-\epsilon}$ implies $x(t+2\pi)=\tilde\sigma x(t)$ and hence $x\in L(\tilde\g,\tilde\sigma)$. In order to stress the dependancy on the isomorphism $\varphi$, let $\epsilon_\varphi:=\epsilon, \lambda_\varphi:=\lambda,$ and $x_\varphi:=x$. Due to Theorem \ref{2.10} (iii) $\epsilon_\varphi,\lambda'_\varphi$ and $x_\varphi$ are well defined.

\begin{proposition}
Let $\check\varphi:\hat L'(\g,\sigma)\to\hat L'(\tilde\g,\tilde\sigma)$ be a linear or conjugate linear map. Then $\check\varphi$ is an isomorphism (of Lie algebras) if and only if there exists a linear (resp. conjugate linear) isomorphism $\varphi:L(\g,\sigma)\to L(\tilde\g,\tilde\sigma)$ such that
\[
\begin{array}{lll}
          &\check\varphi c&=\epsilon_\varphi c\\
\mbox{and}&\check\varphi u&=\varphi u+(x_\varphi,\varphi u)c
\end{array}
\]
for all $u\in L(\g,\sigma)$.
\end{proposition}

\begin{proof}
Since isomorphisms map centers to centers we may restrict our attention to those $\check\varphi$ with
\[
\begin{array}{lll}
\check\varphi c & = & \alpha c\\
\check\varphi u & = &\varphi u+\mu(u)c
\end{array}
\]
where $\alpha\in\F,\ \alpha\ne 0,\ \varphi:L(\g,\sigma)\to L(\tilde\g,\tilde\sigma)$ is a linear (resp. conjugate linear) vector space isomorphism and $\mu:L(\g,\sigma)\to\F$ is linear (resp. conjugate linear). Then $\check\varphi$ is an isomorphism if and only if $[\check\varphi u,\check \varphi v]=\check\varphi[u,v]$ for all $u,v\in L(\g,\sigma)$ or equivalently if $[\varphi u,\varphi v]_0+((\varphi u)',\ \varphi v) c=\varphi [u,v]_0+\mu([u,v]_0)c+\varphi((u',v)c)$. This in turn is equivalent to $\varphi$ being an isomorphism and $((\varphi u)',\varphi v)-\mu([u,v,]_0)=\alpha(u',v)$ (resp. = $\alpha\overline{(u',v)})$. Thus $\varphi u(t)=\varphi_t u(\lambda(t))$ and $(\varphi u)'(t)=\varphi'_t\varphi_t^{-1}(\varphi u(t))+\varphi (u')(t)\cdot\lambda'(t)=([x_\varphi,\varphi u]_0+\lambda'\varphi(u'))(t)$. Hence $((\varphi u)',\varphi v)=(x_\varphi,\varphi[u,v]_o)+(\lambda'\varphi(u'),\varphi v)$. Since $(\psi x,\psi y)_0=(x,y)_0$ (resp. $\overline{(x,y)}_0$) if $\psi:\g\to\tilde\g$ is an isomorphism (resp. conjugate linear isomorphism) we have $\lambda'\cdot(\varphi(u'),\varphi v)=\epsilon_\varphi(u',v)$ (resp. = $\epsilon_\varphi\overline{(u',v)})$
due to ${\int\limits^{2\pi}_0}(u'(\lambda(t)), v(\lambda(t)))_0\lambda'(t)dt={\int\limits^{\lambda(2\pi)}_{\lambda(0)}}(u'(x),v(x))_0\,dx=2\pi\epsilon_\varphi(u',v)$.

Therefore $\check\varphi$ is an isomorphism if and only if $\varphi$ is an isomorphism and
$(x_\varphi,\varphi[u,v]_0)-\mu([u,v]_0)=(\alpha-\epsilon)(u',v)$ (resp.  $\mu([u,v]_0) = (\alpha-\epsilon)\overline{(u',v)})$ for all $u,v\in L(\g,\sigma)$.

Choosing $u\in L(\g,\sigma)$ such that $(u',u)_0$ is not identically zero, and $v:=fu$ for some $f\in C^\infty_{\per}(\R,\F)$ with ${\int\limits^{2\pi}_0}f(u',u)_0\ne 0$, we get $\alpha=\epsilon$. Hence 3.1 (iv) yields $\mu(u)=(x_\varphi,\varphi u)$ for all $u\in L(\g,\sigma)$.
\end{proof}

\begin{corollary}\label{3.3}
The mapping $\check\varphi\mapsto\varphi$ which associates to any isomorphism $\check\varphi:\hat L'(\g,\sigma)\to\hat L'(\tilde\g,\tilde\sigma)$ the induced isomorphism between the loop algebras, is a bijection. In particular $\Aut(\hat L'(\g,\sigma))$ and $\Aut(L(\g,\sigma))$ are isomorphic.
\end{corollary}

\begin{theorem}\label{3.4}
Let $\hat\varphi:\hat L(\g,\sigma)\to\hat L(\tilde\g,\tilde\sigma)$ be a linear or conjugate linear map. Then $\hat\varphi$ is an isomorphism (of Lie algebras) if and only if there exist $\gamma\in\F$ and a linear (resp. conjugate linear) isomorphism $\varphi:L(\g,\sigma)\to L(\tilde\g,\tilde\sigma)$ with $\lambda'_\varphi=\epsilon_\varphi$ constant such that
\[
\begin{array}{lll}
\hat\varphi c & = & \epsilon_\varphi c\\
\hat\varphi d & = & \epsilon_\varphi d-\epsilon_\varphi x_\varphi+\gamma c\\
\hat\varphi u & = & \varphi u+(x_\varphi,\varphi u)c
\end{array}
\]
for all $u\in L(\g,\sigma)$.
\end{theorem}

\begin{proof}
Let $\hat\varphi d=\beta d+u_\varphi+\gamma c$. Then $\hat\varphi$ is an isomorphism if and only if $\beta\ne 0,\hat\varphi$ maps $\hat L'(\g,\sigma)$ to $\hat L'(\tilde\g,\tilde\sigma)$, the induced map $\check\varphi$ is an isomorphism and $[\hat\varphi d,\hat\varphi u]=\hat\varphi[d,u]$ for all $u\in L(\g,\sigma)$. By 3.2, $\hat\varphi c=\epsilon_\varphi c$, and $\hat\varphi u=\varphi u+(x_\varphi,\varphi u)c$ for a (unique) isomorphism $\varphi:L(\g,\sigma)\to L(\tilde\g,\tilde\sigma)$. Thus $[\hat\varphi d,\hat\varphi u]=\hat\varphi[d,u]$ is equivalent to $\beta(\varphi u)'+[u_\varphi,\varphi u]_0-(u_\varphi,(\varphi u)')c=\varphi(u')+(x_\varphi,\varphi(u'))c$ and hence to (i) $[\beta x_\varphi+u_\varphi,\varphi u]_0=(1-\beta\lambda'_\varphi)\varphi(u')$ and (ii) $(u_\varphi,(\varphi u)')=-(x_\varphi,\varphi(u'))$ as $(\varphi u)'=[x_\varphi,\varphi u]_0+\lambda'_\varphi\varphi(u')$.

Let $t_0\in \R$. Then there exists $u\in L(\g,\sigma)$ with $\varphi u(t_0)=0$ and $(\varphi u)'(t_0)\ne 0$ and thus also with $\varphi(u')(t_0)\ne 0$ by the last equation. This shows that (i) is equivalent to $u_\varphi=-\beta x_\varphi$ and $\lambda'_\varphi\equiv {1\over\beta}$ and hence to $u_\varphi=-\epsilon _\varphi x_\varphi$ and $\lambda'_\varphi=\epsilon_\varphi$ as $\epsilon_\varphi={1\over2\pi}{\int\limits^{2\pi}_0}\lambda'_\varphi(t)dt$. Therefore (ii) is a consequence of (i) and the theorem follows.
\end{proof}

We will call the $\varphi$ above to be \textit{induced} by $\hat\varphi$. It is equal to the restriction of $\hat\varphi$ to $L(\g,\sigma)$ followed by the projection $u+\alpha c+\beta d\mapsto u$.

\medskip
The theorem shows in particular that any isomorphism $\varphi:L(\g,\sigma)\to L(\tilde\g,\tilde\sigma)$ with $\lambda'_\varphi$ constant, can be extended to an isomorphism $\hat\varphi:\hat L(\g,\sigma)\to\hat L(\tilde\g,\tilde\sigma)$. Hence Corollary \ref{2.13} and Corollary \ref{2.15} extend immediately to the affine Kac-Moody case:

\begin{corollary}\label{3.5}
\begin{enumerate}
\item[(i)] If $\hat L(\g,\sigma)$ and $\hat L(\tilde \g,\tilde\sigma)$ are isomorphic then $\g$ and $\tilde\g$ are isomorphic.
\item[(ii)] $\hat L(\g,\sigma)$ and $\hat L(\g,\tilde\sigma)$ are isomorphic if and only if $\sigma$ and $\tilde\sigma$ are conjugate in $\Aut\g/\Int\g$. \qed
\end{enumerate}
\end{corollary}

\begin{corollary}
For any $\sigma\in\Aut\g$ there exists an automorphism $\tilde\sigma\in \Aut\g$ of finite order with $\hat L(\g,\sigma)\cong\hat L(\g,\tilde\sigma)$.\qed
\end{corollary}

Specializing 3.4 to automorphisms let
\[
\begin{array}{lll}
\Aut(\hat L(\g,\sigma),\hat L') & := \{\hat\varphi\in\Aut\hat L(\g,\sigma) & \mid\ \hat\varphi = id \mbox{ on } \hat L'(\g,\sigma)\}\\
&= \{\hat\varphi\in\Aut\hat L(\g,\sigma) & \mid\ \hat\varphi =id \mbox{ on } L(\g,\sigma)\},\\
\Aut_{(,)}\hat L(\g,\sigma) & := \{\hat\varphi\in\Aut\hat L(\g,\sigma)& \mid\ (\hat\varphi x,\hat \varphi y)=(x,y)\ \forall\ x,y\in\hat L(\g,\sigma)\}, \mbox{ and}\\
\Aut' L(\g,\sigma) & := \{\varphi\in\Aut L(\g,\sigma) & \mid\ \lambda'_\varphi \mbox{ constant }\}\ .
\end{array}
\]
Note that $\lambda'_\varphi$ constant yields $\lambda'_\varphi=\epsilon_\varphi$.

\begin{corollary}\label{3.7}
\begin{enumerate}
\item[(i)] $\Aut(\hat L(\g,\sigma),\ \hat L')\cong\F$.
\item[(ii)] $\Aut(\hat L(\g,\sigma),\ \hat L')$ is contained in the center of $\Aut\hat L(\g,\sigma)$.
\item[(iii)] $\Aut\hat L(\g,\sigma)\cong\Aut_{(,)}\hat L(\g,\sigma)\times\F$.
\item[(iv)] The mapping $\Aut_{(,)}\hat L(\g,\sigma)\to \Aut'L(\g,\sigma)$ that associates to each $\hat\varphi$ the induced mapping $\varphi$ on $L(\g,\sigma)$ is an isomorphism.
\end{enumerate}
\end{corollary}

\begin{proof}
\begin{enumerate}
\item[(i)] The $\hat\varphi\in\Aut(\hat L(\g,\sigma),\hat L')$ are the automorphisms with $\hat\varphi c=c$ and $\hat\varphi u=u$ and thus with $\hat\varphi d=d+\gamma c$ by 3.4. Hence $\hat\varphi\mapsto\gamma$ defines an isomorphism.
\item[(ii)] Let $\hat\varphi\in\Aut\hat L(\g,\sigma)$ with $\hat\varphi c=c,\ \hat\varphi d=d+\gamma c$ and $\hat\varphi u=u$ for all $u\in L(\g,\sigma)$. Then $\hat\varphi$ commutes with all $\hat\psi\in\End\hat L(\g,\sigma)$ which leave $\F c$ and $\hat L'(\g,\sigma)$ invariant and satisfy $(\hat\psi c,d)=(c,\hat\psi d)$. In particular it commutes with all automorphisms.
\item[(iii)] $\hat\varphi\in\Aut\hat L(\g,\sigma)$ with $\hat\varphi c=\epsilon_\varphi c,\ \hat\varphi d=\epsilon_\varphi d-\epsilon_\varphi x_\varphi+\gamma c$, and $\hat\varphi u=\varphi u+(x_\varphi,\varphi u)c$ leaves the bilinear form $(,)$ invariant if and only if $(\hat\varphi d,\hat\varphi d)=0$ or equivalently if $2\gamma=-\epsilon_\varphi(x_\varphi,x_\varphi)$. Thus the claim follows from (i) and (ii).
\item[(iv)] follows from 3.4 and (iii) above.
\end{enumerate}
\vspace{-27pt}
\end{proof}

3.7 (iii) shows that elements of finite order in $\Aut\hat L(\g,\sigma)$ are contained in $\Aut_{(,)}\hat L(\g,\sigma)$ and are conjugate in $\Aut\hat L(\g,\sigma)$ if and only if they are conjugate in $\Aut_{(,)}\hat L(\g,\sigma)$. Thus we have

\begin{proposition}\label{3.8}
There is a natural bijection between conjugacy classes of elements of finite order of $\Aut\hat L(\g,\sigma)$ and $\Aut'L(\g,\sigma)$.  \qed
\end{proposition}

In case $\F=\C$ let $\overline{\Aut}\hat L(\g,\sigma)$ and $\overline{\Aut}L(\g,\sigma)$ be the sets of \textit{conjugate linear} automorphisms of $\hat L(\g,\sigma)$ and $L(\g,\sigma)$, respectively and $\overline{\Aut}'L(\g,\sigma):=\{\varphi\in\overline{\Aut}L(\g,\sigma)\mid\lambda'_\varphi$ constant$\}$.

\begin{proposition}\label{3.9}
There is a natural bijection between conjugacy classes of elements of finite order of $\overline{\Aut}\hat L(\g,\sigma)$ and $\overline{\Aut}'L(\g,\sigma)$ (where conjugation means conjugation with respect to elements of $\Aut\hat L(\g,\sigma)$ and $\Aut'L(\g,\sigma)$, respectively).
\end{proposition}

\begin{proof}
The proof is in complete analogy to that of \ref{3.7} by considering the enlarged groups $\Aut\hat L(\g,\sigma)\cup\overline{\Aut}\hat L(\g,\sigma)$ and $\Aut'L(\g,\sigma)\cup\overline{\Aut}'L(\g,\sigma)$, respectively. The mapping $\hat\varphi\mapsto\varphi$ is surjective and has kernel $\Aut(\hat L(\g,\sigma),\hat L')\cong\C$ which again splits off as a direct factor (the other factor being $\Aut_{(,)}\hat L(\g,\sigma)\cup\overline{\Aut}_{(,)}\hat L(\g,\sigma)$ where $\overline{\Aut}_{(,)}\hat L(\g,\sigma)=\{\varphi\in\overline{\Aut}\hat L(\g,\sigma) \mid (\varphi x,\varphi y)=\overline{(x,y)}\ \forall\ x,y\})$.
\end{proof}

Later (cf.~Corollaries~\ref{4.12} and \ref{5.9}) we will show that $\Aut'L(\g,\sigma)$ can be replaced by $\Aut L(\g,\sigma)$ in 3.8 (and $\overline{\Aut}'L(\g,\sigma)$ by $\overline{\Aut}L(\g,\sigma)$ in 3.9). Here we will prove part of this statement in Proposition~\ref{3.11}.

\begin{lemma}\label{3.10}
Let $\varphi:L(\g,\sigma)\to L(\g,\sigma)$ and $\psi:L(\g,\sigma)\to L(\g,\tilde\sigma)$ be isomorphisms (possibly conjugate linear) of the form $\varphi u(t)=\varphi_t(u(\lambda(t)))$ and $\psi u(t)=\psi_t(u(\mu(t)))$. Then $\psi\varphi\psi^{-1}(u)(t)=\tilde\varphi_t(u(\tilde\lambda(t)))$ where $\tilde\lambda=\mu^{-1}\circ\lambda\circ\mu$ and $\tilde\varphi_t=\psi_t\circ\varphi_{\mu(t)}\circ(\psi_{\tilde\lambda(t)})^{-1}$.
\end{lemma}

\begin{proof}
Observe that $\psi\varphi(u)(t)=\psi_t(\varphi u(\mu(t)))=\psi_t\varphi_{\mu(t)}(u(\lambda\circ\mu(t)))$ and thus in particular $\psi^{-1}u(t)=(\psi_{\mu^{-1}(t)})^{-1}(u(\mu^{-1}(t)))$.
\end{proof}

\begin{proposition}\label{3.11}
Let $\varphi\in\Aut L(\g,\sigma)$ (resp. $\overline{\Aut}L(\g,\sigma)$) be of finite order. Then there exists  $\psi\in\Aut L(\g,\sigma)$ of the first kind such that $\psi\varphi\psi^{-1}\in\Aut'L(\g,\sigma)$ (resp. $\overline{\Aut}'L(\g,\sigma)$). If $\varphi$ is of the second kind the order of $\varphi$ is even.
\end{proposition}

\begin{proof}
Let $\varphi$ be of order $q$ and $\varphi u(t)=\varphi_t(u(\lambda(t)))$ with $\varphi_t$ and $\lambda(t)$ smooth and $\lambda(t+2\pi)=\lambda(t)+\epsilon 2\pi$. Then $\varphi^{k}u(t)=\varphi_t\varphi_{\lambda(t)}\dots\varphi_{\lambda^{k-1}(t)}(u(\lambda^k(t)))$ for all $k\in\N$ where $\lambda^k:=\lambda\circ\dots\circ\lambda$ denotes the $k^{th}$ iterate. Thus $\lambda^q=id+p2\pi$ for some $p\in\Z$. If $\varphi$ is of the second kind $(\epsilon=-1)$ then $q$ is necessarily even as $\epsilon^k\lambda^k$ is an orientation preserving diffeomorphism of $\R$ for all $k\in\N$.

Let $\nu:={1\over q}{\sum\limits^{q-1}_{k=0}}\epsilon^k\lambda^k$. Then also $\nu$ is an orientation preserving diffeomorphism of $\R$ (note that its derivative is positive everywhere)
with $\nu(t+2\pi)=\nu(t)+2\pi$. It moreover satisfies $\nu(\lambda(t))=\epsilon\nu(t)+{\epsilon\over q}\ p2\pi$ and thus $\nu\circ\lambda\circ\nu^{-1}(t)=\epsilon t+{\epsilon\over q}\ p2\pi$. Hence $\psi u(t):= u(\nu^{-1}(t))$ defines by Theorem \ref{2.10}  an automorphism of the first kind on $L(\g,\sigma)$ and $\psi\varphi\psi^{-1}\in\Aut'L(\g,\sigma)$ (resp. $\overline{\Aut}'L(\g,\sigma))$ by 3.10.
\end{proof}

The main idea of the proof of \ref{3.11} is a variation of the proof that diffeomorphisms of $S^1$ of finite order are conjugate to rotations or reflections. Since this is false for arbitrary diffeomorphisms of $S^1$,  also 3.11 would be false for arbitrary $\varphi\in\Aut L(\g,\sigma)$.

\section{Automorphisms of the first kind of finite order}

Let $\g$ be as before a simple Lie algebra over $\F=\R$ or $\C$ and $\sigma,\tilde\sigma\in\Aut\g$. Let $\Aut^q\g:=\{\varphi\in\Aut\g\mid \ord(\varphi)=q\}$.

\begin{definition}
Two automorphisms $\varphi\in\Aut L(\g,\sigma)$ and $\chi\in\Aut L(\g,\tilde\sigma)$ are called {\rm quasiconjugate} if there exists an isomorphism $\psi:L(\g,\sigma)\to L(\g,\tilde\sigma)$ such that $\chi=\psi\varphi\psi^{-1}$. Similarly $\hat\varphi\in\Aut\hat L(\g,\sigma)$ and $\hat\chi\in\Aut \hat L(\g,\tilde\sigma)$ are called quasiconjugate if there exists an isomorphism $\hat\psi$ with $\hat\chi=\hat\psi\hat\varphi\hat\psi^{-1}$.
\end{definition}

We call $\varphi$ and $\chi$ (resp. $\hat\varphi$ and $\hat\chi$)  only conjugate  if they are quasiconjugate and $\tilde\sigma=\sigma$. It will turn out (Corollary~\ref{4.13}) that any automorphism of the first kind of finite order is quasiconjugate (but not necessarily conjugate) to one with $\varphi_t$ constant. From this it seems to be clear that conjugacy classes of these automorphisms can be parametrized by simple invariants. But we will prove this first and deduce \ref{4.13} as a corollary.

\medskip
Let $\Aut_1 L(\g,\sigma):=\{\varphi\in\Aut L(\g,\sigma)\mid \varphi$ of the first kind$\}$ and $\Aut'_1L(\g,\sigma):=\{\varphi\in\Aut_1L(\g,\sigma)\mid\lambda_\varphi'=1\}$. Recall that any $\varphi\in\Aut_1L(\g,\sigma)$ of finite order is conjugate within $\Aut_1L(\g,\sigma)$ to an element of $\Aut_1'L(\g,\sigma)$ (Proposition \ref{3.11}).

\begin{lemma}
Let $\varphi\in\Aut'_1L(\g,\sigma)$ be of order $q$. Then there exist unique $p\in\{0,\dots,q-1\}$ and $\varphi_t\in\Aut\g\ (t\in\R)$ such that
\[
\varphi u(t)=\varphi_t(u(t+{p\over q}2\pi))
\]
for all $u\in L(\g,\sigma)$, and this $\varphi_t$ depends smoothly on $t$.
\end{lemma}

\begin{proof}
By Theorem \ref{2.10} $\varphi u(t)=\varphi_t(u(\lambda(t)))$ and by assumption $\lambda(t)=t+t_0$ for some $t_0\in\R$. $\varphi^q=\id$ implies $\lambda^q=id+p2\pi$ for some $p\in\Z$ and thus $t_0={p\over q}2\pi$. By replacing $\lambda$ eventually by $\lambda+2m\pi$ for some $m\in\Z$ (and $\varphi_t$ by  $\varphi_t\sigma^{-m}$) we may assume $0\le p<q$. This $p$ is then unique and hence $\varphi_t$ as well. Since $\lambda$ is smooth also $\varphi_t$ is smooth.
\end{proof}

We now associate to each $\varphi\in\Aut'_1L(\g,\sigma)$ of order $q$ an ``invariant'' as follows.

\begin{definition}\label{4.3}
For $q\in\N$ and $p\in\{0,1,\dots,q-1\}$ let $r=r(p,q),\ p'=p'(p,q),\ q'=q'(p,q),\ l=l(p,q)$ and $m=m(p,q)$ be the uniquely determined integers with $r>0,\ p=rp',\ q=rq',\ lp'+mq'=1$ and  $0\le l<q'$. (In particular $r=(p,q)$ is the greatest common divisor of $p$ and $q$).
\end{definition}

Let $\varphi\in\Aut'_1L(\g,\sigma)$ of order $q$ with $\varphi u(t)=\varphi_t(u(t+t_0)),\ t_0={p\over q}\ 2\pi$, and $r,
p',q',l,m$ as defined above. We then have $\varphi^{q'}u(t)=P_t(u(t))$ and $\varphi^lu(t)=\Lambda_t(u(t+{2\pi\over q'}))$ where 
\begin{eqnarray*}
P_t&=&\varphi_t\varphi_{t+t_0}\dots\varphi_{t+(q'-1)t_0}\sigma^{p'}\mbox{ and}\\
\Lambda_t&=&\varphi_t\varphi_{t+t_0}\dots\varphi_{{t+(l-1)t}_0}\sigma^{-m}\ .
\end{eqnarray*}
Since $\varphi$ is of order $q$, this implies that $P_t$ is of order $r$.
Hence $P_t=\alpha_t\varrho_0\alpha_t^{-1}$ for some $\varrho_0\in\Aut^r\g$ and some $\alpha_t\in\Aut\g$ depending smoothly on $t$. For $\g$ has only finitely many automorphisms of order $r$ up to conjugation. From $\varphi^l\varphi^{q'}=\varphi^{q'}\varphi^l$ we get $\Lambda_tP_{t+2\pi/q'}=P_t\Lambda_t$ and thus $\alpha^{-1}_{t+2\pi/q'}\Lambda^{-1}_t\alpha_t\in(\Aut\g)^{\varrho_0}$. In the following we fix for each $r\in\N$ a set ${\cal A}^r(\g)\subset\Aut^r\g$ of representatives of conjugacy classes of automorphisms of $\g$ of order $r$. If $\varrho\in{\cal A}^r(\g)$ we consider in $(\Aut\g)^\varrho$ the equivalence relation $\beta\sim\tilde\beta$ if and only if $\beta$ and $\tilde\beta$ are conjugate in $\pi_0((\Aut\g)^\varrho)=(\Aut\g)^\varrho/((\Aut\g)^\varrho)_0$, i.e. if $\tilde\beta=\gamma\beta\gamma^{-1}\delta$ for some $\gamma\in(\Aut\g)^\varrho$ and $\delta\in((\Aut\g)^\varrho)_0$. We denote by $[\beta]$ the equivalence class containing $\beta$. There are only finitely many equivalence classes as $\pi_0 ((\Aut\g)^\varrho)$ is finite (cf.~Appendix A).

\begin{definition}\ 
\begin{enumerate}
\item[(i)] For $q\in\N$ let $\J^q_1(\g):=\{(p,\varrho,[\beta])\mid p\in\{0,\dots,q-1\},\varrho\in{\cal A}^{(p,q)}(\g),\beta\in(\Aut\g)^\varrho\}$
\item[(ii)] If $\varphi\in\Aut'_1L(\g,\sigma)$ is of order $q$ and $p,\varrho_0,\alpha_t,q'$ and $\Lambda_t$ are chosen as above then
\[
(p,\varrho_0,[\alpha^{-1}_{t+2\pi/q'}\Lambda^{-1}_t\alpha_t])\in \J_1^q(\g)
\]
is called the invariant of $\varphi$.
\end{enumerate}
\end{definition}

\begin{remark}\label{4.5}\rm
$\J_1^q(\g)$ is a finite set.
\end{remark}

{\bf Example}: The invariant of $\varphi\in\Aut'_1L(\g,\sigma)$ with $\varphi_t\equiv\varphi_0$ (note that this requires $\varphi_0\sigma=\sigma\varphi_0)$ is $(p,\alpha\varphi_0^{q'}\sigma^{p'}\alpha^{-1},[\alpha\varphi_0^{-l}\sigma^m\alpha^{-1}])$ where $\alpha$ is chosen  such  that $\alpha\varphi_0^{q'}\sigma^{p'}\alpha^{-1}\in{\cal A}^r(\g)$. If moreover $\varphi$ is an involution $(q=2)$ and $\varphi_0\in{\cal A}^r(\g)$ then its invariant is $(0,\varphi_0,[\sigma])$ or $(1,\id,[\varphi_0^{-1}])$.

\begin{proposition}\label{4.6}
Each element of $\J_1^q(\g)$ occurs as the invariant of some $\varphi\in\Aut'_1L(\g,\sigma)$ of order $q$ for some $\sigma$.
\end{proposition}

\begin{proof}
Let $(p,\varrho_0,[\beta])\in \J_1^q(\g)$. Then $\varrho_0$ and $\beta$ commute and hence also $\varphi_0:=\varrho_0^m\beta^{-p'}$ and $\sigma:=\varrho_0^l\beta^{q'}$ where $p', q', l$ and $m$ are determined by 4.3. Therefore $\varphi u(t):=\varphi_0 u(t+{p\over q}2\pi)$ defines an element of $\Aut_1' L(\g,\sigma)$ of order $q$. By the example above it has invariant $(p,\varrho_0,[\beta])$.\\
\end{proof}

\begin{proposition}\label{4.7}
\begin{enumerate}
\item[(i)] The invariant of each $\varphi\in\Aut'_1L(\g,\sigma)$ of finite order is well defined.
\item[(ii)] Automorphisms $\varphi\in\Aut'_1L(\g,\sigma)$ and $\tilde\varphi\in \Aut'_1L(\g,\tilde\sigma)$ that are of finite order and quasiconjugate by an isomorphism $\psi:L(\g,\sigma)\to L(\g,\tilde\sigma)$ of the first kind have equal invariants.
\end{enumerate}
\end{proposition}

\begin{proof}\ 
\begin{enumerate}
\item[(i)] Only the $\alpha_t$ in the definition of the invariant of $\varphi$ is not determined completely and could be replaced by $\alpha_t\beta_t$ where $\beta_t\in(\Aut\g)^{\varrho_0}$. But this does not affect $[\alpha^{-1}_{t+2\pi/q'}\Lambda^{-1}_t\alpha_t]$.
\item[(ii)] Let $\varphi$ and $\tilde\varphi$ be of order $q$ and of the form $\varphi u(t)=\varphi_t(u(t+t_0))$ and $\tilde\varphi u(t)=\tilde\varphi_t(u(t+\tilde t_0))$ where $0\le p<q$ and $t_0={p\over q} 2\pi, \tilde t_0={\tilde p\over q}2\pi$. Let $\psi u(t)=\psi_t(u(\mu(t)))$ with $\mu(t+2\pi)=\mu (t)+2\pi$. Then $\tilde\varphi_t=\psi_t\varphi_{\mu(t)}
\psi^{-1}_{t+\tilde t_0}$ and $t+\tilde t_0=\mu^{-1}(\mu(t)+t_0)$ by Lemma \ref{3.10}. Therefore $\mu(t+2\pi)=\mu(t)+2\pi$ implies $\tilde p=p$ and $\mu(t+t_0)=\mu(t)+t_0$. Let $r,p', q', l$ and $m$ be defined as in \ref{4.3}. Then
\[
\mu(t+{2\pi/ q'})=\mu(t)+{2\pi/ q'}
\]
as ${1\over q'}=l{p\over q}+m$. Let $P_t$, $\Lambda_t$, $\alpha_t$, $\varrho_0$, $\tilde P_t$, $\tilde\Lambda_t$, and
$\tilde\varrho_0$ be as in the definition of the invariants of $\varphi$ and $\tilde\varphi$ respectively. Then $\tilde P_t=\psi_t\varphi_{\mu(t)}\varphi_{\mu(t)+t_0}\cdots\varphi_{\mu(t)+(q'-1)t_0}\psi^{-1}_{t+p'2\pi}\tilde\sigma^{p'}=\psi_t P_{\mu(t)}\psi_t^{-1}=\psi_t\alpha_{\mu(t)}\varrho_0\alpha_{\mu(t)}^{-1}\psi_t^{-1}$, which implies $\psi^{-1}_{t+p'2\pi}\tilde\sigma^{p'}=\sigma^{p'}\psi_t^{-1}$. Similarly we obtain $\tilde\Lambda_t=\psi_t\Lambda_{\mu(t)}\psi^{-1}_{t+2\pi/q'}$. Therefore we can choose $\tilde\alpha_t:=\psi_t\alpha_{\mu(t)}$ which yields $\tilde\varrho_0=\varrho_0$ and $\tilde\alpha^{-1}_{t+2\pi/q'}\ \tilde\Lambda_t^{-1}\tilde\alpha_t=\alpha^{-1}_{\mu(t)+2\pi/q'}\ \Lambda^{-1}_{\mu(t)}\alpha_{\mu(t)}$. Hence the invariants of $\varphi$ and $\tilde\varphi$ coincide.
\end{enumerate}
\end{proof}

If $\varphi\in\Aut_1L(\g,\sigma)$ has finite order we choose $\psi\in\Aut_1L(\g,\sigma)$ with $\psi\varphi\psi^{-1}\in\Aut_1'L(\g,\sigma)$ (cf. Proposition \ref{3.11}) and define the \textit{invariant }of $\varphi$ to be that of $\psi\varphi\psi^{-1}$. By Proposition~\ref{4.7} (ii) this is well defined and invariant under quasiconjugation with isomorphisms of the first kind.

\medskip
We next describe how the invariant behaves under quasiconjugation with isomorphisms of the second kind.

Let $\imath_q:\J_1^q(\g)\to \J_1^q(\g)$ be the involution with $\imath_q((0,\varrho,[\beta]))=(0,\varrho,[\beta^{-1}])$  and $\imath_q((p,\varrho,[\beta]))=(q-p,\varrho,[\beta^{-1}\varrho])$ for $p\in\{1,\dots,q-1\}$.

\begin{definition}\label{4.8}
We call two elements $a,b\in J^q_1(\g)$ {\rm opposite} if $b=\imath_q(a)$.
\end{definition}

\begin{proposition}\label{4.9}
Let $\varphi\in\Aut_1L(\g,\sigma)$ be of order $q$ and $\psi:L(\g,\sigma)\to L(\g,\tilde\sigma)$ be an isomorphism of the second kind. Then the invariant of $\psi\varphi\psi^{-1}$ is opposite to that of $\varphi$.
\end{proposition}

\begin{proof}
It is enough to consider the special case $\psi u(t)=u(-t)$ and $\tilde\sigma=\sigma^{-1}$ since any other isomorphism of the second kind is a composition of this with an isomorphism of the first kind (which does not change the invariant). Now a direct calculation gives the result.
\end{proof}

\begin{remark}\label{4.10}\rm 
If $q=2$ then $\imath_2((p,\varrho,[\beta]))=(p,\varrho,[\beta^{-1}])$ for all $(p,\varrho,[\beta])\in J^2_1(\g)$ (note $\varrho=\id$ if $p=1$). Hence $\imath_2$ is the identity if for each $\varrho\in\Aut\g$ with $\varrho^2=\id$ each element of $\pi_0((\Aut\g)^{\varrho})$ is conjugate to its inverse. This is for example the case if $\g$ is compact or complex as will be explained in Chapter 6. Hence the invariant of an involution of the first kind does not change in this case under quasiconjugations with arbitrary isomorphisms.
\end{remark}

\begin{theorem}\label{4.11}
Let $\varphi\in\Aut_1L(\g,\sigma)$ and $\tilde\varphi\in\Aut_1L(\g,\tilde\sigma)$ be two automorphisms of order~$q$ with the same (resp. opposite) invariants. Then they are quasiconjugate by an automorphism $\psi$ of the first (resp.~second) kind. Moreover, if $\varphi\in\Aut'_1L(\g,\sigma)$ and $\tilde\varphi\in\Aut'_1L(\g,\tilde\sigma)$ then $\psi$ can be chosen to be of the form $\psi u(t)=\psi_tu(\pm t)$.
\end{theorem}

\begin{proof}
We may assume the invariants of $\varphi$ and $\tilde\varphi$ to coincide since otherwise we could first conjugate $\varphi$ by $\psi_1:L(\g,\sigma)\to L(\g,\sigma^{-1})$ with $\psi_1u(t)=u(-t)$ and apply~\ref{4.9}.
After a first conjugation we may assume $\varphi u(t)=\varphi_tu(t+{p\over q}2\pi),\ \tilde\varphi u(t)=\tilde\varphi_t u(t+{p\over q}2\pi)$ with $0\le p<q$. To find $\psi$ with $\psi u(t)=\psi_t u(t)$ it thus suffices to find a smooth curve $\psi_t$ in $\Aut\g$ with
\begin{enumerate}
\item[(1)] $\psi_{t+2\pi}=\tilde\sigma\psi_t\sigma^{-1}$ and
\item[(2)] $\psi_t\varphi_t\psi^{-1}_{t+t_0}=\tilde\varphi_t$,
\end{enumerate}
where $t_0={p\over q} 2\pi$ as before. Using our standard notations $r,p',q',l,m$ (with $p=rp',\ q=rq',\ lp'+mq'=1,\ 0\le l<q'),\ P_t, \alpha_t,\Lambda_t$ (with $\varphi^{q'}u(t)=P_tu(t),\ P_t=\alpha_t\varrho_0\alpha_t^{-1},\ \varphi^l u(t)=\Lambda_tu(t+2\pi/q'))$ and $\tilde P_t,\tilde\alpha_t,\tilde \Lambda_t$ correspondingly, (1) and (2) are equivalent to
\begin{enumerate}
\item[(1')] $\psi_{t+2\pi/q'}=\tilde\Lambda_t^{-1}\psi_t\Lambda_t$ and
\item[(2')] $\psi_t=\tilde P_t^{-1}\psi_t\ P_t$.
\end{enumerate}
In fact, if we extend the mappings $\psi$, $\varphi$, $P$, $\Lambda$ to \textit{all} smooth $u:\R\to\g$ by the same formulas
and let $\tau u(t):=\sigma u(t-2\pi)$ and define $\tilde\varphi$, $\tilde P$, $\tilde\Lambda$, $\tilde\tau$ correspondingly then $\varphi\tau=\tau\varphi$, $P=\varphi^{q'}\tau^{p'}$, $\Lambda=\varphi^l\tau^{-m}$, $\varphi=\Lambda^{p'}P^m$, $\tau=\Lambda^{-q'}P^l$. Hence (1) and (2) are equivalent to
\begin{enumerate}
\item[(*)] $\tilde\tau\psi=\psi\tau$ and $\tilde\varphi\psi=\psi\varphi$
\end{enumerate}
while (1') and (2') are equivalent to
\[
\tilde\Lambda\psi=\psi\Lambda \mbox{ and } \tilde P\psi=\psi P\ .
\]
and thus to (*).

With $P_t=\alpha_t\varrho_0\alpha^{-1}_t,\ \tilde P_t=\tilde\alpha_t\varrho_0\tilde\alpha_t^{-1}$ and $\chi_t:=\tilde\alpha_t^{-1}\psi_t\alpha_t$, (2') is equivalent to $\chi_t\in(\Aut\g)^{\varrho_0}$ and (1') to
\[
\chi_{t+2\pi/q'}=\tilde\beta_t^{-1}\chi_t\beta_t
\]
where $\beta_t:=\alpha_t^{-1}\Lambda_t\alpha_{t+2\pi/q'}$ and $\tilde\beta_t:=\tilde\alpha_t^{-1}\tilde\Lambda_t\tilde\alpha_{t+2\pi/q'}$ are curves in $(\Aut\g)^{\varrho_0}$. But this equation for $\chi_t\in(\Aut\g)^{\varrho_0}$ can be solved since by assumption $[\tilde\beta_t]=[\beta_t]$, i.e. $\tilde\beta_t=\gamma\beta_t\gamma^{-1}\delta_t$ for some $\gamma\in(\Aut\g)^{\varrho_0}$ and some smooth curve $\delta_t$ in $((\Aut\g)^{\varrho_0})_0$. In fact, for small $t$ we may choose  $\chi_t:\equiv\gamma$ and $\chi_{t+2\pi/q'}:=\tilde\beta_t^{-1}\gamma\beta_t=\delta^{-1}_t\gamma$ and connect these pieces smoothly in $[0,2\pi/q']$. The periodic extension from $[0,2\pi/q']$ to all of $\R$ by $\chi_{t+2\pi/q'}=\tilde\beta_t^{-1}\chi_t\beta_t$ gives then the desired solution.
\end{proof}

For any $\sigma\in\Aut\g$, let
\[
J_1^q(\g,\sigma):=\{(p,\varrho,[\beta]\in J_1^q(\g)\mid\varrho^l\beta^{q'}
\mbox{ is conjugate to }\sigma \mbox{ in }\Aut\g/\Int\g\}
\]
where $l=l(p,q)$ and $q'=q'(p,q)$ are as in~\ref{4.3}. Then $J_1^q(\g,\sigma)$ contains precisely the invariants of those automorphisms of the first kind of order $q$ which are defined on algebras $\hat L(\g,\tilde\sigma)$ and $L(\g,\tilde\sigma)$ isomorphic to $\hat L(\g,\sigma)$ resp. $L(\g,\sigma)$. In fact, if two automorphisms have the same invariant then they are quasiconjugate, and in particular the algebras on which they are defined are isomorphic. Thus the invariant determines the isomorphism type of the algebra. On the other hand, the invariant $(p,\varrho,[\beta])$ can be realized on $\hat L(\g,\varrho^l\beta^{q'}$ (resp. $L(\g,\varrho^l\beta^{\varrho'}))$ by~\ref{4.6}.

We denote by $\Aut^q_1L(\g,\sigma)$ and $\Aut^q_1\hat L(\g,\sigma)$ the sets of automorphism of $L(\g,\sigma)$, resp., $\hat L(\g,\sigma)$ of the first kind of order $q$. The invariant of a $\hat\varphi\in\Aut^q_1\hat L(\g,\sigma)$ is by definition the invariant of the induced $\varphi\in\Aut^q_1L(\g,\sigma)$. Then Proposition~\ref{4.6}, \ref{4.7}, \ref{4.9} and Theorem~\ref{4.11} hold correspondingly for automorphisms of affine Kac-Moody algebras. 
By combining Proposition~\ref{4.6} and Theorem~\ref{4.11} we thus have:

\begin{corollary}\label{4.12}
The mapping that associates to each automorphism its invariant induces the following bijections
\[
\Aut^q_1\hat L(\g,\sigma)/\Aut_1\hat L(\g,\sigma)\to J_1^q(\g,\sigma)
\]
and
\[
\Aut^q_1 L(\g,\sigma)/\Aut_1L(\g,\sigma)\to J_1^q(\g,\sigma)\ ,
\]
where the quotients denote conjugacy classes. Moreover in case $q=2$\quad $\Aut_1\hat L(\g,\sigma)$ (resp. $\Aut_1L(\g,\sigma)$) can be replaced by $\Aut\hat L(\g,\sigma)$ (resp. $\Aut L(\g,\sigma)$). \qed
\end{corollary}

\begin{corollary}\label{4.13}
Any element of $\Aut^q_1\hat L(\g,\sigma)$ (resp. $\Aut^q_1 L(\g,\sigma)$) is quasiconjugate to some $\hat\varphi\in\Aut_1^q \hat L(\g,\tilde\sigma)$ with $\hat\varphi c=c$, $\hat\varphi d=d$ and $\hat\varphi u$ being a twisted loop of the form $\hat\varphi u(t)=\varphi_0(u(t+{p\over q}2\pi))$ for all $u\in L(\g,\tilde\sigma)$ where $\varphi_0\in\Aut \g$ is constant and $p\in\{0,1,\dots,q-1\}$ (resp. to some $\varphi$ with $\varphi u(t)=\varphi_0(u(t+{p\over q}2\pi)$).
\end{corollary}

\begin{remark}\label{4.14}\rm 
In general it is not possible to {\textrm conjugate} an element of $\Aut^q_1 L(\g,\sigma)$ to some $\varphi\in\Aut^q_1 L(\g,\sigma)$ with $\varphi u(t)=\varphi_0(u(t+{p\over q}2\pi))$. For example let $\sigma=\id$ and $(0,\varrho,[\beta]\in J_1^q(\g,\id)$ (i.e. $p\in{\cal A}^q$ and $\beta\in \Int\g$) with $\beta\not\in((\Aut\g)^\varrho)_0$ (such $(\varrho,\beta$ exist, see e.g. Chapter 6). Since any $\varphi\in\Aut^q_1 L(\g,\id)$ with $\varphi u(t)=\varphi_0 u(t)$ has invariant $(0,\alpha\varphi_0\alpha^{-1},[\id])$ (cf. the example after~\ref{4.5}) an automorphism of $L(\g,\id)$ with the above invariant $0,\varrho,[\beta])$ can not be conjugate to such a $\varphi$.
\end{remark}

\section{Automorphisms of the second kind of finite order}

We follow the same strategy as in the last chapter and define also for automorphisms of the \textit{second} kind of finite order an invariant, prove that it parametrizes (quasi)conjugacy classes, and derive from this a series of consequences as in Chapter 4.

Let $\g$ be a simple Lie algebra over $\F=\R$ or $\C$ and $\sigma\in\Aut\g$.
Let $\Aut_2\hat L(\g,\sigma)$, $\Aut_2L(\g,\sigma)$, $\Aut_2^q\hat L(\g,\sigma)$ and $\Aut^q_2(\g,\sigma)$ be the sets of automorphisms of the second kind on $\hat L(\g,\sigma)$ and $L(\g,\sigma)$ of arbitrary order and of order $q$, respectively.

\begin{lemma}\label{5.1}
\begin{enumerate}
\item[(i)] Let $\varphi\in\Aut^q_2L(\g,\sigma)$. Then $q$ is even and there exists $\psi\in\Aut_1L(\g,\sigma)$ with $\psi\varphi\psi^{-1}(u)(t)=\varphi_t(u(-t))$ for some smooth curve $\varphi_t$ in $\Aut\g$.
\item[(ii)] $\varphi\in\Aut_2L(\g,\sigma)$ and $\tilde\varphi\in\Aut_2L(\g,\tilde\sigma)$ are quasiconjugate by an isomorphism of the first kind if and only if they are quasiconjugate by an isomorphism of the second kind. In particular $\Aut_2^qL(\g,\sigma)/\Aut_1L(\g,\sigma)=\Aut_2^qL(\g,\sigma)/\Aut L(\g,\sigma)$. The same is true if $L(\g,\sigma)$ is replaced by $\hat L(\g,\sigma)$.
\end{enumerate}
\end{lemma}

\begin{proof}
\begin{enumerate}
\item[(i)] By 3.11, $q$ is even and there exists $\psi_1\in\Aut_1L(\g,\sigma)$ with $\psi_1\varphi\psi_1^{-1}u(t)=\varphi_t(u(-t+t_0))$ for some smooth curve $\varphi_t\in\Aut\g$ and $t_0\in\R$. A further conjugation by $\psi_2$ with $(\psi_1u)(t)=u(t+t_0/2)$ yields therefore the result.
\item[(ii)] Conjugation of $\varphi$ by an isomorphism $\psi:L(\g,\sigma)\to L(\g,\tilde\sigma)$ or by $\psi\varphi$ is the same. But $\psi\varphi$ is of the second kind if and only if $\varphi$ is of the first kind. The same is true for $\hat\varphi$ and $\hat\psi$.
\end{enumerate}
\vspace{-27pt}
\end{proof}

Let $\varphi\in\Aut^{2q}_2L(\g,\sigma)$ with $\varphi u(t)=\varphi_t(u(-t))$. Then $u(t)=\varphi^{2q}(u)(t)=(\varphi_t\varphi_{-t})^{q}(u(t))$ whence $(\varphi_t\varphi_{-t})^{q}=\id$. Since $\g$ has up to conjugation only finitely many automorphisms of order $\le q$ the order of $\varphi_t\varphi_{-t}$ is constant and thus equal to $q$. Hence there exists
a smooth curve $\alpha_t$ in $\Aut\g$ and $\varrho_0\in\Aut^q\g$ with
\[
\varphi_t\varphi_{-t}=\alpha_t\ \varrho_0\ \alpha_t^{-1}\ .
\]
Let $\varphi_+:=\alpha^{-1}_0\varphi_0\alpha_0$ and $\varphi_-:=\alpha^{-1}_\pi\varphi_\pi\sigma^{-1}\alpha_\pi$. Then $\varphi^2_+=\varphi^2_-=\varrho_0$
as $\varphi_{-\pi}=\sigma^{-1}\varphi_\pi\sigma^{-1}$ which in turn follows from
$\varphi_{t+2\pi}=\sigma\varphi_t\sigma$ (cf. Theorem \ref{2.10}).

\begin{definition}\label{5.2}
Two pairs $(\varphi_+,\varphi_-),\ (\tilde\varphi_+,\tilde\varphi_-)\in(\Aut\g)^2$  with $\varphi_+^2=\varphi_-^2$ and $\tilde\varphi_+^2=\tilde\varphi_-^2$
are called equivalent if there exist $\alpha,\beta\in\Aut\g$ with $\alpha^{-1}\beta\in((\Aut\g)^{\varphi_+^2})_0$ such that $\tilde\varphi_+=\alpha\varphi_+\alpha^{-1}$ and $\tilde\varphi_-=\beta\varphi_-\beta^{-1}$ or $\tilde\varphi_+=\alpha\varphi_-\alpha^{-1}$ and $\tilde\varphi_-=\beta\varphi_+\beta^{-1}$.
This defines an equivalence relation and 
we denote the equivalence class of $(\varphi_+,\varphi_-)$ by $[\varphi_+,\varphi_-]$.
\end{definition}

In particular $(\varphi_+,\varphi_-)\sim(\varphi_-,\varphi_+)\sim(\alpha\varphi_+\alpha^{-1},\alpha\varphi_-\alpha^{-1})\sim(\varphi_+,\beta\varphi_-\beta^{-1})$ for all $\alpha\in\Aut\g$ and $\beta\in((\Aut\g)^{\varphi_+^2})_0$.

\begin{definition}
\begin{enumerate}
\item[(i)] For any $q\in\N$ let
\[
\J_2^{2q}(\g):=\{[\varphi_+,\varphi_-]\mid\varphi_\pm\in\Aut\g,\varphi_+^2=\varphi_-^2,\ord(\varphi_\pm^2)=q\}\ .
\]
(Note that $\ord\varphi^2=q$ is equivalent to $\ord\varphi=2q$ if $q$ is even and to $\ord\varphi\in\{q,2q\}$ if $q$ is odd).
\item[(ii)] If $\varphi\in\Aut^{2q}_2L(\g,\sigma)$ with $\varphi u(t)=\varphi_t(u(-t))$
 and $\varphi_\pm$ are defined as above we call $[\varphi_+,\varphi_-]\in\J_2^{2q}(\g)$ the invariant of $\varphi$.
\end{enumerate}
\end{definition}

\begin{remark}\rm
\begin{enumerate}
\item[(i)] The invariant of $\varphi$ is well defined. For $\varphi_t\varphi_{-t}=\alpha_t\varrho_0\alpha_t^{-1}=\tilde\alpha_t\tilde\varrho_0\tilde\alpha_t^{-1}$ implies $\beta_t\beta_0^{-1}\in((\Aut\g)^{\varrho_0})_0$ where $\beta_t=\alpha_t^{-1}\tilde\alpha_t$. Hence $[\tilde\varphi_+,\tilde\varphi_-]=[\beta_0^{-1}\varphi_+\beta_0, \beta_\pi^{-1}\varphi_-\beta_\pi]=[\beta_\pi\beta_0^{-1}\varphi_+\beta_0\beta_\pi^{-1},\varphi_-]=[\varphi_+,\varphi_-]$ since $\varphi_+^2=\varrho_0$.
\item[(ii)] $\J_2^{2q}(\g)$ is a finite set.
\end{enumerate}
\end{remark}

\begin{proposition}\label{5.5}
Any $[\varphi_+,\varphi_-]\in\J_2^{2q}(\g)$ is the invariant of some $\varphi\in\Aut^{2q}_2(L(\g,\sigma))$ with $\varphi u(t)=\varphi_t(u(-t))$ for some $\sigma\in\Aut\g$. In fact, one may take $\sigma:=\varphi_-^{-1}\varphi_+$ and $\varphi_t\equiv\varphi_+$.
\end{proposition}

\begin{proof}
Let $\varphi_t:\equiv\varphi_+$ and $\sigma:=\varphi_-^{-1}\varphi_+$.
Then $\varphi_{t+2\pi}=\sigma\varphi_t\sigma$ for all $t$ and hence $\varphi$ with $\varphi u(t):=\varphi_t(u(-t))$ is contained in $\Aut^{2q}_2L(\g,\sigma)$.
Since $\varphi_t\varphi_{-t}=\varphi^2_+$ we may take $\alpha_t=\id$. Thus $\varphi$ has invariant $[\alpha_0^{-1}\varphi_0\alpha_0,\alpha_\pi^{-1}\varphi_\pi\sigma^{-1}\alpha_\pi]=
[\varphi_+,\varphi_+(\varphi^{-1}_-\varphi_+)^{-1}]=[\varphi_+,\varphi_-]$.
\end{proof}

\begin{proposition}\label{5.6}
The invariant of $\varphi\in\Aut^{2q}_2L(\g,\sigma)$ with $\varphi u(t)=\varphi_t(u(-t))$ does not change under quasiconjugation. More precisely, if $\psi:L(\g,\sigma)\to L(\g,\tilde\sigma)$ is an isomorphism such that $\tilde\varphi:=\psi\varphi\psi^{-1}$ is also of the form $\tilde\varphi u(t)=\tilde\varphi_t(u(-t))$
 then $\varphi$ and $\tilde\varphi$ have the same invariant.
\end{proposition}

\begin{proof}
Let $\varphi u(t)=\varphi_t(u(-t))$ and $\varphi_t\varphi_{-t}=\alpha_t\varrho_0\alpha_t^{-1}$.
By 5.1 (ii) we may assume $\psi$ to be of the first kind and hence by \ref{2.10} of the form $\psi u(t)=\psi_t u(\mu(t))$ where $\psi_t$ and $\mu(t)$ are smooth,  $\psi_{t+2\pi}=\tilde\sigma\psi_t\sigma^{-1}$, and $\mu(t+2\pi)=\mu(t)+2\pi$ for all $t$. Lemma 3.10 implies $\tilde\varphi_t=\psi_t\varphi_{\mu(t)}\psi_{-t+2k\pi}^{-1}\tilde\sigma^k$ and $-t=\mu^{-1}(-\mu(t))-2k\pi$ for some $k\in\Z$. Thus $\mu(-t)=-\mu(t)-2k\pi$ and in particular $\mu(0)=-k\pi,\ \mu(\pi)=-(k-1)\pi$. Moreover $\tilde\varphi_t\tilde\varphi_{-t}=\psi_t\varphi_{\mu(t)}\varphi_{-\mu(t)}\psi_t^{-1}=\tilde\alpha_t\varrho_0\tilde\alpha_t^{-1}$ with $\tilde\alpha_t:=\psi_t\alpha_{\mu(t)}$. This yields $\tilde\varphi_+:=\tilde\alpha^{-1}_0\tilde\varphi_0\tilde\alpha_0=\alpha^{-1}_{-k\pi}\varphi_{-k\pi}\sigma^k\alpha_{-k\pi}$ and $\tilde\varphi_-:=\tilde\alpha^{-1}_\pi\tilde\varphi_\pi\tilde\sigma^{-1}\tilde\alpha_\pi=\alpha^{-1}_{-(k-1)\pi}\varphi_{-(k-1)\pi}\sigma^{k-1}\alpha_{-(k-1)\pi}$. From $\varphi_{t+2\pi}\varphi_{-t-2\pi}=\sigma\varphi_t\varphi_{-t}\sigma^{-1}$ we get $\alpha^{-1}_{t+2\pi}\sigma\alpha_t\in(\Aut\g)^{\varrho_0}$ and by conjugating $\tilde\varphi_+$ by $\alpha^{-1}_{-k\pi+2\pi}\sigma\alpha_{-k\pi}$ and $\tilde\varphi_-$ by $\alpha^{-1}_{-(k-1)\pi+2\pi}\sigma\alpha_{-(k-1)\pi}$ we see that $[\tilde\varphi_+,\tilde\varphi_-]$ depends only on $k$ modulo $2$.
Thus it suffices to consider the cases $k=0$ and $k=1$. If $k=0$ then $\tilde\varphi_+=\varphi_+,\ \tilde\varphi_-=\varphi_-$ and the invariants of $\varphi$ and $\tilde\varphi$ coincide. If $k=1$ then $(\tilde\varphi_+,\tilde\varphi_-)=(\alpha^{-1}_{-\pi}\varphi_{-\pi}\sigma\alpha_{-\pi},\alpha^{-1}_0\varphi_0\alpha_0)=(\alpha^{-1}_{-\pi}\sigma^{-1}\varphi_\pi\alpha_{-\pi},\varphi_+)=(\phi^{-1}_\pi\varphi_-\phi_\pi,\phi^{-1}_0\varphi_+\phi_0)$ where $\phi_t:=\alpha^{-1}_t\varphi_t\alpha_{-t}$. Note that $\varphi_+^2=\varrho_0$ and  $\phi_t\in(\Aut\g)^{\varrho_0}$ as $\phi_t\phi_{-t}=\phi_{-t}\phi_t=\varrho_0$. Hence $[\tilde\varphi_+,\tilde\varphi_-]=[\varphi_+,\varphi_-]$.
\end{proof}

If $\varphi\in\Aut_2^{2q} L(\g,\sigma)$ is arbitrary we may choose $\psi\in\Aut L(\g,\sigma)$ by \ref{5.1} with $\psi\varphi\psi^{-1}(u)(t)=\tilde\varphi_t u(-t)$ and define the invariant of $\varphi$ to be that of $\psi\varphi\psi^{-1}$. By the last proposition this is well defined and invariant under quasiconjugation. The invariant of $\hat\varphi\in\Aut^{2q}_2\hat L(\g,\sigma)$ is that of the induced $\varphi$.

\begin{theorem}\label{5.7}
Two automorphisms $\varphi\in\Aut^{2q}_2L(\g,\sigma)$ and
$\tilde\varphi\in\Aut_2^{2q}L(\g,\tilde\sigma)$ are quasiconjugate if
and only if they have the same invariant. Moreover, if in addition $\varphi\in\Aut'L(\g,\sigma)$ and
$\tilde\varphi\in\Aut'L(\g,\tilde\sigma)$ the conjugating isomorphism
$\psi:L(\g,\sigma)\to L(\g,\tilde\sigma)$ can be chosen to be of the form $\psi u(t)=\psi_t(u(t))$.
\end{theorem}

\begin{proof}
Let $\varphi$ and $\tilde\varphi$ have the same invariants, the other direction being clear by the above remarks.
After a first conjugation (see Lemma~\ref{5.1}) we may assume $\varphi u(t)=\varphi_t u(-t)$ and $\tilde\varphi u(t)=\tilde\varphi_t u(-t)$ with $\varphi_{t+2\pi}=\sigma\varphi_t\sigma$, $\tilde\varphi_{t+2\pi}=\tilde\sigma\varphi_t\tilde\sigma$ and $\varphi_t\varphi_{-t}=\alpha_t\varrho_0\alpha^{-1}_t$, $\tilde\varphi_t\tilde\varphi_{-t}=\tilde\alpha_t\tilde\varrho_0\tilde\alpha_t^{-1}$ for some $\varrho_0,\tilde\varrho_0\in\Aut^q(\g)$.
Let $\varphi_+:=\alpha_0^{-1}\varphi_0\alpha_0,\  \varphi_-:=\alpha^{-1}_\pi\varphi_\pi\sigma^{-1}\alpha_\pi,\ \tilde\varphi_+:=\tilde\alpha^{-1}_0\tilde\varphi_0\tilde\alpha_0$, and $\tilde\varphi_-:=\tilde\alpha^{-1}_\pi\tilde\varphi_\pi\tilde\sigma^{-1}\tilde\alpha_\pi$. Then $\varphi_\pm^2=\varrho_0$ and $\tilde\varphi_\pm^2=\tilde\varrho_0$.
Equality of the invariants is equivalent to $\tilde\varphi_+=\alpha\varphi_+\alpha^{-1}$ and $\tilde\varphi_-=\beta\varphi_-\beta^{-1}$ or $\tilde\varphi_+=\alpha\varphi_-\alpha^{-1}$ and $\tilde\varphi_-=\beta\varphi_+\beta^{-1}$ for some $\alpha,\beta\in\Aut\g$ with $\alpha^{-1}\beta\in((\Aut\g)^{\varrho_0})_0$. We may assume  the first possibility holds. For otherwise we could first conjugate $\varphi$ by $\psi_1 u(t):=u(t-\pi)$ which would replace $\varphi_t$ by $\varphi_{t-\pi}\sigma$ and $(\varphi_+,\varphi_-)$ by $(\phi^{-1}_\pi\varphi_-\phi_\pi,\ \phi^{-1}_0\varphi_+\phi_0)$ where $\phi_t=\alpha^{-1}_t\varphi_t\alpha_{-t}\in(\Aut\g)^{\varrho_0}$ is as above. Since $\varphi_\pm^2$ and $\tilde\varphi_\pm^2$ are conjugate we may assume $\tilde\varrho_0=\varrho_0$. This implies $\alpha,\beta\in(\Aut\g)^{\varrho_0}$.

We now try to find an automorphism $\psi:L(\g,\sigma)\to L(\g,\tilde\sigma)$ of the form $\psi u(t)=\psi_t(u(t))$ with $\psi\varphi\psi^{-1}=\tilde\varphi$. This amounts to find a smooth curve $\psi_t$ in $\Aut\g$ with
\begin{center}
$(1)\ \psi_{t+2\pi}=\tilde\sigma\psi_t\sigma^{-1}$ and (2) $\psi_{-t}=\tilde\varphi_{t}^{-1}\psi_t\varphi_t$. 
\end{center}
From (2) we get $\psi_t=\tilde\varphi_{-t}^{-1}\tilde\varphi_t^{-1}\psi_t\varphi_t\varphi_{-t}$ and hence $\chi_t:=\tilde\alpha_t^{-1}\psi_t\alpha_t\in(\Aut\g)^{\varrho_0}$. Working with $\chi_t$ instead of $\psi_t$, (1) and (2) are equivalent to 
\begin{center}
(1') $\chi_{t+2\pi}=\tilde\sigma_t\chi_t\sigma_t^{-1}$ and (2') $\chi_{-t}=\tilde\phi_t^{-1}\chi_t\phi_t$
\end{center}
 where $\sigma_t:=\alpha^{-1}_{t+2\pi}\sigma\alpha_t$ and $\phi_t=\alpha^{-1}_t\varphi_t\alpha_{-t}$, and $\tilde\sigma_t,\tilde\phi_t$ are correspondingly defined. The advantage of these equations over (1) and (2) is that (2') at $t$ implies (2') at $-t$ because of $\phi_t\phi_{-t}=\varrho_0=\tilde\phi_t\tilde\phi_{-t}$. The idea now is to define $\chi_t$ in $[0,\pi]$, to extend it by (2') to $[-\pi,\pi]$ and then by (1') to all of $\R$. To make this well defined we have to ensure that $\chi_0=\tilde\phi_0^{-1}\chi_0\phi_0$ and $\chi_\pi=\tilde\sigma_{-\pi}\chi_{-\pi}\sigma^{-1}_{-\pi}$, the last equation being equivalent to $\chi_\pi=\tilde\sigma_{-\pi}\tilde\phi^{-1}_\pi\chi_\pi\phi_\pi\sigma^{-1}_{-\pi}$. Both conditions can be matched by choosing $\chi_0:=\alpha$ and $\chi_\pi:=\beta$ where $\alpha,\beta\in(\Aut\g)^{\varrho_0}$ are as above. Since $\alpha$ and $\beta$ lie in the same connected component we can smoothly connect them by a curve in $(\Aut\g)^{\varrho_0}$. By the above described extension we thus get $\chi_t$ in $(\Aut\g)^{\varrho_0}$ for all $t\in\R$ and this satisfies (1') and (2'). Moreover it is smooth except possibly at the integer multiples of $\pi$. But this can also be achieved by choosing $\chi_t$ more carefully around $0$ and $\pi$ according to the following lemma (at $\pi$ one has to apply it to  $\phi_{\pi-t}\sigma^{-1}_{t-\pi}$ and $\tilde\phi_{\pi-t}\tilde\sigma^{-1}_{t-\pi}$
  in order to get a smooth solution of
$\chi_{\pi+t}=\tilde\sigma_{t-\pi}\tilde\phi^{-1}_{\pi-t}\chi_{\pi-t}\phi_{\pi-t}\sigma^{-1}_{t-\pi}$).
\end{proof}

\begin{lemma}
Let $\phi_t$ and $\tilde\phi_t$ be smooth curves in $(\Aut\g)^{\varrho_0}$ with $\phi_t\phi_{-t}=\varrho_0=\tilde\phi_t\tilde\phi_{-t}$ for small $t$ and $\tilde\phi_0=\alpha\phi_0\alpha^{-1}$ for some $\alpha\in(\Aut\g)^{\varrho_0}$. Then $\chi_{-t}=\tilde\phi_t^{-1}\chi_t\phi_t$ has a smooth solution near $t=0$ with $\chi_0=\alpha$.
\end{lemma}

\begin{proof}
Let $\phi_t:=\phi_0e^{\ad x_t},\ \tilde\phi_t:=\tilde\phi_0e^{\ad \tilde x_t}$ near $t=0$ with smooth $x_t,\tilde x_t\in\g^{\varrho_0}$ and $x_0=\tilde x_0=0$. Then $\phi_t\phi_{-t}=\varrho_0=\tilde\phi_t\tilde\phi_{-t}$ implies $x_{-t}=-\phi_0x_t$ and $\tilde x_{-t}=\tilde\phi_0\tilde x_t$. Now $\chi_t:=e^{(-{1\over2}\ad \tilde x_{-t}}\cdot\alpha\cdot e^{{1\over2}\ad x_{-t}}$ is a solution.
\end{proof}

For any $\sigma\in\Aut\g$, let
\[
\J_2^{2q}(\g,\sigma):=\{[\varphi_+,\varphi_-]\in\J_2^{2q}(\g)\mid\varphi_-^{-1}\varphi_+ \mbox{ is conjugate to  $\sigma$ in } \Aut\g/\Int\g\}.
\]
Then  the elements of $\J_2^{2q}(\g,\sigma)$ are due to \ref{5.5} precisely the invariants of elements of $\Aut_2^{2q}L(\g,\sigma)$ and as in the case of automorphisms of the first kind we have:

\begin{corollary}\label{5.9}
There are natural bijections
\[
\Aut^{2q}_2\hat L(\g,\sigma)/\Aut \hat L(\g,\sigma)\leftrightarrow \Aut^{2q}_2 L(\g,\sigma)/\Aut L(\g,\sigma)\leftrightarrow \J^{2q}_2(\g,\sigma)\ .
\]
\end{corollary}

\begin{corollary}\label{5.10}
Any element of $\Aut^{2q}_2\hat L(\g,\sigma)$, resp. $\Aut^{2q}_2 L(\g,\sigma)$ is quasiconjugate to some $\hat\varphi\in\Aut_2^{2q}\hat L(\g,\tilde\sigma)$ with $\hat\varphi c=-c,\ \hat\varphi d=-d,\ \hat\varphi u\in L(\g,\tilde\sigma)$ and $\hat\varphi u(t)=\varphi_0(u(-t))$, resp.~to some $\varphi\in\Aut_2^{2q}L(\g,\tilde\sigma)$ with $\varphi u(t)=\varphi_0(u(-t))$ where $\varphi_0\in\Aut\g$ is constant.
\end{corollary}

\begin{remark}\label{5.11}\rm
Also in this case it is not always possible to {\rm conjugate} $\varphi\in\Aut^{2q}_2L(\g,\sigma)$ (or $\hat\varphi\in\Aut_2^{2q}\hat L(\g,\sigma)$) to one with $\varphi_t$ constant. For example an involution $\varphi$ of the second kind on $L(\g)$ with invariant $[\varphi_+,\varphi_-]$ where $\varphi_\pm\in\Int\g$ and $\varphi_+$ is not conjugate to $\varphi_-$ in $\Int\g$ is not conjugate to a $\tilde\varphi\in\Aut^2_2L(\g)$ with $\tilde\varphi_t$ constant.
\end{remark}

\begin{remark}\rm
The mapping $\varphi\mapsto\varphi^2$ from $\Aut^{2q}_2 L(\g,\sigma)$ to $\Aut^{q}_1 L(\g,\sigma)$ induces the mapping $\J_2^{2q}(\g,\sigma)\to\J_1^{q}(\g,\sigma),\ [\varphi_+,\varphi_-]\mapsto (0,\varrho,[\varphi_-^{-1}\varphi_+])$ where $\varrho\in{\cal A}^q$ is conjugate to $\varphi_+^2=\varphi^2_-$.
\end{remark}

\section{Involutions}

In this chapter we specialize our results to involutions (automorphisms of order $2$) and derive an explicit classification of all involutions of $L(\g,\sigma)$ and $\hat L(\g,\sigma)$ up to conjugation (more generally quasiconjugation) in case $\g$ is compact or complex. Actually we may restrict ourselves to the compact case and moreover to the case of loop algebras. For the classifications in the other cases  are in a natural bijection with this by Theorem~\ref{7.5} of the next chapter, Corollary~\ref{3.3} and Proposition~\ref{3.8}. Thus let $\g$ be a real compact simple Lie algebra and hence of type $\a_n(n\ge 1)$, $\b_n(n\ge 2)$, $\c_n(n\ge 3)$, $\d_n(n\ge 4)$, $\e_6$, $\e_7$, $\e_8$, $\f_4$ or $\g_2$. The group $\pi_0(\Aut\g)=\Aut\g/\Int\g$ of connected components of $\Aut\g$ is isomorphic to the group of isomorphisms of the Dynkin diagram of $\g$ and therefore isomorphic to $1,\Z_2$ or the symmetric group $S_3$. In particular the conjugacy classes of $\pi_0(\Aut\g)$ are determined by their order. Hence according to \ref {2.13}, $L(\g,\sigma)$ and $L(\g,\tilde\sigma)$ are isomorphic if and only if $\bar o(\sigma)=\bar o(\tilde\sigma)$ where $\bar o(\sigma)$ denotes the order of $\sigma$ in $\Aut\g/\Int\g$ (i.e. the smallest $k$ such that $\sigma^k\in\Int\g)$. We denote by $\g^{(k)}$ (resp. $\hat\g^{(k)})$ any $L(\g,\sigma)$ (resp. $\hat L(\g,\sigma))$ with $\bar o(\sigma)=k$.

\subsection{Involutions of the first kind}
By the general results of chapter 4, involutions of the first kind are up to quasiconjugation with isomorphisms of the first kind (actually with arbitrary isomorphisms, see below) in bijective correspondence with the set of triples
\[
(p,\varrho,[\beta])
\]
where $p\in\{0,1\}$ and $\varrho\in\Aut^r\g$ represents a conjugacy class of automorphisms of order $r:=(p,2)$, the greatest common divisor of $p$ and $2,$ where $\beta\in(\Aut\g)^\varrho$ and $[\beta]$ denotes the conjugacy class of $\beta\cdot((\Aut\g)^\varrho)_0$ in $\pi_0((\Aut\g)^\varrho)$. We call involutions of the first kind with $p=0$ (resp. $p=1)$ of type 1a (resp. 1b). In case 1a, $r=2$ and $\varrho$ is an involution while in case 1b, $r=1$ and  $\varrho=\id$.

Thus the list of conjugacy classes of \textbf{involutions of type 1a} is  a certain refinement of Cartan's list of involutions $\varrho$ by the conjugacy classes of $\pi_0((\Aut\g)^\varrho)$.
We therefore start by recalling Cartan's list, thereby fixing  a list ${\cal A}^2(\g)$  of representatives of conjugacy classes of involutions on each simple $\g$. This will be used throughout this chapter and will be called the standard list of involutions. For later use we also indicate  the outer involutions in this list. Note that $\varrho$ is outer if and only if the rank of $\k:=\g^\varrho$ is less than the rank of $\g$.

\fontsize{12}{3ex}\selectfont
\vskip 1 cm
\begin{center}
{STANDARD LIST OF INVOLUTIONS ${\cal A}^2(\g)$}\\
\end{center}
\vskip -0,5cm
\begin{center}
\begin{tabular}{||ll|l|l||}
\hline\hline
\hspace{2cm}$\g$& &\hspace{3cm}$\varrho$	&  $\hspace{1cm}\varrho$ outer\\
\hline
$\a_1=\su(2)$ && $\varrho_1:=\Ad  \tau_1$ & -\\
$\a_{2n}=\su(2n+1)$ & $(n\ge 1)$ & $\varrho_p:=\Ad \tau_p\ (1\le p\le n),\varrho_{n+1}:=\mu$ & $\varrho_{n+1}$\\
$\a_{2n-1}=\su(2n)$ & $ (n\ge 2)$ & $\varrho_p:=\Ad \tau_p\ (1\le p\le n), \varrho_{n+1}:=\mu,$& $\varrho_{n+1}, \varrho_{n+2}$\\
&&$\varrho_{n+2}:=\mu \Ad J$ &\\
$\b_n=\so(2n+1)$ & $ (n\ge 2)$ & $\varrho_p:=\Ad \tau_p\ (1\le p\le n)$ & -\\
$\c_{n}=\sp(n) $ & $(n\ge 3)$ & $\varrho_p:=\Ad \tau_p\ (1\le p\le [{n\over 2}]),\varrho_{[{n\over2}]+1}:=\Ad  i E_n$ & -\\
$\d_4=\so(8)$ && $\varrho_p:=\Ad \tau_p\ (1\le p\le 4)$ & $\varrho_1,\varrho_3$\\
$\d_n=\so(2n)$ & $ (n\ge 5)$ & $\varrho_p:=\Ad \tau_p\ (1\le p\le n),\varrho_{n+1}:=\Ad J$ & $\varrho_p,(p$ odd, $1\le p\le n)$\\
$\e_6$ && $\varrho_1,\varrho_2,\varrho_3,\varrho_4$ & $\varrho_1,\varrho_4$\\
$\e_7$ && $\varrho_1,\varrho_2,\varrho_3$ & -\\
$\e_8$ && $\varrho_1,\varrho_2$ & -\\
$\f_4$ && $\varrho_1,\varrho_2$ & -\\
$\g_2$ && $\varrho_1$ & -\\
\hline\hline
\end{tabular}
\end{center}
\fontsize{12}{3ex}\selectfont
\vskip 1 cm

The notations are as follows. $\tau_p$ denotes the diagonal matrix of appropriate size whose first $p$ diagonal elements are $-1$ and whose other diagonal elements are $1$. $J$ denotes the matrix
$\left(
\begin{array}{cc}
&E_n\\
-E_n&
\end{array}\right)$
where $E_n$ is the $n\times n$ unit matrix and $\mu\in\Aut(\su(n))$ is complex conjugation. Note that in case $\a_1$, $\mu=\Ad J$ and $\Ad J$ is conjugate to $\Ad\tau_1$ while in case $\d_4$, $\Ad J$ is conjugate to $\Ad\tau_2$ (corresponding to Cartan's isomorphisms of symmetric spaces in low dimensions).
 The symplectic algebra $\sp(n)$ is viewed as a subalgebra of the $n\times n$ quaternionic matrices. If $\g=\e_6$ the involutions $\varrho_1,\dots,\varrho_4$ are chosen to commute, which is possible by Lemma~\ref{A.8} of Appendix~A. In case of the exceptional algebras the order of succession of the $\varrho_i$ is chosen in the standard way corresponding to Cartan's symmetric spaces E\ I - E\ XII. In particular - and only this will be used  - $\varrho_1$ and $\varrho_4$ are outer in case of $\e_6$ and the fixed point algebra of $\varrho_2$  in case of $\e_7$ is isomorphic to $\so(12)+\su(2)$.
\\
The group $\pi_0((\Aut\g)^\varrho)$ of connected components of $(\Aut\g)^\varrho$ has been determined by Cartan \cite{Car} and Takeuchi \cite{Tak}. It is isomorphic to $1,\Z_2,\Z_2\times\Z_2, \D_4$ (the dihedral group with 8 elements) or $S_4$ (the symmetric group in 4 letters). $\D_4$, which may be identified with the symmetries of the standard unit square, has 5 conjugacy classes: $\{\id\},\{-\id\}$, the reflections along the two axes, the reflections along the two diagonals, and the rotations by angles~$\pm \pi/2$. Also $S_4$ has 5 conjugacy classes, namely the sets of cycles of order $1$ to $4$ and the set of products $(a,b)(c,d)$ of two cycles with $\{a,b,c,d\}=\{1,\dots,4\}$. In particular any element in $\pi_0((\Aut\g)^\varrho)$ is conjugate to its inverse. According to \ref{4.9} and \ref{4.11} the classification of involutions of the first kind up to conjugation with \textit{arbitrary} automorphisms is therefore the same as up to conjugation only with automorphisms \textit{of the first kind}.

The group $\pi_0((\Aut\g)^\varrho)$ contains the group $\pi_0((\Int\g)^\varrho)$ of connected components of $(\Int\g)^\varrho$ as a normal subgroup and is actually the semidirect product of this with a subgroup $F$ of $\Aut\g/\Int\g$. These groups are listed in \cite{Loo}\ p. 156, and by this list one can find representatives of the conjugacy classes of $\pi_0((\Aut\g)^\varrho)$ in most cases easily. The results are contained in Table~1.
\begin{table}
\fontsize{10}{2ex}\selectfont
\begin{minipage}{\linewidth}

\begin{center}
{TABLE 1}\\
REPRESENTATIVES $\sigma$ OF CONJUGACY CLASSES OF $\pi_0((\Aut\g)^\varrho)$
\end{center}

\begin{center}
\begin{tabular}{{||l|l|l|l|l|l|l||}}
\hline\hline
\hspace{1cm} $\g$        & \hspace{1cm} $\varrho$                  &$\pi_0((\Int\g)^\varrho)$	& $\pi_0((\Aut\g)^\varrho)$ & \quad$\sigma(k=1)$ & \quad  $\sigma(k=2)$ & \ $\sigma(k=3)$\\
\hline
$\a_1$                 & $\varrho_1$                & $\Z_2$ & $\Z_2$     &  $\id ,\mu$  &    -            & \hspace{0,6cm}- \\
$\a_{2n}$\hspace{0,3cm}   $(n\ge1)$ & $\varrho_p\ (1\le p\le n+1)$ & $1$    & $\Z_2$     &  $\id $      & $\varrho_{n+1}$ &\hspace{0,6cm}- \\
$\a_{2n-1}$  $(n\ge2)$ & $\varrho_p\ (1\le p<n)$     & $1$    & $\Z_2$     &  $\id $      & $\varrho_{n+1}$ & \hspace{0,6cm}- \\
&           $\varrho_n$                & $\Z_2$ & $\Z_2\times\Z_2$        & $\id ,\Ad J$  & $\varrho_{n+1},\varrho_{n+2}$ & \hspace{0,6cm}- \\
&           $\varrho_{n+1}$ & $\Z_2$ & $\Z_2\times\Z_2$ & $\id ,\varrho_1$ & $ \varrho_{n+1},\varrho_1\varrho_{n+1}$ & \hspace{0,6cm}-\\
&           $\varrho_{n+2}$            & $1$    &$\Z_2$                  & $\id $       & $\varrho_{n+2}$ & \hspace{0,6cm}- \\
$\b_n$   \hspace{0,4cm}  $(n\ge 2)$  & $\varrho_p\ (1\le p\le n)$ & $\Z_2$  & $\Z_2$   & $\id ,\varrho_1 \Ad \tau_{p+1}$ &- &\hspace{0,6cm}- \\
$\c_{2n-1}$  $(n\ge 2)$ & $\varrho_p \ (1\le p<n)$ & $1$ & $1$ & $\id $ & - & \hspace{0,6cm}-\\
&           $\varrho_n$ & $\Z_2$ & $\Z_2$ & $\id ,\Ad jE$ & - & \hspace{0,6cm}- \\
$\c_{2n}$ \hspace{0,3cm} $(n\ge 2) $ &$\varrho_p\ (1\le p<n)$ & $1 $ & $1$ & $\id $ & - &\hspace{0,6cm}- \\
&           $\varrho_n$ & $\Z_2$ & $\Z_2$ & $\id ,\Ad J$ & - & \hspace{0,6cm}- \\
&           $\varrho_{n+1}$ & $\Z_2$ & $ \Z_2$ & $\id ,\Ad jE$ & - & \hspace{0,6cm}-\\
$\d_4$       &$\varrho_p$ ($p=1$ or $3$) & $1$ & $\Z_2$ & $\id $ &$\varrho_p$ & \hspace{0,6cm}- \\
&           $\varrho_2$ & $\Z_2$ & $\Z_2\times\Z_2$ & $\id ,\varrho_1\varrho_3$ & $\varrho_1,\varrho_3$ & \hspace{0,6cm}-\\
&           $\varrho_4$ & $\Z_2\times\Z_2$ & $S_4$ & $\id ,\Ad J$ & $\varrho_1,\varrho_1 \Ad J$ &\hspace{0,5cm}$\vartheta$\\
$\d_{2n}$ \hspace{0,2cm}   $(n\ge 3)$  & $\varrho_p\ (1\le p<2n)$ &  &&&&\\
&           \quad $p$ even  &$\Z_2$ &  $\Z_2\times \Z_2$ &$\id ,\varrho_1\varrho_{p+1}$ & $\varrho_1,\varrho_{p+1}$ & \hspace{0,6cm}-\\ 
&           \quad $p$ odd   & $1$ & $\Z_2$ & $\id $ & $\varrho_p$ & \hspace{0,6cm}- \\
&           $\varrho_{2n}$ & $\Z_2\times\Z_2$ & $\D_4$ & $\id ,\varrho_1\Ad \tau_{2n+1},$&$\varrho_1,\varrho_1 \varrho_{2n+1}$ & \hspace{0,6cm}-\\
&           &&&$\varrho_{2n+1}$ & &\\
&           $\varrho_{2n+1}$ & $\Z_2$ & $\Z_2$ & $\id ,\varrho_{2n}$ & - &\hspace{0,6cm}- \\
$\d_{2n-1}$ $(n\ge 3)$ & $\varrho_p\ (1\le p <2n-1)$ & & & & &\\
&           \quad $p$ even & $\Z_2$ & $\Z_2\times\Z_2$ & $\id ,\varrho_1\varrho_{p+1}$ & $\varrho_1,\varrho_{p+1}$ &\hspace{0,6cm}- \\
&           \quad $p$ odd & $1$ & $\Z_2$ & $\id $ & $\varrho_p$ & \hspace{0,6cm}-\\
&           $\varrho_{2n-1}$ & $\Z_2$ & $\Z_2\times\Z_2$ & $\id ,\varrho_{2n}$ & $\varrho_{2n-1},\varrho_{2n-1}\varrho_{2n}$ & \hspace{0,6cm}-\\
&           $\varrho_{2n}$ & $1$ & $\Z_2$ & $\id $ & $\varrho_{2n-1}$ & \hspace{0,6cm}-\\
$\e_6$       & $\varrho_1,\varrho_2,\varrho_3,\varrho_4$ &$ 1 $ & $\Z_2$ & $\id $ & $\varrho_1$ & \hspace{0,6cm}-\\
$\e_7$       & $\varrho_p$ ($p=1$ or $3$) & $\Z_2$ & $\Z_2$ & $\id ,\sigma_p$ &- &\hspace{0,6cm}-\\
&           $\varrho_2$ & $1$ & $1$ & $\id $ & - &\hspace{0,6cm}- \\
$\e_8$       & $\varrho_1,\varrho_2$ & $1$ & $1$ & $\id $ & - &\hspace{0,6cm}-\\
$\f_4$       & $\varrho_1,\varrho_2$ & $1$ &  $1$ & $\id $ & - & \hspace{0,6cm}-\\
$\g_2$       & $\varrho_1$ & $1$ & $1$ & $\id $ & - & \hspace{0,6cm}- \\
\hline\hline
\end{tabular}
\end{center}

\fontsize{9}{2ex}\selectfont
{\it {The notations are that of the Standard List. In case $\e_7,\ \sigma_p=e^{adX}$ for some $X\in\g$ with $\varrho_pX=-X, X\ne 0$, and $(ad X)^3=-\pi^2(adX)$.}}
\end{minipage}
\fontsize{12}{3ex}\selectfont
\end{table}
For more details (including a simplified calculation of these groups) we refer to Appendix~A. We denote the representatives of the conjugacy classes of $\pi_0((\Aut\g)^\varrho)$ in this list by $\sigma$ and their order in $\Aut\g/\Int\g$ by $k$ (i.e. $k=\bar o(\sigma))$. In case of $\g=\e_7$ we let $\sigma_p:=e^{\ad  X_p}$ where $X_p\in\g$ is an arbitrary non zero element with $\varrho_p X_p=-X_p$ and $(\ad X_p)^3=-\pi^2\ \ad X_p$ (cf.~Appendix~A). 

An involution corresponding to $(\varrho,\sigma)$  is given by $u(t)\mapsto\varrho(u(t))$ on $L(\g,\sigma)\cong\g^{(k)}$. It has fixed point algebra $L(\k,\sigma_{\vert_\k})$, where $\k=\g^\varrho$, and extends to $\hat L(\g,\sigma)$ by $c\mapsto c, d\mapsto d$.

By simply rearranging Table 1 we get the classification of conjugacy classes of involutions of type 1a given in Table 2. It lists under $\g^{(k)}$ all pairs $(\varrho,\sigma)$ from Table 1 where $\varrho$ is from the standard list of involutions  on $\g$ and $\sigma$ represents a conjugacy class of $\pi_0((\Aut\g)^\varrho)$ with $\bar o(\sigma)=k$.
\begin{table}
\fontsize{12}{3ex}\selectfont
\fontsize{10}{2ex}\selectfont
\begin{minipage}{\linewidth}
\begin{center}
{TABLE 2}\\
{INVOLUTIONS OF THE FIRST KIND}\\
\end{center}
\begin{center}

\begin{tabular}{{||l|l|l|l|l||}}
\hline\hline
\hspace{0,5cm} $\g^{(k)}$ & \hspace{0,5cm} $(\varrho,\sigma)$ (type 1a) &$\beta$ (type 1b)& \hspace{0,3cm} number\\
\hline
$\a^{(1)}_1$ & $(\varrho_1,\id ),(\varrho_1,\mu)$ & $\id $ & $2+1$\\
$\a^{(1)}_{2n}\hspace{0,3cm} (n\ge 2)$ & $(\varrho_p,\id ),\ 1\le p\le n+1$ & $\id ,\mu$ & $(n+1)+2$\\
$\a^{(1)}_{2n-1}\ (n\ge2)$ & $(\varrho_p, \id ),\ 1\le p\le n+2$ & $\id ,\mu$ &$(n+4)+2$\\
&         $(\varrho_n,\Ad J),(\varrho_{n+1},\varrho_1)$ &&\\
$\b^{(1)}_n\hspace{0,3cm} (n\ge2)$   & $(\varrho_p,\id ),\ 1\le p\le n$ & $\id $ &$2n+1$\\
&                     $(\varrho_p,\varrho_1\Ad \tau_{p+1}),\ 1\le p\le n$&&\\
$\c^{(1)}_{2n-1}\ (n\ge 2)$ & $(\varrho_p,\id ),\ 1\le p\le n$ & $\id $ & $(n+1)+1$\\
&                      $(\varrho_n,\Ad jE)$&&\\
$\c^{(1)}_{2n}\hspace{0,3cm} (n\ge 2)$ & $(\varrho_p,\id ),\ 1\le p\le n+1$ & $\id $ & $(n+3)+1$\\
&                      $(\varrho_n,\Ad J), (\varrho_{n+1}, \Ad jE)$&&\\
$\d^{(1)}_4$ & $(\varrho_p,\id ), \hspace{0,3cm} 1\le p\le 4$ & $\id ,\varrho_1$ & $6+2$\\
&                    $(\varrho_2,\varrho_1\varrho_3), (\varrho_4, \Ad J)$ &&\\
$\d^{(1)}_{2n}\hspace{0,3cm}(n\ge3)$  & $(\varrho_p,\id ),\ 1\le p\le 2n+1$ & $\id ,\varrho_1$ & $(3n+3)+2$\\
&                    $(\varrho_{2 l},\varrho_1 \Ad \tau_{2l+1}),\ 1\le l\le n$ & &\\
&                    $(\varrho_{2n},\varrho_{2n+1}), (\varrho_{2n+1},\varrho_{2n})$ & &\\
$\d_{2n-1}^{(1)}\hspace{0,1cm}(n\ge3)$ & $(\varrho_p,\id ),\ 1\le p\le 2n$ & $\id ,\varrho_1$ & $3n+2$\\
&                    $(\varrho_{2l},\varrho_1\varrho_{2l+1}),\ 1\le l\le n-1$ &&\\
&                    $(\varrho_{2n-1},\varrho_{2n})$&&\\
$\e^{(1)}_6$         & $(\varrho_p,\id ),\ 1\le p\le 4$ & $\id ,\varrho_1$ & $4+2$\\
$\e^{(1)}_7$         & $(\varrho_p,\id ),\ 1\le p\le 3$ & $\id $& $5+1$\\
&                     $(\varrho_1,\sigma_1),(\varrho_3,\sigma_3)$&&\\
$\e^{(1)}_8$         & $(\varrho_1,\id ), (\varrho_2,\id )$ & $\id $ & $2+1$\\
$\f_4^{(1)}$         & $(\varrho_1,\id ), (\varrho_2,\id )$ & $\id $ & $2+1$\\
$\g_2^{(1)}$         & $(\varrho_1,\id )$ & $\id $ & $1+1$\\
\hline
$\a^{(2)}_{2n}$      & $(\varrho_p,\varrho_{n+1}),\ 1\le p\le n+1$ & - & $(n+1)+0$\\
$\a^{(2)}_{2n-1}$    & $(\varrho_p,\varrho_{n+1}),\ 1\le p\le n+1$ & - & $(n+4)+0$\\
&                     $(\varrho_n,\varrho_{n+2}), (\varrho_{n+1},\varrho_1\varrho_{n+1}), (\varrho_{n+2},\varrho_{n+2})$ & &\\
$\d^{(2n)}_{2n}\hspace{0,3cm}(n\ge2)$ & $(\varrho_{2l-1},\varrho_{2l-1}),\ 1\le l\le n$ & - & $3n+0$\\
&                       $(\varrho_{2l},\varrho_1),\ 1\le l\le n$ & &\\
&                       $(\varrho_{2l},\varrho_{2l+1}),\ 1\le l\le n-1$ &&\\
&                       $(\varrho_{2n},\varrho_1 \varrho_{2n+1})$ &&\\
$\d^{(2)}_{2n-1}\hspace{0,2cm}(n\ge3)$ & $(\varrho_{2l-1},\varrho_{2l-1}),\ 1\le l\le n$ & - & $3n+0$\\
&                        $(\varrho_{2l},\varrho_1),\ 1\le l\le n-1)$ & &\\
&                        $(\varrho_{2l},\varrho_{2l+1}),\ 1\le l\le n-1$ &&\\
&                        $(\varrho_{2n-1},\varrho_{2n-1}\varrho_{2n}),(\varrho_{2n},\varrho_{2n-1})$&&\\
$\e^{(2)}_6$             & $(\varrho_p,\varrho_1),\ 1\le p\le 4$ &-&$4+0$\\
\hline
$\d^{(3)}_4$             & $(\varrho_4,\vartheta)$ & $\vartheta$ & $1+1$\\
\hline\hline
\end{tabular}
\end{center}
\end{minipage}
\fontsize{12}{3ex}\selectfont
\end{table}

\textbf{Involutions of type 1b} are classified by the conjugacy classes $[\beta]$ of $\Aut\g/\Int\g$ and are represented by $u(t)\mapsto\beta u(t+\pi)$ on $L(\g,\beta^{-2})\cong\g^{(k)}$ where $k=\bar o(\beta^{-2})$. In particular $k=1$ or $3$, and $3$ only occurs once, namely in case $\g=\so(8)$ and $\beta$ the triality automorphism. On $L(\g)\cong\g^{(1)}$ there are up to conjugation one or two involutions of type 1b depending on whether $\g$ admits no or one outer involution. The results are also listed  in Table 2.

The fixed point algebra of $u(t)\mapsto\beta u(t+\pi)$ is $\{u(t)\mid u(t+\pi)=\beta^{-1} u(t)\}\subset L(\g,\beta^{-2})$ and is thus isomorphic to $\g^{(l)}$ with $l=\bar o(\beta)$.

\begin{remark}{\rm By means of the fixed point algebras it can be checked easily that the above classification of involutions of the first kind is in bijection with that of Bausch and Rousseau (\cite{BR}, Tables p. 133 - 138) although the latter was obtained in the algebraic case, i.e. by working with algebraic instead of smooth loops. But we will give an a priori proof in chapter 8 that both classifications  coincide, thus obtaining  in particular a simplified proof of their classification.}
\end{remark}

\subsection{Involutions of the second kind}

By the results of chapter 5 the quasiconjugacy classes of involutions of the second kind are in bijection with the equivalence classes $[\varrho_+,\varrho_-]$ of pairs $(\varrho_+,\varrho_-)$ of automorphisms of $\g$ with $\varrho^2_\pm=\id$. Two equivalence classes $[\varrho_+,\varrho_-]$ and $[\tilde\varrho_+,\tilde\varrho_-]$ coincide if and only if $\tilde\varrho_+=\alpha\varrho_+\alpha^{-1}$ and $\tilde\varrho_-=\beta\varrho_-\beta^{-1}$ or $\tilde\varrho_+=\alpha\varrho_-\alpha^{-1}$ and $\tilde\varrho_-=\beta\varrho_+\beta^{-1}$ for some $\alpha,\beta\in\Aut\g$ with $\alpha^{-1}\beta\in\Int\g$. An involution corresponding to $[\varrho_+,\varrho_-]$ is given for example by $u(t)\mapsto\varrho_+(u(-t))$ on $L(\g,\sigma)\cong\g^{(k)}$ where $\sigma=\varrho_-\varrho_+$ and $k=\bar o(\varrho_-\varrho_+)$ is the order of $\varrho_-\varrho_+$ in $\Aut\g/\Int\g$. It extends to $\hat L(\g,\sigma)\cong\hat\g^{(k)}$ by $c\mapsto -c, d\mapsto -d$.

To determine the equivalence classes $[\varrho_+,\varrho_-]$ more explicitly we first describe  the involutions on each $\g$ up to conjugation with \emph{inner} automorphisms.

\begin{proposition}\label{6.2}
A list of representatives of involutions on each $\g$ up to conjugation with \emph{inner} automorphisms is obtained from the standard list as follows.
\begin{enumerate}
\item [(i)] If $\g$ is not isomorphic to $\d_{2m}\cong\so(4m)$ then one can take the same list (two involutions which are conjugate are also conjugate by an inner automorphism in this case).
\item[(ii)] If $\g=\so(4m)$ and $m\ge 3$ then one gets a complete list by adding $\varrho'_{2m+1}:=\varrho_1\varrho_{2m+1}\varrho_1$.
\item[(iii)] If $\g=\so(8)$ a list is given by $\varrho_p,\varrho_p',\varrho_p''\ (1\le p \le 3)$ and $\varrho_4$ where $\varrho'_p=\vartheta\varrho_p\vartheta^{-1},\varrho''_p=\vartheta^2\varrho_p\vartheta^{-2}$ and $\vartheta$ denotes the triality automorphism.
\end{enumerate}
\end{proposition}

\begin{proof} If $\varrho_1,\dots,\varrho_k\in\Aut(\g)$ are the involutions from the standard list and $\alpha_1:=id,\\ \alpha_2,\dots,\alpha_l\in\Aut\g$ are representatives of the elements of $\Aut\g/\Int\g$ then any involution is conjugate by an inner automorphism to at least one of the $\alpha_i\varrho_j\alpha_i^{-1}\ (1\le i\le l, 1\le j\le k)$. This proves in particular (i) since in all cases there, either $l=1$, or $l=2$ and $\alpha_2$ can be chosen to commute with the $\varrho_j$. In fact, in the cases $\su(n),\ \so(4m+2)$ and $\e_6$ one may take $\alpha_2=\mu, \Ad\tau_{2m+1}$, and $\varrho_1$, respectively.\\
(ii) If $\g=\so(4m)$ and $m\ge 3$ we also have $l=2$ and may take $\alpha_2:=\Ad\tau_1$, which commutes with all $\varrho_p$ except $\varrho_{2m+1}=\Ad J$. Moreover $\alpha_2\varrho_{2m+1}\alpha_2^{-1}=\Ad\tau_1J\tau_1$ is not conjugate to $\Ad J$ by an inner automorphism as $\{A\in O(4m) \mid AJA^{-1}=\pm J\}=U(2m)\cup
\left({\scriptsize
\begin{array}{cc}
E_{2m}&\\
&-E_{2m}
\end{array}}
\right)
\cdot U(2m)\subset SO(4m)$.

(iii) Finally, if $\g=\so(8)$ then $\Aut\g/\Int\g\cong S_3$ and for any $x,y\in\Aut\g/\Int\g$ with $\ord\ x=2$ and $\ord\ y=3$ one has $xyx^{-1}=y^{-1}$ and $\Aut\g/\Int\g=\{1,x,yx=y^2xy^{-2},y^2x=yxy^{-1},y,y^2\}$. This applies to $x=\bar\varrho_1=\bar\varrho_3$ and $y=\bar\vartheta$ where $\bar\alpha$ denotes the image of $\alpha\in\Aut\g$ in $\Aut\g/\Int\g$ and shows in particular that $\varrho_p,\varrho_p', \varrho_p''$ are pairwise not conjugate by an inner automorphism in case $p=1$ and $3$. Since $\varrho_4$ commutes with $\vartheta$ (cf. \cite{Loo}) a complete list of representatives can therefore be found among the $\varrho_p, \varrho_p',\varrho_p''\ (1\le p\le 3)$ and $\varrho_4$, and the only question is whether $\varrho'_2$ and $\varrho''_2$ have to be deleted. But this is not the case. For $\varrho_2=\Ad\tau_2$ and $\Ad J$ (which are conjugate by an outer automorphism) are not conjugate by an inner automorphism as $\tau_2$ and $\pm J$ have different eigenvalues. Thus $\varrho'_2$ and $\varrho_2''$ can not both   be deleted and from this it follows easily that neither of them can be deleted.
\end{proof}

\begin{theorem}
The conjugacy classes of involutions of the second kind on each $\g^{(k)}$ are in bijection with the pairs $[\varrho_+,\varrho_-]$ of Table 3. A representative corresponding to $[\varrho_+,\varrho_-]$ is given by $u(t)\mapsto\varrho_+(u(-t))$ on $L(\g,\varrho_-\varrho_+)$ and this extends to an automorphism of $\hat L(\g,\varrho_-\varrho_+)$ by $c\mapsto -c, d\mapsto-d$.
\end{theorem}
\begin{table}
\fontsize{10}{2ex}\selectfont
\begin{minipage}{\linewidth}
\begin{center}
{TABLE 3}\\
INVOLUTIONS OF THE SECOND KIND
\end{center}

\begin{center}

\begin{tabular}{||ll|l|l||}
\hline\hline
$\g^{(k)}$&	& $[\varrho_+,\varrho_-]$	& number\\
\hline
$\a_1^{(1)}$	& & $[\varrho_p,\varrho_q]\  (0\le p\le q\le 1)$	& 3\\
$\a_{2n}^{(1)}$ & $(n\ge 1)$	&$[\varrho_p,\varrho_q] \ (0\le p\le q\le n$ or $p=q=n+1)$	&${1\over2}n(n+3)+2$\\
$\a^{(1)}_{2n-1}$ & $ (n\ge 2)$	&$[\varrho_p,\varrho_q]\ (0\le p\le q\le n$ or $n+1\le p\le q\le n+2)$ & ${1\over2}n(n+3)+4$\\
$\b_n^{(1)}$ & $ (n\ge 2)$	& $[\varrho_p,\varrho_q]\ (0\le p\le q\le n)$ & ${1\over2}(n+1)(n+2)$\\
$\c_{2n}^{(1)}$ & $ (n\ge 2)$ & $[\varrho_p,\varrho_q]\ (0\le p\le q\le n+1)$ & ${1\over2}(n+2)(n+3)$\\
$\c_{2n-1}^{(1)}$ & $ (n\ge 2)$ & $[\varrho_p,\varrho_q]\ (0\le p\le q\le n)$ & ${1\over2}(n+1)(n+2)$\\
$\d_4^{(1)}$& & $[\varrho_p,\varrho_q]\ (0\le p\le q\le 4, p+q$ even), $[\varrho_2,\varrho_2']$ & $10$\\
$\d_{2n}^{(1)}$ & $ (n\ge 3)$ & $[\varrho_p,\varrho_q] \ (0\le p\le q\le 2n, p+q$ even), & $n^2+3n+4$\\
&&$[\varrho_p,\varrho_{2n+1}]\ (0\le q\le 2n, p$ even), &\\
&&$[\varrho_{2n+1},\varrho_{2n+1}], [\varrho_{2n+1}, \varrho'_{2n+1}]$ &\\
$\d^{(1)}_{2n-1}$ & $ (n\ge 3)$ & $[\varrho_p,\varrho_q]\ (0\le p\le q\le 2n+1, p+q$ even) & $(n+1)^2$\\
$\e^{(1)}_6$ && $[\varrho_p,\varrho_p]\ (0\le p\le 4),[\varrho_0,\varrho_2],[\varrho_0,\varrho_3],[\varrho_1,\varrho_4],[\varrho_2,\varrho_3]$&$9$\\
$\e_7^{(1)}$ && $[\varrho_p,\varrho_q]\ (0\le p\le q\le 3)$& $10$\\
$\e_8^{(1)}$ && $[\varrho_p,\varrho_q]\ (0\le p\le q\le 2)$& $6$\\
$\f_4^{(1)}$ && $[\varrho_p,\varrho_q]\ (0\le p\le q\le 2)$ & $6$\\
$\g_2^{(1)}$ && $[\varrho_p,\varrho_q]\ (0\le p\le q\le 1)$ & $3$\\
\hline
$\a_{2n}^{(2)}$ & $ (n\ge 1)$ & $[\varrho_p,\varrho_{n+1}] \ (0\le p\le n)$ & $n+1$\\
$\a_{2n-1}^{(2)}$ & $ (n\ge 2)$ & $[\varrho_p,\varrho_{n+1}]\ (0\le p\le n),[\varrho_p,\varrho_{n+2}](0\le p\le n)$ & $2(n+1)$\\
$\d_4^{(2)}$&& $[\varrho_p,\varrho_q]\ (0\le p\le q\le 4,p+q$ odd), $[\varrho_1,\varrho_2'],[\varrho_2,\varrho_3']$& $8$\\
$\d_{2n}^{(2)}$ & $ (n\ge 3)$ & $[\varrho_p,\varrho_q]\  (0\le p\le q\le 2n, p+q$ odd), & $n(n+2)$\\
&&$[\varrho_p,\varrho_{2n+1}]\ (0\le p\le 2n,p$ odd) &\\
$\d_{2n-1}^{(2)}$ & $ (n\ge 3)$ & $[\varrho_p,\varrho_q] \ (0\le p\le q\le 2n, p+q$ odd) & $n(n+1)$\\
$\e_6^{(2)}$ && $[\varrho_0,\varrho_1],[\varrho_0,\varrho_4],[\varrho_1,\varrho_2][\varrho_1,\varrho_3],[\varrho_2,\varrho_4][\varrho_3,\varrho_4]$& $6$\\
\hline
$\d_4^{(3)}$& & $[\varrho_1,\varrho_1'],[\varrho_1,\varrho_3'],[\varrho_3,\varrho_3']$ & $3$\\
\hline\hline
\end{tabular}
\end{center}

\fontsize{9}{2ex}\selectfont
{The notations are that of the standard list with the additional conventions $\varrho_0:=\id$, $\varrho'_{2n+1}=\varrho_1\varrho_{2n+1}\varrho_1$ if $\g=\d_{2n}$, $n\ge 3$, and $\varrho'_p=\vartheta\varrho_p\vartheta^{-1}$ if $\g=\d_4$. A conjugacy class corresponding to $[\varrho_+,\varrho_-]$ is represented by $u(t)\mapsto \varrho_+(u(-t))$ on $L(\g,\varrho_-\varrho_+)\cong\g^{(k)}$ where $k$ is the smallest positive integer such that $(\varrho_-\varrho_+)^k\in \Int\g$. This representative extends to $\hat\g^{(k)}$ by $c\mapsto -c$, $d\mapsto -d$.}
\end{minipage}
\fontsize{12}{3ex}\selectfont
\end{table}
\begin{proof} By the discussion above we have   to determine the equivalence classes $[\varrho_+,\varrho_-]$ of the $\varrho_\pm\in\Aut\g$ with $\varrho_\pm^2=\id$, and $k:=\bar o(\varrho_-\varrho_+)$. This can be done easily in most cases by means of Cartan's list in combination with Proposition \ref{6.2}. Each such $[\varrho_+,\varrho_-]$ then corresponds to a unique conjugacy class of involutions of the second kind on $\g^{(k)}$ and vice versa.

More precisely, if $\g$ is not isomorphic to $\so(4m)$ then $\varrho_+$ and $\varrho_-$ may be taken from the standard list (say $\{\varrho_1,\dots,\varrho_r\})$ enlarged by $\varrho_0:= \id$. Hence the $[\varrho_p,\varrho_q]$ with $0\le p\le q\le r$ represent all equivalence classes $[\varrho_+,\varrho_-]$, and without repetition. Moreover $k=1$ if $\varrho_-\varrho_+$ is inner and $k=2$ otherwise.
\\
If $\g=\so(4m)$ and $m\ge3$ then $\varrho_\pm$ may be chosen  from $\{\varrho_p=\Ad \tau_p\mid 0\le p\le 2m\}\cup\{\varrho_{2m+1}=AdJ, \varrho'_{2m+1}=\Ad \tau_1J\tau_1\}$. Since $\Ad \tau_1$ commutes with  $\varrho_p$ if $p\le 2m$ (and hence $[\varrho_p,\varrho'_{2m+1}]=[\varrho_p,\varrho_{2m+1}]$), a complete list of equivalence classes is given by the $[\varrho_p,\varrho_q]$ with $0\le p\le q\le 2m+1$ \ together with $[\varrho_{2m+1},\varrho'_{2m+1}]$. Again, $k=1$ if $\varrho_-\varrho_+$ is inner and $k=2$ otherwise.
\\
The most interesting case is $\g=\so(8)$. Here we may take $\varrho_\pm$ from $\{\varrho_0=id, \varrho_4\}\cup{\bigcup\limits^3_{p=1}}\{\varrho_p,\varrho_p',\varrho_p''\}$, and $\varrho_+$ actually from $\{\varrho_p\mid 0\le p\le 4\}$ as we may conjugate $\varrho_+,\varrho_-$ simultaneously by $\vartheta$ or $\vartheta^2$. Since $\varrho_0$ and $\varrho_4$ commute with $\vartheta$, since $[\varrho_p,\varrho'_q]=[\vartheta^2\varrho_p\vartheta^{-2},\vartheta^2\varrho_q'\vartheta^{-2}]=[\varrho''_p,\varrho_q]=[\varrho_q,\varrho_p'']$, and $[\varrho_p,\varrho_q']=[\vartheta\varrho_1\varrho_p\varrho_1\vartheta^{-1},\vartheta\varrho_1\vartheta\varrho_q\vartheta^{-1}\varrho_1\vartheta^{-1}]=[\varrho'_p,\varrho_q]$ (note $\varrho_1\vartheta=\vartheta^{-1}\alpha\varrho_1$ for some $\alpha\in\Int\ \so(8))$, a complete list of representatives is given by the $[\varrho_p,\varrho_q]$ with  $0\le p\le q\le 4$ together with the $[\varrho_p,\varrho'_q]$ with $1\le p\le q\le 3$. Also note that $[\varrho_p,\varrho_q]\ne [\varrho_p,\varrho'_q]$ as $\varrho_q$ and $\varrho'_q$ are not conjugate by an inner automorphism. If $[\varrho_+,\varrho_-]=[\varrho_p,\varrho_q]$ then $k=1$ or $k=2$ according to $p+q$ being even or odd. Since $\varrho_2$ and $\varrho_2'$ are inner the $k$ of $[\varrho_2,\varrho'_p]$ and $[\varrho_p,\varrho_2']$ is $1$ if $p$ is even and $2$ if $p$ is odd. Finally the $[\varrho_p,\varrho_q']$ with $p,q\in\{1,3\}$ have $k=3$ as $\bar\vartheta\bar\varrho_q\bar\vartheta^{-1}\bar\varrho_p=\bar\vartheta^2\bar\varrho_q\bar\varrho_p=\bar\vartheta^2$ in $\Aut\g/\Int\g$.
\end{proof}

\begin{remark}{\rm There is a close connection between involutions of the second kind on affine Kac-Moody algebras and Hermann examples of hyperpolar actions on compact Lie groups. An Hermann example of a hyperpolar action is the action of $K_+\times K_-$ on $G$ by $(k_+,k_-). g=k_+gk_-^{-1}$ where $G$ is a compact Lie group and $K_\pm$ are symmetric subgroups, i.e. (open subgroups of) fixed point groups of involutions $\varrho_\pm$ on $G$. This action is hyperpolar in the sense that there exists a torus in $G$ which meets every orbit and always orthogonally. Kollross \cite{Kol} has classified hyperpolar actions on compact simple, simply connected Lie groups $G$ and proved that they are either Hermann examples, $\sigma$ actions (the action of $\{(g,\sigma(g))\mid g\in G\}$ on $G$ for some $\sigma\in\Aut \ G)$ or cohomogeneity $1$ actions. Moreover he has classified Hermann actions up to a natural equivalence (cf. also \cite{MatT}) and his list coincides with our list of equivalence classes $[\varrho_+,\varrho_-]$ with $\varrho_\pm\ne \id$. Thus there is (almost) a bijection between Hermann examples and involutions of the second kind. In \cite{HPTT} this has been already observed for a special class of involutions, namely those which leave the subspace $\R c+\R d$ of a Kac-Moody algebra $\hat L(\g,\sigma)$ invariant; but as we proved above, any involution is quasiconjugate to such a special one. An explanation for this surprising bijection lies in Terng's construction of $P(G,H)$ actions on Hilberts spaces which associates to each hyperpolar action of $H\subset G\times G$ on $G$ the action of the group $P(G,H)$ of $H^1$-curves in $G$ with endpoints in $H$
on the Hilbert space $L^2([0,1],\g)$ by gauge transformations. Under this mapping the Hermann example of $K_+\times K_-$ on $G$ with $K_\pm=G^{\varrho_\pm}$ corresponds to the isotropy representation of the ''Kac-Moody symmetric space`` defined by the involution of the second kind with invariant $[\varrho_+,\varrho_-]$ while the $\sigma$-actions correspond to the group case, i.e. to the adjoint action of the associated Kac-Moody group on $\hat L(\g,\sigma)$ (cf. \cite{HPTT} and \cite{Hei1}). In \cite{Gro} it has been shown that the isotropy representation of any Kac-Moody symmetric space is (hyper)polar.}
\end{remark}

\begin{remark}{\rm Involutions of the second kind have been classified in the algebraic case by~\cite[Tables I, II, p.~85 -- 94]{B$_3$R}. These authors obtained exactly the same number of involutions of the second kind in each case as we do, except for $D^{(1)}_4$ where they seem to have overlooked some redundancies and obtained $14$ instead of $10$. Actually they classified almost split real forms of complex affine Kac-Moody algebras. But as we will see later these are in bijection with involutions of the second kind (Chapter~7) and the classifications coincide in the smooth and algebraic case (Chapter~8).}
\end{remark}

%%%%%%%%%%%%%%%%%%%%%%%%%%%%%%%%%%%%%%%%%%%%%%%%%%%%%%%%%%
\section{Real forms}
In this chapter let $\g$ be a {\sl complex} simple Lie algebra, $\sigma\in\Aut\g$ and $\GC$ either the loop algebra $L(\g,\sigma)$ or the affine Kac-Moody algebra $\hat L(\g,\sigma)$.

Our purpose is to show that real forms of $\GC$ correspond, like in finite dimensions, to involutions of a ''compact real form`` $\U$ of $\GC$ and are hence classified by the results of the last chapter, and furthermore that each real form of $\GC$ has a Cartan decomposition which is unique up to conjugation. To this end we prove that extensions of automorphisms of finite order of $\U$ to linear, resp. conjugate linear automorphisms of $\GC$ induce bijections  between their conjugacy classes, also as in finite dimensions. Some relevant finite dimensional background material is given in Appendix~B.

We first note that in the {\sl complex} case (as well as in the {\sl compact case}), $\J_i^q(\g,\sigma)$ only depends on the order $\bar o(\sigma)$ of $[\sigma]$ in $\Aut\g/\Int\g$ (cf.~\ref{2.14} and the introduction of Chapter~6). We therefore define $\J_i^q(\g,k)$ to be $\J_1^q(\g,\sigma)$ if $\bar o(\sigma)=k$.

Now let $\varphi\in\overline{\Aut}L(\g,\sigma)$ or $\hat\varphi\in\overline{\Aut}\hat L(\g,\sigma)$ with induced $\varphi$ be a conjugate linear automorphism. Then $\varphi$ is standard by \ref{2.17} and \ref{3.4}, that is of the form $\varphi u(t)=\varphi_t(u(\lambda(t)))$ where $\varphi_t\in\overline{\Aut}\g$ depends smoothly on $t$ and $\lambda:\R\to\R$ is a diffeomorphism with $\lambda(t+2\pi)=\lambda(t)+\epsilon 2\pi$ (actually$\lambda(t)=\epsilon t+t_0$ if $\varphi$ is induced by $\hat\varphi$) for some $\epsilon\in\{\pm 1\}$ and $\varphi_{t+2\pi}=\sigma\varphi_t\sigma^{-\epsilon}$. To avoid confusion with the existing literature we call $\varphi$ or $\hat\varphi$ to be of type $1$ (resp.~type $2$) if $\epsilon=1$ (resp.~$\epsilon=-1$) (in [R...] those of type $1$ are called of $2^{nd}$ kind and those of type $2$ of $1^{st}$ kind). We denote by $\overline{\Aut}_i\GC$ the set of conjugate linear automorphisms of $\GC$ of type $i$ and by $\overline{\Aut}^{2q}_i\GC$ the subset of those of order $2q$. Recall that the order of a conjugate linear automorphism is always even (if finite).

In the following let $\u\subset\g$ be a fixed compact real form of $\g$ and $\omega$ the conjugation with respect to $\u$. Furthermore let for each $r\in\N,\ \tilde{\cal A}^r$ be a fixed set of representatives of the conjugacy classes of $\Aut^r\g\cup\overline{\Aut}^r\g$. Of course $\tilde {\cal A}^r={\cal A}^r\cup\bar{\cal A}^r$ where ${\cal A}^r\subset\Aut^r\g$ and $\bar{\cal A}^r\subset\overline{\Aut}^r\g$ are sets of representatives of conjugacy classes of $\Aut^r\g$ and $\overline {\Aut}^r\g$, respectively.

\begin{definition}
Let $q\in\N$ and $k:=\bar o(\sigma)$.
\begin{enumerate}
\item[(i)] $\bar \J_1^{2q}(\g,k):=\{(p,\varrho_0,[\beta])\mid p\in\{0,1,\dots,2q-1\},\varrho\in\tilde{\cal A}^r\cap\ \omega^{q'}\Aut\g,\beta\in(\omega^l\Aut\g)^{\varrho}$,
$\bar o(\varrho^l\beta^{q'})=k\}$ where $r=r(p,2q),q'=q'(p,2q)$ and $l=l(p,2q)$ are as in \ref{4.3} and $[\beta]$ denotes the conjugacy class of the image of $\beta$ in $\pi_0((\Aut\g\cup{\overline\Aut}\g)^{\varrho})$.
\item[(ii)] $\bar \J_2^{2q}(\g,k):=\{[\varphi_+,\varphi_-]\mid\varphi_\pm\in\overline{\Aut}\g,\ \varphi^2_+=\varphi_-^2,\  \ord (\varphi^2_\pm)=q,\bar o(\varphi_-^{-1}\varphi_+)=k\}$, where $(\varphi_+,\varphi_-)$ and $(\tilde\varphi_+,\tilde\varphi_-)$ like in \ref{5.2} are called equivalent if $\tilde\varphi_+=\alpha\varphi_\pm\alpha^{-1}$ and $\tilde\varphi_-=\beta\varphi_\pm\beta^{-1}$ for some $\alpha,\beta\in\Aut\g$ with $\alpha^{-1}\beta\in((\Aut\g)^{\varphi_+^2})_0$ and $[\varphi_+,\varphi_-]$ denotes the equivalence class.
\end{enumerate}
\end{definition}
Then we can associate to each $\varphi\in\overline{\Aut}^{2q}_iL(\g,\sigma)$ an invariant in $\bar\J_i^{2q}(\g,k)$ in the same manner as in the complex linear case. If $\hat\varphi\in\overline{\Aut}_i^{2q}\hat L(\g,\sigma)$, the invariant of $\hat\varphi$ is by definition that of the induced $\varphi\in\overline{\Aut}^{2q}_i L(\g,\sigma)$. Like in the complex linear case the following result holds (cf.~\ref{5.9}, \ref{4.9}).

\begin{theorem}\label{7.2}
For any $q\in\N$ and $i\in\{1,2\}$, the mapping $\overline{\Aut}^{2q}_i\GC\to\bar\J_i^{2q}(\g,k)$ that associates to each $\hat\varphi$ or $\varphi$ its invariant, induces bijections
\[
\begin{array}{ccc}
\overline{\Aut}^{2q}_1 \GC/\Aut_1\GC&\to&\bar\J_1^{2q}(\g,k)\\
\overline{\Aut}^{2q}_2 \GC/\Aut\GC&\to&\bar\J_2^{2q}(\g,k)
\end{array}
\]
where the quotients denote the sets of conjugacy classes with respect to $\Aut_1 \GC$ and $\Aut\GC$, respectively. Moreover $\overline{\Aut}^{2q}_1\GC/\Aut\GC$  is in bijection with $\bar\J_1^{2q}(\g,k)/{\sim}$ where $\sim$ denotes an equivalence relation which is trivial in case $q=1$. \qed
\end{theorem}

If $\u\subset\g$ is a compact real form of $\g$ invariant under $\sigma$ we call $L(\u,\sigma_\vert)$ or $\hat L(\u,\sigma_\vert)$ a compact real form of $L(\g,\sigma)$, resp., $\hat L(\g,\sigma)$, where $\sigma_\vert$ denotes the restriction of $\sigma$ to $\u$. More generally we make the following definition:

\begin{definition}\label{7.3}
$\U\subset \GC$  is called a \emph{compact real form} of $\GC$ if there exists $\tilde\sigma\in\Aut\g$, a compact $\tilde\sigma$-invariant real form $\u$ of $\g$, and an isomorphism $\Phi:
\tilde\GC \to \GC$ with $\Phi\tilde\U=\U$ where $\tilde\GC=L(\g,\tilde\sigma)$ and $\tilde\U=L(\u,\tilde\sigma_\vert)$ (resp. $\hat L(\g,\tilde\sigma)$ and $\hat L(\u,\tilde\sigma_\vert)$).
\end{definition}

\begin{proposition}\label{7.4}
$\GC$ has a compact real form and this is unique up to conjugation: If $\U_1,\U_2$ are two compact real forms of $\GC$ then there exists $\Phi\in\Aut\GC$ with $\Phi\U_1=\U_2$.
\end{proposition}

\begin{proof}
For simplicity we restrict to the case $\GC =L(\g,\sigma)$. Since we may assume $\sigma$ to be of finite order, $\g$ has a compact $\sigma$-invariant real form~$\u$ (cf.~(B) of Appendix~B) and thus $\GC$ a compact real form $\U_1=L(\u,\sigma_\vert)$. If $\U_2=\Psi L(\tilde\u,\tilde\sigma_\vert)$ is a second, where $\tilde\sigma\in\Aut\g$,
$\Psi:L(\g,\tilde\sigma)\to L(\g,\sigma)$ is an isomorphism, and $\tilde\u$ is a $\tilde\sigma$-invariant compact real form of $\g$ we may assume $\tilde\u=\u$ after eventually conjugating $\tilde\u$ and $\tilde\sigma$.

Since $\GC$ and $\tilde\GC$  are isomorphic, $\bar o(\sigma)=\bar o(\tilde\sigma)$, and hence also $\bar o(\sigma_\vert)=\bar o(\tilde\sigma_\vert)$ as $\Aut\u/\Int\u$ is naturally isomorphic with $\Aut\g/\Int\g$ (cf. Appendix B). Thus there exists an isomorphism from $L(\u,\sigma_\vert)$ to $L(\u,\tilde\sigma_\vert)$ and this extends to an isomorphism $\chi:\GC$
$\to L(\g,\tilde\sigma)$. Therefore $\Phi:=\Psi\circ\chi:\GC\to\GC$ is an isomorphism with $\Phi(\U)=\U_1$.
\end{proof}

In the following we fix one compact real form $\U$ of $\GC$ and denote by $\Omega\in\overline{\Aut}\GC$ the conjugation in $\GC$ with respect to $\U$. The complex linear and conjugate linear extensions of automorphisms $\Phi\in\Aut\U$ to $\Phi_\C\in\Aut\GC$ and $\Phi_\C\Omega\in\overline{\Aut}\GC$, respectively, induce mappings between the conjugacy classes. Similar as in finite dimensions (cf. Appendix B) we have the following result.

\begin{theorem}\label{7.5}
Complex linear resp. conjugate linear extensions of automorphisms of $\U$ from $\U$ to $\GC$ induce the following bijections between conjugacy classes.
\begin{enumerate}
\item[(i)]
$\Aut^q_i\U/\Aut_1\U\to\Aut_i^q\GC /\Aut_1\GC$\quad ($q\in\N$, if $i=1$, resp., $q$ even if $i=2$)
\item[(ii)]
$\Aut^{2q}_i\U/\Aut_1\U\to\overline{\Aut}^{2q}_i \GC / \Aut_1\GC$\quad ($q$ even, if $i=1$, resp., $q\in\N$, if $i=2$)
\item[(iii)]
$(\Aut_1^q\U\cup\Aut_1^{2q}\U)/\Aut_1\U\to\overline{\Aut}^{2q}_1\GC/\Aut_1 \GC$\quad ($q$ odd)
\end{enumerate}
where $i\in\{1,2\}$.
\end{theorem}

\begin{remark}\label{7.6}\rm
The theorem also holds if $\Aut_1\U$ and $\Aut_1\GC$ in the denominators are replaced by $\Aut\U$ and $\Aut\GC$. This is clear if $i=2$ by~\ref{5.1}, and also if $i=1$ and the automorphisms have order $\le 2$ by~\ref{4.9}. It follows otherwise by a slight extension of the proof below.
\end{remark}

{\sl Proof of Theorem~\ref{7.5}.}
Since conjugacy classes of automorphisms of order $q$ are parametrized by their invariants it is enough to show that the above mappings induce bijections between the corresponding sets of invariants.
\begin{enumerate}
\item[1.] We begin with automorphisms of the first kind $(i=1)$. The induced mappings between the invariants are
\begin{enumerate}
\item[(i)] $(p,\varrho,[\beta])\mapsto (p,\varrho_\C, [\beta_\C])$\ ,
\item[(ii)] $(p,\varrho,[\beta])\mapsto (p,\varrho_\C\omega^{q'},[\beta_\C\omega^l])$\ ,
\item[(iii)] $\J_1^q(\u,k)\ni(\tilde p,\varrho,[\beta])\mapsto (p:=2\tilde p,\varrho_\C\omega^{q'},[\beta_\C\omega^l])$\ ,\\
$\J_1^{2q}(\u,k)\ni (p,\varrho,[\beta])\mapsto(p,\varrho_\C\omega^{q'},[\beta_\C\omega^l])$
\end{enumerate}
where $q'=q'(p,2q)$ and $l=l(p,2q)$ are determined from $p$ and $2q$ by \ref{4.3}. Here we have chosen the representatives of the conjugacy classes of automorphisms of finite order on $\u$ and $\g$ in such a way that
\begin{enumerate}
\item[(a)] ${\cal A}^r(\g)=\{\varphi_\C\mid\varphi\in{\cal A}^r(\u)\}$
\item[(b)] $\bar{\cal A}^{2r}(\g)=\{\varphi_\C\omega\mid\varphi\in{\cal A}^{2r}(\u)\}$ if $r$ is even, and
\item[(c)] $\bar{\cal A}^{2r}(\g)=\{\varphi_\C\omega\mid\varphi\in{\cal A}^r(\u)\}\cup\{\varphi_\C\omega\mid\varphi\in{\cal A}^{2r}(\u)\}$ if $r$ is odd.
\end{enumerate}

This is possible due to Proposition B.1 from the Appendix.

Thus for fixed $p$, the bijectivity $\varrho\leftrightarrow\varrho_\C$ is clear in (i) from (a) and that of $\varrho\leftrightarrow\varrho_\C\omega^{q'}$ in (ii) from (a) if $q'=2q/(p,2q)$ even and from (b) if $q'$ is odd as the order $r=(p,2q)$ of $\varrho$ is divisible by 4 in this case. Similarly the bijectivity $\varrho\leftrightarrow\varrho_\C\omega^{q'}$ follows in (iii) from (a) if $p$ is odd (and hence $q'$ even) and from (c) if $p$ is even (hence $q'$ odd). Note that in the first line of (iii) above the order of $\varrho$, i.e. the greatest common divisor of $\tilde p$ and $q$ is odd, while in the second line it is even if $p$ is even.

The bijectivity of the mappings between the sets of invariants follows from Proposition B.2.

\item[2.] We now consider automorphisms of the second kind $(i=2)$. The induced mappings $\J^{2q}_2(\u,k)\to\J_2^{2q}(\g,k)$ resp. $\bar\J_2^{2q}(\g,k)$ are given by $[\varphi_+,\varphi_-]\mapsto[\varphi_{+_\C},\varphi_{-_\C}]$ resp. $[\varphi_{+_\C}\omega,\varphi_{-_\C}\omega]$. The proof of surjectivity amounts to show that for any $[\psi_+,\psi_-]\in\J_2^{2q}(\g,k)\cup\bar\J_2^{2q}(\g,k)$ there exist $\alpha,\beta\in\Aut\g$ with $\alpha^{-1}\beta\in((\Aut\g)^{\psi_+^2})_0$ such that $\alpha\psi_+\alpha^{-1}$ and $\beta\psi_-\beta^{-1}$ leave $\u$ invariant. After a first conjugation of $\psi_+$ and $\psi_-$ with the same automorphism $(\alpha=\beta)$ we may assume $\psi_-\u=\u$ as $\psi_-$ has finite order (cf. (B) of Appendix B). Thus $\psi_+^2=\psi_-^2$ leaves $\u$ invariant and we are looking for an $\alpha\in((\Aut\g)^{\psi^2_+})_0$ with $\alpha\psi_+\alpha^{-1}(\u)=\u$. Let $\psi_+=\psi_0 e^{\ad X}$ with $\psi_0(\u)=\u$ and $X\in i\u$ (cf. (A) of Appendix B). Then $\psi_+^2 e^{-\ad X}=\psi_0e^{\ad X}\psi_0=\psi_0^2 e^{\ad \psi_0^{-1}X}$ and hence $\psi_+^2=\psi_0^2$ and $\psi_0X=-X$ by the uniqueness of this decomposition (cf. (A) of Appendix B). Thus we may choose $\alpha:=e^{\ad  1/2 X}$ which conjugates $\psi_+$ to $\psi_0$ and is contained in $((\Aut\g)^{\psi^2_+})_0$ as $\psi_+^2X=\psi_0^2X=X$.

To prove injectivity, let $[\varphi_+,\varphi_-],[\tilde\varphi_+,\tilde\varphi_-]\in\J^{2q}_2(\u,k)$ with $[\psi_+,\psi_-]=[\tilde\psi_+,\tilde\psi_-]$ where $\psi_\pm$ and $\tilde\psi_\pm$ are the complex resp. conjugate linear extensions of $\varphi_\pm,\tilde\varphi_\pm$ to $\g$. After eventually interchanging $\varphi_+$ with $\varphi_-$, there hence exist $\alpha_\pm\in\Aut\g$ with $\tilde\psi_\pm=\alpha_\pm\psi_\pm\alpha_\pm^{-1}$ and $\alpha_+^{-1}\alpha_-\in((\Aut\g)^{\psi_+^2})_0$. For each $\alpha\in\Aut\g^\R=\Aut\g\cup\overline{\Aut}\g$ there exist unique $\alpha'\in\Aut\g^\R$ and $X\in i\u$ with $\alpha'(\u)=\u$ and $\alpha=\alpha'e^{\ad X}$ (cf. Appendix B). Let $\alpha_\pm=\alpha'_\pm e^{\ad X_\pm}$ with $\alpha'_\pm(\u)=\u$ and $X_\pm\in i\u$. Then $(\tilde\psi_\pm\alpha'_\pm)\cdot e^{\ad X_\pm}=(\alpha'_\pm\psi_\pm)e^{\ad \psi^{-1}_\pm X_\pm}$. Hence $\tilde\psi_\pm=\alpha'_\pm\psi_\pm\alpha_\pm'^{-1}$, and $\psi_\pm X_\pm=X_\pm$, implying $e^{\ad X_\pm}\in((\Aut\g)^{\psi_+^2})_0$ and  $\alpha_+'^{-1}\alpha'_-\in((\Aut\g)^{\psi_+^2})_0$. In other words we may assume that $\alpha_\pm$ leave $\u$ invariant. If $\psi\in\Aut\g$ leaves $\u$ invariant then it follows from the decomposition of the elements of $\Aut\g^\R$ described above that $(\Aut\g)^\psi=\{\alpha\in(\Aut\g)^\psi\mid\alpha\u=\u\}\cdot\{e^{\ad X}\mid X\in i\u^\psi\}$
 and hence that $\{\alpha\in(\Aut\g)^\psi_0\mid\alpha\u=\u\}=\{\alpha\in(\Aut\g)^\psi\mid\alpha\u=\u\}_0$. Thus the restriction of $\alpha^{-1}_+\alpha_-$ to $\u$ is contained in $((\Aut\u)^{\varphi_+^2})_0$. This finally shows $[\varphi_+,\varphi_-]=[\tilde\varphi_+,\tilde\varphi_-]$ and finishes the proof of the injectivity of $\J^{2q}_2(\u,k)\to\J^{2q}_2(\g,k)$ (resp. $\bar\J_2^{2q}(\g,k))$.
\end{enumerate}
\vspace{-27pt}
\qed

We now restrict to the case of conjugate linear involutions.  Their fixed point sets are the real forms of $\GC$. The real forms of $\GC$ corresponding to conjugate linear involutions of type 1 are called {\sl almost compact} and those corresponding to involutions of type 2 {\sl almost split}.

According to \ref{7.5} (iii) the almost compact real forms are in bijection with $\{\id_\U\}\cup\Aut^2_1\U/\Aut\U$, where  of course $\id_U$ corresponds to the compact real form. The conjugacy classes of almost split real forms are by \ref{7.5} (ii) in bijection with $\Aut^2_2\U/\Aut\U$. Thus we have the following result, completely analogous to the finite dimensional case:

\begin{corollary}\label{7.7}
Let $\U$ be a compact real form of $\GC$. Then the conjugacy classes of noncompact real forms of $\GC$ are in bijection with the conjugacy classes of involutions on $\U$. The correspondence is given by $\U={\cal K}+{\cal P}\mapsto{\cal K}+i\cal P$ where $\cal K$ and $\cal P$ are the $(+1)$- and $(-1)$-eigenspaces of the involution. \qed
\end{corollary}

\begin{corollary}\label{7.8}
Each element of $\bar\J_i^2(\g,k)$ can be represented by a conjugate linear involution of the form $\varphi u(t)=\varphi_0 u(\epsilon t+t_0)$ with $\varphi_0\in\overline{\Aut}\g,$ (and $t_0=0$ if $\epsilon=-1\ ,\ t_0\in\{0,\pi\}$ if $\epsilon=+1)$ on some $L(\g,\sigma)$. 
\end{corollary}

\begin{proof}
The corresponding statement is true for $\J_i^2(\u,k)$ and $\J_1^1(\u,k)$ and thus follows for $\bar\J_i^2(\g,k)$ from \ref{7.5} (ii) and (iii).
\end{proof}

There are two obvious candidates of real forms of $L(\g,\sigma)$ (and similarly for $\hat L(\g,\sigma))$, namely $L(\g^*,\sigma^*)$ where $\g^*$ is a $\sigma$-invariant real form of $\g$ and $\sigma^*$ denotes the restriction of $\sigma$ to $\g^*$, and $L_\pi(\g,\bar\varphi):=\{u:\R\to\g\mid u(t+\pi)=\bar\varphi u(t),\ u\in C^\infty\}\cong L(\g,\bar\varphi)$ where $\bar\varphi\in\overline{\Aut}\g$ satisfies $\bar\varphi^2=\sigma$. Note that in the last case any $u\in L(\g,\sigma)$ can be uniquely decomposed as $u_++u_-$ with $u_\pm(t+\pi)=\pm\bar\varphi u_\pm(t)$ (by taking $u_\pm={1\over2}(u(t)\pm\bar\varphi u(t-\pi))$) and thus as $u_1+i u_2:=u_++u_-$ with $u_1,u_2\in L_\pi(\g,\bar\varphi)$.

\begin{proposition}
Up to quasiconjugation the real forms $L(\g^*,\sigma^*)$  and $L_\pi(\g,\bar\varphi)$ of $L(\g,\sigma)$ are precisely the real forms which correspond to involutions of type 1a and 1b, respectively.
\end{proposition}

Of course the analogous statements hold for $\hat L(\g,\sigma)$.

\begin{proof}
\begin{enumerate}
\item[(i)] $L(\g^*,\sigma^*)$ is the fixed point set of the conjugate linear involution $u(t)\mapsto\omega^*u(t)$ where $\omega^*\in\overline{\Aut}\g$ denotes conjugation with respect to $\g^*$. Conversely, a conjugate linear involution of type 1a is quasiconjugate by \ref{7.8} to one of the form $u(t)\mapsto\omega^* u(t)$ on some $L(\g,\tilde\sigma)$ where $\omega^*\in\overline{\Aut}\g$ is an involution which commutes with $\tilde\sigma$. Its fixed point set is $\{u\in L(\g,\tilde\sigma)\mid u(t)\in\g^*\}=L(\g^*,\sigma^*)$ where $\g^*$ is the real form corresponding to $\omega^*$ and $\sigma^*$ is the restriction of $\tilde\sigma$.
\item[(ii)] $L_\pi(\g,\bar\varphi)$ is the fixed point set of the conjugate linear involution $u(t)\mapsto\bar\varphi^{-1}u(t+\pi)$ of $L(\g,\bar\varphi^2)$. Conversely, a conjugate linear involution of type 1b is quasiconjugate by \ref{7.8} to one of the form $u(t)\mapsto\bar\varphi^{-1}u(t+\pi)$ on some $L(\g,\sigma)$ where $\bar\varphi\in\overline{\Aut}\g$ and $\bar\varphi^{-2}\sigma=\id$, i.e. $\sigma=\bar\varphi^2$.
\end{enumerate}
\vspace{-27pt}
\end{proof}

Hence, given $\g$ and $\sigma\in\Aut\g$ the following objects are in bijective correspondence:
\begin{enumerate}
\item[(i)] Non compact real forms of $\hat L(\g,\sigma)$ (or $L(\g,\sigma))$ of type 1 a up to isomorphism
\item[(ii)] Pairs $(\varrho,[\beta])$ where $\varrho\in{\cal A}^2(\g)$, $\beta\in(\Aut\g)^\varrho$, $\beta$ is conjugate to $\sigma$ in $\pi_0(\Aut\g)$, and $[\beta]$ denotes the conjugacy class of $\beta$ in $\pi_0((\Aut\g)^\varrho)$
\item[(iii)] Affine Kac-Moody algebras $\hat L(\g^*,\sigma^*)$ (or loop algebras $L(\g^*,\sigma^*)$) up to isomorphism where $\g^*$ is a non compact real form of $\g$, $\sigma^*\in\Aut\g^*$, and $\sigma^*_\C$ is conjugate to $\sigma$ in $\pi_0(\Aut\g)$.
\end{enumerate}
Note that the bijection between (ii) and (iii) also follows from Corollaries~\ref{3.5} and \ref{2.13}, respectively, since $\varrho\in{\cal A}^2(\g)$ corresponds to an isomorphism class of a non compact real form $\g^*$ of $\g$ and $\hat L(\g^*,\sigma^*)$ is isomorphic to $\hat L(\g^*,\tilde\sigma^*)$ if and only if $\sigma^*$ and $\tilde\sigma^*$ are conjugate in $\pi_0(\Aut\g^*)$ which is isomorphic to $\pi((\Aut\g)^\varrho)$ (cf.~Proposition B.2 (i) of the Appendix).

The almost split real forms of $L(\g,\sigma)$ are in bijection with $\J_2^2(\g,k)$ and thus with $\J_2^2(\u,k)$ where $\u$ is a compact real form of $\g$ and $k=\bar o(\sigma)$. If $[\varrho_+,\varrho_-]\in\J_2^2(\u,k)$ then the corresponding real form of $L(\g,\varrho_-\varrho_+)\cong L(\g,\sigma)$ is
\begin{eqnarray*}
\GC^*&=&\{u\in L(\u,\varrho_-\varrho_+)\mid\varrho_+u(t)=u(-t)\}
\oplus i\{ u\in L(\u,\varrho_-\varrho_+)\mid\varrho_+u(t)=-u(-t)\}\\
&=&\{u:\R\to\u\mid\varrho_+u(t)=u(-t),\varrho_-u(t+\pi)=u(-t+\pi),u\in C^\infty\}\oplus\\
&&\oplus\ i\{ u:\R\to\u\mid\varrho_+u(t)=-u(-t),\varrho_-u(t+\pi)=-u(-t+\pi),u\in C^\infty\}\ .
\end{eqnarray*}
\begin{remark}\label{7.10}\rm
If $\varrho_-\varrho_+$ has finite order, which we may assume by \ref{8.14}, and $\ord\varrho_-\varrho_+=l$ then $\GC^*$ may be described by developing each $u$ into its Fourier series also as follows.
\begin{eqnarray*}
\GC^*&=&\{{\sum\limits_{n\in\Z}}u_ne^{int/l}\mid u_n\in\g,\varrho_-\varrho_+u_n=e^{2\pi in/l}u_n,\varrho_+u_n=\omega u_n\}\\
&=&\{{\sum\limits_{n\in\Z}}u_ne^{int/l}\mid u_n\in\g^*_+,u_ne^{n\pi i/l}\in\g_-^*\}
\end{eqnarray*}
where the sums are supposed to represent $C^\infty$-functions, $\omega:\g\to\g$ denotes complex conjugation with respect to $\u$, and $\g_\pm^*:=\g^{\varrho_\pm\omega}$. In particular the finite sums ${\sum\limits_{\vert n\vert\le N}}u_ne^{int/l}\in\GC^*$, which can be viewed by replacing $e^{int/l}$ by $z$ as functions on $\C$, take values in $\g_+^*$ and $\g_-^*$ when restricted to the real line and the line $\R e^{\pi i/l}$, respectively.
\end{remark}

As an example let $\g=\s\l(2,\C)$. Since $\g$ has no outer automorphisms all $\hat L(\g,\sigma)$, resp., $L(\g,\sigma)$ are isomorphic to $\hat L(\g)$, resp., $L(\g)$. Up to conjugation $\g$ has only one involution, and this may be represented  by $\tau:=\Ad 
\left(
\begin{array}{cc}
    1&\\

    &-1
    \end{array}
    \right)$
or $\mu$ with $\mu A=-A^t$. A compact real form of $\g$ is $\u:=\su(2)$ and this is invariant under $\mu$ and $\tau$. Let $\omega$ be the conjugation with respect to $\u$, i.e. $\omega A=-\bar A^t$. Up to conjugation, the only noncompact real form of $\g$ is $\g^{\mu\omega}=\s\l(2,\R)$. Moreover $(\Aut\ \s\l(2,\C))^\mu$ has two connected components, represented e.g. by $\id$ and $\tau$. Thus the almost compact real forms of $L(\s\l(2,\C))$ are up to quasconjugation
\begin{enumerate}
\item[1 a)] $L(\su(2)), L(\s\l(2,\R)), L(\s\l(2,\R),\tau)$
\item[1 b)] $L_\pi(\s\l(2,\C),\omega)$.
\end{enumerate}

There are three almost split real forms of $L(\s\l(2,\C))$ corresponding to $\{[\id,\id],[\mu,\mu],[\mu,\id]\}=J^2_2(\s\l(2,\C))$, and these may be described as\\
$\{\Sigma u_n e^{int}\mid u_n\in\su(2)\}$\\
$\{\Sigma u_n e^{int}\mid u_n\in\s\l(2,\R)\}$\\
$\{\Sigma u_ne^{int/2}\mid u_n\in\s\l(2,\R)$ and $(i)^nu_n\in\su(2)\}$\\
where all sums are assumed to represent $C^\infty$-functions.

We finally consider Cartan decompositions. Let $\GC:=\hat L(\g,\sigma)$  or $L(\g,\sigma)$ be as above where $\g$ is complex (and simple) and $\sigma\in\Aut\g$.

\begin{definition}\label{7.11}
Let $\GC^*$ be a real form of $\GC$. Then a vector space decomposition $\GC^*=\KC+\PC$ is called
a Cartan decomposition of $\GC^*$ if there exists a compact real
form $\U$ of $\GC$ such that $\KC=\GC^*\cap\U$ and $\PC=\GC^*\cap
i\U$.
\end{definition}

As in finite dimensions it follows that $\KC$ and $\PC$ are the
$(+1)$- and $(-1)$-eigenspaces of an involution on $\GC^*$ and that
conversely such an eigenspace decomposition $\GC^*=\KC+\PC$ with
respect to some involution is a Cartan decomposition if and only if
$\U:=\KC+i\PC$ is a compact subalgebra of $\GC$.

\begin{corollary}\label{7.12}
Any real form $\GC^*$ of $\GC$ has a Cartan decomposition and this
is unique up to conjugation.
\end{corollary}

\begin{proof}
{\sl Existence} follows from the fact that noncompact real forms
$\GC^*$ correspond to involutions on a compact real form $\U$ and thus
arise as $\KC+i\tilde\PC$ where $\U=\KC+\tilde\PC$ is the eigenspace
decomposition of $\U$ with respect to some involution. Hence
$\GC^*=\KC+i\tilde \PC$ is a Cartan decomposition. If $\GC^*$ is compact
we take $\KC:=\GC^*$ and $\PC:=\{0\}$.

To prove {\sl uniqueness} let $\GC^*=\KC_1+\PC_1=\KC_2+\PC_2$ be two
Cartan decompositions. Then $\U_1:=\KC_1+i\PC_1$ and
$\U_2:=\KC_2+i\PC_2$ are two compact real forms of $\GC$ and hence
isomorphic by~\ref{7.4}. Let $\alpha:\U_1\to\U_2$ be an isomorphism and $\KC_1':=\alpha(\KC_1)$, $\PC_1':=\alpha_\C(\PC_1)$. The two decompositions of $\U_2$ into $\KC_2+i\PC_2$ and $\KC_1'+i\PC_1'$ yield isomorphic real forms $\GC^*=\KC_2+\PC_2$ and $\alpha_\C\GC^*=\KC_1'+\PC_1'$. Hence injectivity of the mapping in~\ref{7.5} (iii) yields an isomorphism $\beta$ of $\U_2$ with $\beta(\KC_1')=\KC_2$ and $\beta(i\PC_1')=i\PC_2$. Thus $\beta\circ\alpha:\U_1\to\U_2$ maps $\KC_1$ to $\KC_2$ and $i\PC_1$ to $i\PC_2$. Its complex linear extension therefore preserves $\GC^*$ and conjugates the two Cartan decompositions.
%\vspace{-27pt}
\end{proof}

\section{The algebraic case}

Involutions, finite order automorphisms and real forms of affine Kac-Moody algebras have been studied in the algebraic setting (cf. below) by numerous authors, starting with a paper by F. Levstein \cite{Lev} and culminating in the classification of involutions and real forms in \cite{B$_3$R} and \cite{BMR}. Our aim here is to show that our elementary methods also work in this case and lead to the same results as in the $C^\infty$-case if suitably modified and combined with a basic result of Levstein. In particular conjugacy classes of automorphisms of finite order as well as real forms are also classified in this situation by the invariants introduced above and are thus in bijection with their smooth counterparts.

\subsection{Preliminaries}

Let $\g$ be an arbitrary simple Lie algebra over $\F=\R$ or $\C$. Only from \ref{8.4} onwards we will restrict $\g$ to be compact if $\F=\R$. Let $\sigma\in\Aut\g$ be an automorphism of finite order with $\sigma^l=\id$ ($l$ not necessarily the order of $\sigma$). Then we call
\[
L_\alg(\g,\sigma):=\{u:\R\to\g\mid u(t+2\pi)=\sigma u(t),\ u(t)={\sum\limits_{\vert n\vert\le N}}u_ne^{{int}/l},\ 
N\in\N,\ u_n\in\g_\C\}
\]
the algebraic loop algebra where $\g_\C$ denotes the complexification of $\g$ if $\F=\R$ and $\g_\C=\g$ if $\F=\C$. Note that for $\F=\R$ the sum is contained in $\g$ if and only if $u_{-n}=\overline{u_n}$ where the bar denotes conjugation with respect to $\g$. The algebraic loop algebra is a Lie algebra with the pointwise bracket $[u,v]_0(t):=[u(t),v(t)]$.

\begin{remark}\label{8.1}\rm
Usually one embeds $L_\alg(\g,\sigma)$ into $L_\alg(\g):=L_\alg(\g,\id)$ as $\{u:\R\to\g\mid u(t+2\pi/l)=\sigma u(t), u(t)={\sum\limits_{\vert n\vert\le N}}u_ne^{int}\in\g, N\in\N, u_n\in\g_\C\}$ by $u(t)\mapsto u(l\cdot t)$. But for our purposes the above definition is more convenient. Note also that the embedding depends on $l$ while $L_\alg(\g,\sigma)$ is independent of $l$. In fact
\[
L_\alg(\g,\sigma)=\{u:\R\to\g\mid u(t+2\pi)=\sigma u(t),\ u(t)={\sum\limits_{\vert n\vert\le N}}u_ne^{iq_nt},\ N\in\N,\ u_n\in\g_\C,\ q_n\in\Q\}
\]
as $u(t+2\pi l)=u(t)$ implies $q_n\in{1\over l}\Z$.
\end{remark}

\begin{remark}\label{8.2}\rm
The last description of $L_\alg(\g,\sigma)$ would make sense also for a $\sigma\in\Aut\g$ which is not necessarily of finite order. But since the sums are finite there would exist for each $u$ an $l=l(u)$ with $u(t+2\pi l)=u(t)$ implying $u(t)\in\s:=\{x\in \g\mid \sigma^kx=x$ for some $k\in\N\}$. Thus $\g$ could be replaced by $\s$ which is a $\sigma$-invariant subalgebra on which $\sigma$ has finite order.
\end{remark}

\begin{definition}\label{8.3}
If $\sigma\in\Aut\g$ and $l\in\N$ with $\sigma^l=\id$ we let
\[
\g_n:=\g_n(\sigma,l):=\{x\in\g_\C\mid\sigma x=e^{2\pi in/l}x\}
\]
for all $n\in\Z$.
\end{definition}

Of course, $\g_{n+l}=\g_n$ and $\g=\g_0\oplus\dots\oplus\g_{l-1}$. Moreover, $u(t)=\sum u_ne^{int/l}$ satisfies $u(t+2\pi)=\sigma u(t)$ if and only if $u_n\in\g_n$ for all $n$. Since $\sigma$ leaves the Killing form $(.,.)_0$ of $\g$ invariant, $(\g_n,\g_m)_0=0$ unless $l$ divides $n+m$.

We now extend $L_\alg(\g,\sigma)$ to the affine Kac-Moody algebra $\hat L_\alg(\g,\sigma)$ in the usual way by
\[
\hat L_\alg(\g,\sigma)=L_\alg(\g,\sigma)+\F c+\F d
\]
with
\[
\begin{array}{ll}
&[c,x]:=[x,c]=0\\
&[d,u]:=-[u,d]:=u'\\
&[u,v]:=[u,v]_0+(u',v)c
\end{array}\]
for all $x\in\hat L_\alg(\g,\sigma)$ and $u,v\in L_\alg(\g,\sigma)$, where $u'$ denotes the derivative of $u$ and $(u,v)={1\over2\pi}{\int\limits^{2\pi}_0}(u(t),v(t))_0dt$ is the averaged Killing form. Then $\hat L_\alg(\g,\sigma)$ is a Lie algebra and the {following} well known result is easily proved.

\begin{proposition}\label{8.4}\
\begin{enumerate}
\item[(i)] The derived algebra $\hat L'_\alg(\g,\sigma)$ of $\hat L_\alg(\g,\sigma)$ is equal to $L_\alg(\g,\sigma)\oplus\F c$
\item[(ii)] $\F c$ is the center of $\hat L_\alg(\g,\sigma)$ and $\hat L'_\alg(\g,\sigma)$
\item[(iii)] $L_\alg(\g,\sigma)$ is isomorphic to $\hat L'_\alg(\g,\sigma)/\F c$
\item[(iv)] $L_\alg(\g,\sigma)$ is equal to its derived algebra.
\qed
\end{enumerate}
\end{proposition}
$\hat L_\alg(\g,\sigma)$ carries the natural symmetric biinvariant form $(.,.)$ which extends the biinvariant form on $L_\alg(\g,\sigma)$ by the requirements $c,d\ \bot\ L_\alg(\g,\sigma)$, $(c,c)=(d,d)=0$, and $(c,d)=1$.

\subsection{Isomorphisms between loop algebras}
Let $\tilde\g$ be another simple Lie algebra over $\F$ and $\tilde\sigma\in\Aut\tilde\g$ with $\tilde\sigma^{\tilde l}=\id$. Then we have the following examples of isomorphisms:
\begin{enumerate}
\item[(i)] $\tau_r:L_\alg(\g,\sigma)\to L_\alg(\g,\sigma)$ with
\[
\tau_r({\sum\limits_{\vert n\vert\le N}}u_ne^{int/l})={\sum\limits_{\vert n\vert\le N}} u_nr^{n/l}e^{int/l}
\]
where $r>0$ and $\F=\C$. Note that this example does not exist if $\F=\R$ (unless $r=1$, i.e.~$\tau_r=\id)$ since $u_{-n}=\overline{u_n}$ in that case. Note also that the definition of $\tau_r$ does not depend on $l$ (with $\sigma^l=\id)$.
\item[(ii)] $\varphi: L_\alg(\g,\sigma)\to L_\alg(\tilde\g,\tilde\sigma)$ \textit{standard}, i.e.~
\[
\varphi u(t)=\varphi_t(u(\epsilon t+t_0))
\]
where $\epsilon\in\{\pm 1\},\ t_0\in\R,\ \varphi_t:\g\to\tilde\g$ is an isomorphism for all $t\in\R$ such that $\varphi_t$ depends ``algebraically'' on $t$, that is $\varphi_t={\sum\limits_{\vert n\vert\le N}}\varphi_ne^{int/L}$ for some $L,N\in\N$ and $\varphi_n\in\Hom(\g,\tilde\g)$, and $\varphi_{t+2\pi}=\tilde\sigma\varphi_t\sigma^{-\epsilon}$ (the ``periodicity condition'') holds. The periodicity condition is equivalent to $\varphi u\in L_\alg(\tilde\g,\tilde\sigma)$. The $L$ in the description of $\varphi_t$ can be chosen to be the smallest common multiple of $l$ and $\tilde l$ (if $\sigma^l=\tilde\sigma^{\tilde l}=\id$) and in particular to be $l$ if $\sigma=\tilde\sigma$.

Note that with $\varphi_t$ also $\varphi_t^{-1}$ depends algebraically on $t$. For the $\varphi_t$ may be viewed as invertible matrices whose entries are finite Laurent series and whose determinants are constant as the $\varphi_t$ preserve the nondegenerate Killing forms of $\g$ and $\tilde\g$. The inverse of $u(t)\mapsto\varphi_t(u(\epsilon t+t_0))$ is thus also standard and given by
\[
u(t)\mapsto\varphi^{-1}_{\epsilon t-\epsilon t_0}(u(\epsilon t-\epsilon t_0))\ .
\]
\end{enumerate}

\begin{theorem}\label{8.5}
Any isomorphism
\[
\psi: L_\alg(\g,\sigma)\to L_\alg(\tilde\g,\tilde\sigma)
\] is of the form $\psi=\varphi\circ\tau_r$ where $\varphi$ is standard and $r>0$, and this decomposition is unique. Moreover, $\psi$ is standard (i.e.~$r=1$) if $\F=\R$.
\end{theorem}

\begin{proof}
We first assume $\F=\C$.
\begin{enumerate}
\item[(i)] (Uniqueness) Let $\varphi$ be standard, $r>0$ and $\varphi=\tau_r$. From $u:=u_ne^{int}$ we obtain $\varphi_t(u_n)e^{in(\epsilon t+t_0)}=u_nr^{n/l} e^{int/l}$ for all $u_n\in\g_n$. Since $\g_n=\g_{n+l}$ we may replace $n$ by $n+l$ and get $\vert r\vert=1$ and hence $r=1$ and thus $\tau_r=\id=\varphi$. Since $\tau_r\circ\tau_s=\tau_{r+s}$ and the composition and inverses of standard isomorphisms are standard, uniqueness follows.
\item[(ii)] (Existence) Let $\psi:L_\alg(\g,\sigma)\to L_\alg(\tilde\g,\tilde\sigma)$ be an isomorphism. 
Let $C^\alg_\per(\R,\C):=\{{\sum\limits_{\vert n\vert\le N}}a_ne^{int}\mid N\in\N,\ a_n\in\C\}$.
The essential point is to prove for each $t\in\R$ the existence  of a mapping $\alpha_{t}:C^\alg_\per(\R,\C)\to\C$ such that
\[
\psi(fu)(t)=\alpha_{t}(f)\cdot\psi(u)(t)
\]
for all $f\in C^\alg_\per(\R,\C)$ and $u\in L_\alg(\g,\sigma)$. This follows exactly as in the smooth case, cf.~the proof of \ref{2.9}. From this we get 
\[
\psi(fu)=\alpha(f)\cdot\psi(u)
\]
for all $f\in C^\alg_\per(\R,\C)$ and $u\in L_\alg(\g,\sigma)$ where $\alpha:C^\alg_\per(\R,\C)\to C^\alg_\per(\R,\C)$ is defined by $\alpha (f)(t):=\alpha_t(f)$ since we can choose for example $\psi(u)$ to be a constant different from zero to see that $\alpha(f)$ is periodic and algebraic. The mapping $\alpha$ is necessarily an algebra homomorphism with $\alpha(1)=1$. From $\alpha(e^{it})\cdot\alpha(e^{-{it}})=1$ we conclude $\alpha(e^{it})=be^{i\epsilon t}$ for some $b\in \C^*$ and $\epsilon\in \Z$. Since $\psi$ and hence $\alpha$ are isomorphisms we actually get $\epsilon\in\{\pm 1\}$. Let $a$ be a complex number with $a^l=b$ and $\varphi_t:\g\to\tilde\g$ be defined by
\[
\psi(u_ke^{ikt/l})=\varphi_t(u_k)a^k e^{i\epsilon kt/l}
\]
where $u_k\in\g_k$ and $k\in\{0,\dots,l-1\}$. If $k\in\Z$ is arbitrary and $u_k\in\g_k$ we write $k$ as $k=\tilde k+ml$ with $\tilde k\in\{0,\dots,k-1\}$ and $m\in\Z$ and obtain
\[
\begin{array}{rclcl}
\psi(u_ke^{ikt/l})&=&\psi(u_ke^{i\tilde k t/l}\cdot e^{imt})&=&(b e^{i\epsilon t})^m\cdot\varphi_t(u_k)a^{\tilde k}e^{i\epsilon\tilde k t/l}\\
&=&\varphi_t(u_k)a^ke^{i\epsilon kt/l}
\end{array}
\]
as well. Let $a=r^{1/l}e^{{it_0}/l}$ for some $r>0$ and $t_0\in\R$. Then $\psi(u)=\varphi_t(\tau_r(u(\epsilon t+t_0))$ for all $u\in L_\alg(\g,\sigma)$ and the theorem follows.

If $\F=\R$ we decompose the complexification $\psi_\C$ of $\psi$ as $\psi_\C=\varphi\circ\tau_r$ by means of (i). Since $\psi_\C$ commutes with the conjugation we have $\psi_\C(\overline{u_k}e^{-ikt/l})=\overline{\psi_\C(u_ke^{ikt/l}})$ for any $u_k\in\g_k$ and hence $r^{-k/l}\varphi(\overline{u_k}e^{-ikt/l})=r^{k/l}\overline{\varphi(u_ke^{ik/l}})$. This implies $\varphi_t(\overline{u_k})=r^{2k/l}\overline{\varphi_t(u_k})$ and thus $r=1$, since we may replace $k$ by $k+l$ in the last equation without changing $u_k$.
\end{enumerate}
\vspace{-27pt}
\end{proof}

A direct calculation gives:
\begin{lemma}\label{8.6}
Let $\F=\C$ and $\varphi:L_\alg(\g,\sigma)\to L_\alg(\g,\sigma)$ be a standard automorphism with $\varphi u(t)=\varphi_t(u(\epsilon t+t_0))$ and $\varphi_t={\sum\limits_{\vert n\vert\le N}}\varphi_me^{imt/L}$ for some $L, N\in\N$. Let $r>0$. Then
\[
\tau_r\circ\varphi=\ {^r\varphi}\circ\tau_{r^\epsilon}
\]
where $^r\varphi$ is the standard automorphism with the same $\epsilon$ and $t_0$ but with $^r\varphi_t=\Sigma\varphi_m r^{m/L}e^{imt/L}$.
\end{lemma}

We call an isomorphism $\psi=\varphi\circ\tau_r:L_\alg(\g,\sigma)\to L_\alg(\tilde\g,\tilde\sigma)$ with $\varphi u(t)=\varphi_t(u(\epsilon t+t_0))$ be of the first (second) kind if $\epsilon=1$ (resp. $\epsilon=-1$). By the last lemma and our previous results about standard isomorphisms~(\ref{3.10}) this notion is invariant under conjugation.

\begin{corollary}\label{8.7}
\begin{enumerate} 
\item[(i)] Any automorphism  $\varphi$ of $L_\alg(\g,\sigma)$ of finite order and of the first kind is standard. Moreover $\tau_r\varphi\tau_r^{-1}={^r\varphi}$.
\item[(ii)] For each automorphism $\varphi$ of the second kind there exists a unique $r>0$ such that $\tau_r\varphi\tau_r^{-1}$ is standard.
\end{enumerate}
\end{corollary}

\begin{proof}
\begin{enumerate}
\item[(i)] Let $\varphi=\tilde\varphi\circ\tau_r$ be of order $k$ with $\tilde\varphi$ standard and $r>0$. Then $\id=\varphi^k=\psi\circ\tau_{r^k}$ where $\psi$ is standard. Hence $r^k=1$ and thus $r=1$. The second statement follows from~\ref{8.6}.
\item[(ii)] Let $\varphi=\tilde\varphi\circ\tau_s$ with $\tilde\varphi$ standard and $s>0$. Then $\tau_r\tilde\varphi \tau_r^{-1}={^r}\tilde\varphi\circ\tau_{r^2\cdot s}$ by~\ref{8.6} and the claim follows.
\end{enumerate}
\vspace{-27pt}
\end{proof}

\subsection{Isomorphisms between affine Kac-Moody algebras}
We show in this subsection that  there is a bijection between the sets of automorphisms of finite order of $\hat L_\alg(\g,\sigma)$ and $L_\alg(\g,\sigma)$ also in the algebraic case. This is known (and follows from \cite{PK}, Theorem~2 and Corollary~11 and \cite{Rou2}, Proposition~\ref{2.5}). We include an elementary proof also to make the paper self contained as much as possible.

Any isomorphism $\hat\varphi:\hat L_\alg(\g,\sigma)\to\hat L_\alg(\tilde\g,\tilde\sigma)$ induces an isomorphism $\check\varphi:\hat L'_\alg(\g,\sigma)\to\hat L'_\alg(\tilde\g,\tilde\sigma)$ and any isomorphism $\check\varphi$ between the derived algebras induces an isomorphism $\varphi$ between the loop algebras by ignoring elements in the center. In particular we have natural homomorphisms $\Aut\hat L_\alg(\g,\sigma)\to\Aut\hat L'_\alg(\g,\sigma)$ and $\Aut\hat L'_\alg(\g,\sigma)\to\Aut L_\alg(\g,\sigma)$ which we simply denote by $\hat\varphi\mapsto\check\varphi$ and $\check\varphi\mapsto\varphi$, respectively.

\begin{theorem}\label{8.8}
\begin{enumerate}
\item[(i)] Any isomorphism $\varphi:L_\alg(\g,\sigma)\to L_\alg(\tilde\g,\tilde\sigma)$ extends to an isomorphism $\hat\varphi:\hat L_\alg(\g,\sigma)\to\hat L_\alg(\tilde\g,\tilde\sigma)$, that is $\hat\varphi\mapsto\varphi$ is surjective.
\item[(ii)] The mapping $\Aut\hat L'_\alg(\g,\sigma)\to\Aut L_\alg(\g,\sigma)$ given by $\check\varphi\mapsto\varphi$ is an isomorphism.
\item[(iii)] The kernel of the mapping $\Aut\hat L_\alg(\g,\sigma)\to\Aut\hat L'_\alg(\g,\sigma)$ given by $\hat\varphi\mapsto\check\varphi$ is
\[\{\hat\varphi\in\Aut\hat L_\alg(\g,\sigma)\mid\exists\nu\in\F:\ \hat\varphi c=c,\ \hat\varphi d=d+\nu c,\ \hat\varphi u=u,\mbox{ for all } u\in L_\alg(\g,\sigma)\}
\]
and thus isomorphic to $(\F,+)$.
\item[(iv)] $\Aut\hat L_\alg(\g,\sigma)\cong\Aut_{(,)}\hat L_\alg(\g,\sigma)\times\F$ where $\Aut_{(,)}$ denotes the subgroup of automorphisms preserving the natural bilinear form $(.,.)$ of $\hat L_\alg(\g,\sigma)$.
\end{enumerate}
\end{theorem}

To prepare the proof we note that any isomorphism $\hat\varphi:\hat L_\alg(\g,\sigma)\to\hat L_\alg(\tilde\g,\tilde\sigma)$ is of the form
\begin{eqnarray}
\hat\varphi c&=&\lambda c\\
\hat\varphi d&=&\mu d+x_\varphi+\nu c\\
\hat \varphi u&=&\varphi u+\alpha(u)c
\end{eqnarray}
since it preserves the center and the derived algebra. Here $\lambda,\mu\in\F\setminus\{0\}$, $\nu\in\F$, $\alpha:L_\alg(\g,\sigma)\to\F$ is linear, $x_\varphi\in L_\alg(\tilde\g,\tilde\sigma)$, and $\varphi$ is the induced isomorphism between the loop algebras. Conversely, given $\lambda,\mu,\nu,\alpha, x_\varphi$ and $\varphi$ as above, the linear mapping $\hat\varphi$ described by (8.1) - (8.3) is an isomorphism if and only if
\begin{eqnarray}
\varphi(u)'&=&-{1\over\mu}[x_\varphi,\varphi u]_0+{1\over\mu}\varphi(u')\\
\alpha(u')&=&-{1\over\mu}(x_\varphi,\varphi u')\\
((\varphi u)',\varphi v)&=&\lambda (u',v)+\alpha([u,v]_0)
\end{eqnarray}
for all $u,v\in L_\alg(\g,\sigma)$.

By means of (8.4), (8.6) is equivalent to
\[
(\varphi(u'),\varphi(v))-\lambda\mu(u',v)=\mu\alpha([u,v]_0)+(x_\varphi,\varphi[u,v]_0) \ .
\]
Furthermore $(\hat\varphi x,\hat\varphi y)=(x,y)$ for all $x\in\hat L'_\alg(\g,\sigma)$ and $y\in\hat L_\alg(\g,\sigma)$ if and only if $\lambda\mu=1$, $(\varphi u,\varphi v)=(u,v)$, and $\alpha(u)=-{1\over\mu}(x_\varphi,\varphi u)$ for all $u,v \in L_\alg(\g,\sigma)$. Moreover $\hat\varphi$ preserves the bilinear form if and only if $\nu=-{1\over2\mu}(x_\varphi,x_\varphi)$ in addition .

\begin{lemma}\label{8.9}
Let $\F=\C$, $\alpha:L_\alg(\g,\sigma)\to\C$ be linear and $\beta\in\C$. If $\beta(u',v)=\alpha([u,v]_0)$ for all $u,v\in L_\alg(\g,\sigma)$ then $\alpha=0$ and $\beta=0$.
\end{lemma}

\begin{proof}
Let $\sigma^l=\id$. The equation yields for $u:=u_me^{imt/l}$ and $v:=v_ne^{int/l}$ with $u_m\in\g_m$ and $v_n\in\g_n$: $\beta(u_m,v_n)_0{im\over l}{1\over2\pi}{\int\limits^{2\pi}_0}e^{i(m+n)t/l}=\alpha([u_m,v_n]e^{i(m+n)t/l})$ and in particular
\[
\beta(u_m,v_{-m})_0{im\over l}=\alpha([u_m,v_{-m}])\ .
\]
If we replace $m$ by $m+l$ without changing $u_m$ and $v_{-m}$ and take the difference we get ${\beta(u_m,v_{-m})_0=0}$ and thus $\beta=0$ which finally implies $\alpha=0$ by~\ref{8.4}~(iv).
\end{proof}

{\sl Proof of Theorem~\ref{8.8}.}
We may assume $\F=\C$. The real case follows from this essentially by complexification.
\begin{enumerate}
\item[(i)] It is enough to consider the cases $\varphi=\tau_r$ and $\varphi$ standard. Now $\tau_r$ can be extended to a $\hat\tau_r$ as described in (8.1) - (8.3) with $\lambda=\mu=1$, $\nu=0$, $x_\varphi=0$, and $\alpha=0$. Since $\tau_r(u)'=\tau_r(u')$ and $\tau_r$ preserves the natural bilinear form, the equations (8.4) - (8.6) are satisfied. Actually the extension is equal to $e^{s \ad d}$, where $s:=-i\ \log r$. It also preserves the natural bilinear form.

If $\varphi$ is standard with $\varphi u(t)=\varphi_t(u(\epsilon t+t_0))$ then we define $\hat\varphi$ by $\lambda:=\mu:=\epsilon$ and $\alpha(u):=-\epsilon(x_\varphi,\varphi u)$ where $x_\varphi$ is determined by $\varphi'_t\varphi_t^{-1}=-\epsilon\,\ad x_\varphi$. The conditions (8.4) - (8.6) are obviously satisfied (for any $\nu$). If in addition we choose $\nu:=-{\epsilon\over2}(x_\varphi,x_\varphi)$ then $\hat\varphi$ also preserves the natural bilinear form.
\item[(ii)] Surjectivity follows from (i) and injectivity from Lemma~\ref{8.9}.
In fact, any $\check\varphi$ from the kernel is of the form
\[
\begin{array}{lll}
\check\varphi c&=&\lambda c\\
\check\varphi u&=&u+\alpha(u)c
\end{array}
\]
and therefore satisfies $(1-\lambda)(u',v)=\alpha([u,v]_0)$ for all $u,v\in L_\alg(\g,\sigma)$ as follows from $[\check\varphi u,\check\varphi v]=\check\varphi[u,v]$.
\item[(iii)] If $\hat\varphi\in\Aut\hat L_\alg(\g,\sigma)$ restricts to the identity on $\hat L'(\g,\sigma)$ then $\hat\varphi c=c,\hat\varphi d=\mu d+x_\varphi+\nu c$ and $\hat\varphi u=u$ and hence
\[
\begin{array}{lll}
(1-\mu)u'&=&[x_\varphi, u]_0\ ,\\
(x_\varphi,u')&=&0
\end{array}
\]
for all $u\in L_\alg(\g,\sigma)$. This implies $\mu=1$ and $x_\varphi=0$ by Lemma~\ref{8.9}.
\item[(iv)] Since the above constructed extensions of automorphisms from $L_\alg(\g,\sigma)$ to $\hat L_\alg(\g,\sigma)$ as well as the elements of the kernel of $\hat\varphi\to\varphi$ satisfy $\lambda=\mu$ and $(\hat\varphi x,\hat\varphi y)=(x,y)$ for all $x\in\hat L'_\alg(\g,\sigma)$ and $y\in\hat L_\alg(\g,\sigma)$ (equivalently $\lambda=\mu\in\{\pm 1\}$, $(\varphi u,\varphi v)=(u,v)$, $\alpha(u)=-{1\over\mu}(x_\varphi,\varphi u)$ for all $u,v\in L_\alg(\g,\sigma))$ this is true for all $\hat\varphi\in\Aut\hat L_\alg(\g,\sigma)$. In particular the kernel of $\hat\varphi\mapsto\check\varphi$ lies in the center of $\Aut\hat L_\alg(\g,\sigma)$ and $\hat\varphi\in\Aut\hat L_\alg(\g,\sigma)$ preserves the natural bilinear form if and only if $\nu=-\epsilon (x_\varphi,x_\varphi)$ in (\ref{8.2}). Hence the product decomposition of $\Aut\hat L_\alg(\g,\sigma)$ follows readily. \qed
\end{enumerate}

By the last part of the proof, a linear mapping $\hat\varphi:\hat L_\alg(\g,\sigma)\to\hat L_\alg(\tilde\g,\tilde\sigma)$ as given by (\ref{8.1}) - (\ref{8.3}) is an isomorphism if and only if $\varphi$ is an isomorphism between the loop algebras, $\lambda=\mu=:\epsilon\in\{\pm 1\}$ and
\begin{eqnarray}\label{e8.7}
\varphi(u')&=&-\epsilon[x_\varphi,\varphi u]_0+\epsilon\varphi(u')\\
\alpha(u)&=&-\epsilon(x_\varphi,\varphi u)
\end{eqnarray}
for all $u\in L_\alg(\g,\sigma)$.

The theorem implies that $\hat\varphi\mapsto\check\varphi\mapsto\varphi$ induces isomorphisms $\Aut_{(,)}\hat L_\alg(\g,\sigma)\to\Aut L'_\alg(\g,\sigma)\to\Aut L_\alg(\g,\sigma)$ and shows that automorphisms of finite order of $\hat L_\alg(\g,\sigma)$ are contained in $\Aut_{(,)}\hat L_\alg(\g,\sigma)$. Thus we have:

\begin{corollary}\label{8.10}
The mappings $\hat\varphi\mapsto\check\varphi$ and $\check\varphi\to\varphi$ induce bijections between the sets of elements of finite order in $\Aut\hat L_\alg(\g,\sigma)$, $\Aut\hat L'_\alg(\g,\sigma)$ and $\Aut L_\alg(\g,\sigma)$, respectively, as well as between their conjugacy classes.
\end{corollary}

\subsection{Automorphisms of finite order}
From now on we assume $\g$ to be compact if $\F=\R$. Thus $\g$ is either complex or compact (and simple). In this subsection we are going to attach to each automorphism of finite order an invariant like in the smooth case and to prove surjectivity of the map from the set of conjugacy classes to the set of invariants. Injectivity (the hard part) will be shown in the following subsection.

Thus let $\varphi\in\Aut L_\alg(\g,\sigma)$ be of finite order. If $\varphi$ is of the first kind then it is standard by Corollary \ref{8.7} and thus extends to a unique automorphism of $L(\g,\sigma)$ whose invariant we define to be the invariant of $\varphi$. If $\varphi$ is of the second kind there exists a unique $r>0$ such that $\tau_r\varphi\tau_r^{-1}$ is standard and we define the invariant of $\varphi$ to be the invariant of (the extension of) $\tau_r\varphi \tau_r^{-1}$. It follows from our previous results (Propositions~\ref{4.7} and \ref{5.6}) that this invariant indeed is invariant under quasiconjugation with isomorphisms of the first kind (and arbitrary isomorphisms if $\varphi$ is of the second kind or $\varphi$ is an involution). For the only problem is conjugation with $\tau_r$ (if $\varphi$ is of the first kind). But the invariant of $\tau_r\varphi\tau_r^{-1}=\ {^r\varphi}$ varies by construction smoothly with $r$ and hence is constant as the set of invariants is discrete.

If $\hat\varphi\in\Aut\hat L_\alg(\g,\sigma)$ is of finite order we define the invariant of $\hat\varphi$ to be that of the induced $\varphi$.

We have already observed in the smooth case that each element of  $\J^q_i(\g)$ can be realized on $L(\g,\sigma)$ for some $\sigma\in\Aut\g$ as the invariant of an automorphism $\varphi$ of the form $\varphi u(t)=\varphi_0(u(\epsilon t+t_0))$ where $\varphi_0\in\Aut\g$ is constant (Propositions~\ref{4.6} and \ref{5.5}). Since this $\varphi$ is algebraic the only problem is to show that $\sigma$ can be chosen to be of finite order.

\begin{lemma}\label{8.11}
Let $G$ be a compact Lie group and $g_\pm\in G$ with $g_+^2=g_-^2$. Then there exists $h\in G_0$ such that $hg^{-1}_-h^{-1}g_+$ is of finite order.
\end{lemma}

\begin{proof}
The compact abelian group $\overline{\{(g_-^{-1}g_+)^n\mid n\in\Z\}}$ is isomorphic to the product of a torus $T$ and a finite group $F$. In particular $g_-^{-1}g_+=e^{X}\cdot f$ for some $f\in F$ and some $X$ from the Lie algebra of $T$. Now, $g^2_+=g_-^2$ is equivalent to
\[
g_-(g^{-1}_-g_+)g_-^{-1}=(g_-^{-1}g_+)^{-1}
\]
whence conjugation with $g_-$ induces the inverse mapping on $T$. In particular $g_-e^{tX}g_-^{-1}=e^{-tX}$ for all $t$. Let $h:=e^{-{1\over2}X}$. Then $hg_-^{-1}h^{-1}g_+=h^2g_-^{-1}g_+=f$ has finite order.
\end{proof}

\begin{remark}\label{8.12}\rm
If $\g$ is complex (i.e. $\F=\C$) and $\u$ is a compact real form of $\g$, recall that
\[
\begin{array}{rcl}
\Aut_\u\g \times i\u&\to&\Aut\g\\
(\varphi,X)&\mapsto&e^{\ad X}\cdot\varphi
\end{array}
\]
is a diffeomorphism where $\Aut_\u\g=\{\varphi\in\Aut\g\mid\varphi\u=\u\}$ and  that moreover any compact subgroup of $\Aut\g$ is conjugate to a subgroup of $\Aut_\u\g$ (cf. Appendix B).
\end{remark}

\begin{lemma}\label{8.13}
Let $\F=\C$ and $\varphi_\pm\in\Aut\g$ with $\varphi_+^2=\varphi_-^2$ and of finite order. Then there exists $\alpha\in((\Aut\g)^{\varphi_+^2})_0$ and a compact real form of $\g$ that is invariant under $\varphi_+$ and $\alpha\varphi_-\alpha^{-1}$.
\end{lemma}

\begin{proof}
Since $\varphi_\pm$ are of finite order there are compact real forms $\u_\pm$ invariant under $\varphi_+$ and $\varphi_-$, respectively. Let $\u:=\u_+$. Then $\u_-=\alpha^{-1}\u$ for some $\alpha\in\Aut\g$ as compact real  forms are conjugate, and thus $\alpha\varphi_-\alpha^{-1}(\u)=\u$. By the remark above, we may assume $\alpha=e^{\ad X}$ for some $X\in i\u$. Since $\u$ is invariant under $\varphi_+$ and $e^{\ad X}\varphi_-e^{-\ad X}$ and hence also under their squares $\psi:=\varphi_+^2$ and $\chi:=e^{\ad X}\psi e^{-\ad X}$. Now $e^{-\ad X}\cdot\chi=e^{-\ad\psi X}\cdot \psi$ implies by Remark~\ref{8.12} above $\psi X=X$ and hence $\alpha\in((\Aut\g)^\psi)_0$.
\end{proof}

\begin{theorem}\label{8.14}
Any element of $\J_i^q(\g)$ (where $i\in\{1,2\}$ and $q\in\N$ with $q$ even if $i=2$) can be realized on $L_\alg(\g,\sigma)$ for some $\sigma$ of finite order as the invariant of a standard automorphism of the form $\varphi u(t)=\varphi_0(u(\epsilon t+t_0))$ where $\varphi_0\in\Aut\g$ is constant. In particular the mappings
\[
\Aut_i^q L_\alg(\g,\sigma)\to \J_i^q(\g,\sigma)
\]
are surjective.
\end{theorem}

\begin{proof}
\begin{enumerate}
\item[(i)] The last statement follows from the claim before as any element of $\J_i^q(\g,\sigma)$ can be realized as the invariant of an automorphism on some $L_\alg(\g,\tilde\sigma)$ which is isomorphic to $L_\alg(\g,\sigma)$. But the corresponding automorphism on $L_\alg(\g,\sigma)$ has the same invariant.
\item[(ii)] Let $(p,\varrho,[\beta])\in \J_1^q(\g)$ where $\varrho$ has finite order and $\beta\in(\Aut\g)^\varrho$. We realize this invariant as in the $C^\infty$-case by $\varphi u(t)=\varphi_0(u(t+{p\over q}2\pi))$ on $L_\alg(\g,\sigma)$ where $\varphi_0 $ and $\sigma$ are certain products of powers of $\varrho$ and $\beta$. Therefore it suffices to show that $\beta$ can be replaced in its equivalence class by an element of finite order, i.e. to show the existence of an $\alpha\in((\Aut\g)^\varrho)_0$ with $\beta\cdot\alpha$ of finite order. If $\g$ is compact this can be achieved by Lemma~\ref{2.16} directly. If $\g$ is complex we first remark that $\varrho$ leaves a compact real form $\u$ invariant and that hence $(\Aut\g)^\varrho=(\Aut_\u\g)^\varrho\cdot\{e^{\ad X}\mid X\in i\u^\varrho\}$ by \ref{8.12}. Therefore Lemma~\ref{2.16} with $G:=(\Aut\g)^\varrho$ and $H:=(\Aut_\u\g)^\varrho$ yields the result.
\item[(iii)] Let $[\varphi_+,\varphi_-]\in\J_2^q(\g)$. An automorphism with this invariant is given by $\varphi u(t):=\varphi_+(u(-t))$ on $L_\alg(\g,\varphi_-^{-1}\varphi_+)$ provided $\varphi_-^{-1}\varphi_+$ has finite order. Since $(\varphi_+,\varphi_-)$ is equivalent to $(\varphi_+,h\varphi_-h^{-1})$ for any $h\in((\Aut\g)^{\varphi^2_+})_0$ it is therefore enough to find such an $h$ with $h\varphi_-^{-1}h^{-1}\varphi_+$ of finite order. If $\g$ is compact we can apply \ref{8.11} directly. If $\g$ is complex we may assume by \ref{8.13} that $\varphi_\pm$ leave a compact real form $\u$ invariant and can then apply \ref{8.11} to $G:=\Aut_\u\g$ and $\varphi_\pm\in G$.
\end{enumerate}
\end{proof}

\subsection{Injectivity of $\Aut^q_i L_\alg(\g,\sigma)/\Aut_1 L_\alg(\g,\sigma)\to \J_i^q(\g,\sigma)$}
To prove injectivity, we use Levstein's result that in case $\g$ complex, any automorphism of finite order of $\hat L_\alg(\g,\sigma)$ leaves a Cartan subalgebra invariant \cite{Lev}. This also holds if $\g$ is compact by almost the same reasoning as was shown in \cite{Rou2}, Theorem 3.8. From these results we get:

\begin{proposition}\label{8.15}
Any automorphism of $L_\alg(\g,\sigma)$ of finite order is conjugate to a standard automorphism $\varphi$ with $\varphi u(t)=\varphi_t(u(\epsilon t+t_0))$ and $\varphi_t=e^{\ad tX}\cdot\varphi_0$ where $X$ is contained in some Cartan subalgebra $\a$ of $\g^\sigma$ and $\varphi_0\in\Aut\g$ leaves $\a$ invariant.
\end{proposition}

\begin{proof}
Let $\hat\varphi\in\Aut\hat L_\alg(\g,\sigma)$ be the (unique) extension of finite order of the given automorphism which we may assume to be standard of the form $\varphi u(t)=\varphi_t(u(\epsilon t+t_0))$ by~\ref{8.7}. Due to Levstein and Rousseau, $\hat\varphi$ leaves a Cartan subalgebra invariant. Since Cartan subalgebras are conjugate (see 4.4 of \cite{Rou2} in case $\g$ compact) we may assume that $\hat\varphi$ leaves $\hat\h:=\{u\in L_\alg(\g,\sigma)\mid u(t)\in\a$ constant$\}\oplus\F c\oplus\F d$ invariant. Since $\hat\varphi d=\epsilon d+x_\varphi+\nu c$ is contained in this Cartan subalgebra, $x_\varphi$ is constant and lies in $\a$. Equation~(\ref{e8.7}) therefore implies $\varphi'_t\varphi_t^{-1}=-\epsilon\ad x_\varphi$ and hence $\varphi_t=e^{\ad tX}\cdot\varphi_0$ for $X=-\epsilon x_\varphi\in\a$ and some $\varphi_0\in\Aut\g$. Moreover $\varphi_0(\a)=\a$ as $\hat\varphi$ leaves $\hat\h$ invariant.
\end{proof}

In the next step we even get rid of the $e^{\ad tX}$ factor by a quasiconjugation (i.e. a change of $\sigma$).

\begin{lemma}\label{8.16}
A curve $\varphi_t$  of automorphisms of $\g$ of the form $e^{\ad tX}\cdot\varphi_0$ is algebraic if and only if $e^{\ad 2\pi X}$ has finite order. This condition is satisfied if
$\varphi_{t+2\pi}=\tilde\sigma\varphi_t\sigma^{-\epsilon}$ holds for all $t\in\R$ for some $\sigma,\tilde\sigma\in\Aut\g$ of finite order and some $z\in\Z$.
\end{lemma}

\begin{proof}
\begin{enumerate}
\item[(i)] Recall that $\varphi_t$ is algebraic if $\varphi_t={\sum\limits_{\vert m\vert\le N}}\varphi_m e^{i m t/L}$ for some $N,L\in\N$ and $\varphi_m\in\mbox{End}(\g)$ in case $\g$ is complex, and $\varphi_t$ is algebraic if $(\varphi_t)_\C$ is algebraic in case $\g$ is real. In particular we may assume $\g$ complex in what follows.
\item[(ii)] If $\varphi_t=e^{\ad tX}\cdot\varphi_0$ is algebraic like in (i) then $\varphi_{t+2\pi L}=\varphi_t$ and thus $(e^{\ad 2\pi X})^L=\id$.
\item[(iii)] Conversely, if $e^{\ad 2\pi LX}=\id$ for some $L\in\N$ then $\ad X$ is semisimple and $X$ is contained in some Cartan subalgebra. On each root space $\g_\alpha$, $e^{\ad tX}$ acts as $e^{t\alpha(X)}\cdot \id$ and $e^{\ad 2\pi LX}=\id$ implies $\alpha(X)\in{i\over L}\Z$. Thus $e^{\ad tX}$ and hence $\varphi_t$ are algebraic.
\item[(iv)] $\varphi_{t+2\pi}=\tilde\sigma\varphi_t\sigma^{-\epsilon}$ is equivalent to $e^{\ad 2\pi X} =\tilde\sigma\varphi_0\sigma^{-\epsilon}\varphi_0^{-1}$ and $\tilde\sigma X=X$. Therefore $e^{\ad 2\pi X}$ and $\tilde\sigma$ commute. Hence also $\tilde\sigma$ and $\varphi_0\sigma^{-\epsilon}\varphi_0^{-1}$ commute and $e^{\ad 2\pi X}$ is of finite order if $\sigma$ and $\tilde\sigma$ are of finite order.
\end{enumerate}
\end{proof}

\begin{proposition}\label{8.17}
Any automorphism of $L_\alg(\g,\sigma)$ of finite order is quasiconjugate (by a standard automorphism) to an isomorphism of the form $\varphi u(t)=\varphi_0(u(\epsilon t+t_0))$ where $\varphi_0\in\Aut\g$ is constant.
\end{proposition}

\begin{proof}
By \ref{8.15} we may assume that the given automorphism $\tilde\varphi$ of finite order is of the form $\tilde\varphi u(t)=\tilde\varphi_t(u(\epsilon t+t_0))$ with $\tilde\varphi_t=e^{\ad tX}\cdot\varphi_0$ where $X$ is contained in some Cartan subalgebra $\a$ of $\g^\sigma$ and $\varphi_0(\a)=\a$. We are then looking for a $\tilde\sigma\in\Aut\g$ of finite order and an isomorphism $\psi:L_\alg(\g,\sigma)\to L_\alg(\g,\tilde\sigma)$ such that $\varphi:=\psi\tilde\varphi\psi^{-1}$ is of the form $\varphi u(t)=\varphi_0(u(\epsilon t+t_0))$. With the ansatz $\psi u(t)=\psi_t(u(t))$ this is equivalent to find $\tilde\sigma,\psi_t\in\Aut\g$ with
\begin{enumerate}
\item[(i)] $\psi_t\tilde\varphi_t\psi_{\epsilon t+t_0}^{-1}$ is constant
\item[(ii)] $\psi_{t+2\pi}=\tilde\sigma\psi_t\sigma^{-1}$ and
\item[(iii)] $t\mapsto\psi_t$ is algebraic.
\end{enumerate}
Assuming $\psi_t=e^{\ad tY}$ for some $Y\in\a$ these equations are in turn equivalent to
\begin{enumerate}
\item[(i')] $Y-\epsilon \varphi_0 Y+X=0$
\item[(ii')] $\tilde\sigma=e^{\ad 2\pi Y}\sigma$
\item[(iii')] $e^{\ad 2\pi Y}$ has finite order.
\end{enumerate}

These can be solved by $Y:={1\over q}{\sum\limits_{j=1}^{q-1}}\epsilon^j j\varphi_0^j X$ and $\tilde\sigma:=e^{\ad 2\pi Y} \sigma$ where $q$ is the order of $\tilde\varphi$. In fact, since $\tilde\varphi_t$ is algebraic, $e^{\ad 2\pi X}$ as well as $e^{\ad 2\pi\varphi_0^{j}X}=\varphi_0^je^{\ad 2\pi X}\varphi_0^{-j}$ have finite order by \ref{8.16}. Therefore also $e^{\ad 2\pi Y}$ and $\tilde\sigma$ are of finite order and (ii') and (iii') are satisfied. Finally, $Y-\epsilon\varphi_0 Y={1\over q}{\sum\limits^q_{j=1}}\epsilon^j\varphi_0^j X-\epsilon^q\varphi_0^q X=-X$ as $\tilde\varphi^q=\id$ implies ${\sum\limits^{q-1}_{j=0}}\epsilon^j\varphi_0^j X=0$ and hence ${\sum\limits^q_{j=1}}\epsilon^j\varphi_0^j X=0$ and $\epsilon^q\varphi_0^q X=X$.
\end{proof}

\begin{lemma}\label{8.18}
Let $\g$ be complex $(\F=\C)$ and $\varphi\in\Aut L_\alg(\g,\sigma)$ be an automorphism of finite order with $\varphi u(t)=\varphi_0(u(\epsilon t+t_0))$ for some $\varphi_0\in\Aut\g$, $\epsilon\in\{\pm 1\}$, and $t_0\in\R$. Then $\g$ has a compact real form invariant under $\varphi_0$ and $\sigma$.
\end{lemma}

\begin{proof}
Since $\varphi$ has finite order there exists $q\in\N$ and $p\in\Z$ such that $\varphi_0^q\sigma^p=\id$. Therefore $\varphi_0$ has finite order since $\sigma$ has finite order by assumption. Now, $\varphi_0=\sigma\varphi_0\sigma^{-\epsilon}$ implies that the group generated by $\varphi_0$ and $\sigma$ is finite. It therefore leaves a compact real form of $\g$ invariant (cf. Appendix B).
\end{proof}

The next lemma is well known and essentially says that the so-called $\sigma$- and Hermann actions on a compact Lie group are hyperpolar (\cite{HPTT}).

\begin{lemma}\label{8.19}
Let $\tilde G$ be a compact connected Lie group with Lie algebra $\tilde\g$ and $\sigma,\varrho_+,\varrho_-\in\Aut\tilde\g$ with $\varrho_\pm^2=\id$. Then there exist for any $g\in\tilde G$.
\begin{enumerate}
\item[(i)] $h\in \tilde G$ and $X\in\tilde\g^\sigma$ such that
\[
h g\sigma(h)^{-1}=e^{X}
\]
\item[(ii)] $k_+\in(\tilde G^{\varrho_+})_0$ and $X\in\tilde\g$ such that $\varrho_+X=\varrho_- X=-X$ and
\[
k_+ g\varrho_-(k_+g)^{-1}=e^{X}\ .
\]
\end{enumerate}
\end{lemma}

\begin{proof}
We endow $\tilde G$ with a biinvariant metric which is also invariant under automorphism by using a multiple of the Killing form. If a compact group $H$ acts isometrically on $\tilde G$ then the image of the normal space to the orbit $H(e)$ at $e$ under the exponential mapping meets all other orbits. In fact, a shortest connection from an arbitrary $g\in \tilde G$ to $H(e)$ is a geodesic that meets $H(e)$ orthogonally and thus can be moved by the group action to a geodesic that hits $H(e)$ orthogonally at $e$ and starts on $H(g)$.

We now consider the following two isometric actions on $\tilde G$:
\begin{enumerate}
\item[(i')] The $\sigma$-action of $\tilde G$ on itself by
\[
h.g:=h g\sigma(h)^{-1}
\]
and
\item[(ii')] The action of $K_+ \times K_-$, where $K_\pm:=(\tilde G^{\varrho_\pm})_0$, by
\[
(k_+,k_-).g:=k_+g k_-^{-1}\ .
\]
\end{enumerate}
The tangent spaces to their orbits at $e$ are $\{X-\sigma X\mid X\in\tilde\g\}$ and $\tilde\g^{\varrho_+}+\tilde\g^{\varrho_-}$ and hence their normal spaces are $\tilde\g^\sigma$ and $\{X\in\tilde\g\mid\varrho_\pm X=-X\}$, respectively. Thus (i) follows immediately and in case (ii) we find for any $g\in\tilde G$ elements $k_\pm\in K_\pm$ and $Y\in\{ X\in\tilde\g\mid \varrho_\pm X=-X\}$ with $k_+ gk_-^{-1}=e^{Y}$. The last equation implies $(k_+g)\cdot\varrho_-(k_+g)^{-1}=e^{2Y}$.
\end{proof}

\begin{proposition}\label{8.20}
Two automorphisms $\varphi\in\Aut L_\alg(\g,\sigma)$ and $\tilde\varphi\in\Aut L_\alg(\g,\tilde\sigma)$ of order $q$ with the same invariants are quasiconjugate by an isomorphism of the first kind.
\end{proposition}

\begin{proof}
Since invariants do not change under quasiconjugation we may assume by~\ref{8.17} $\varphi$ and $\tilde\varphi$ to be of the form $\varphi u(t)=\varphi_0u(\epsilon t+t_0)$ resp. $\tilde\varphi u(t)=\tilde\varphi_0 u(\epsilon t+t_0)$ where $\varphi_0,\tilde\varphi_0\in\Aut\g$, $\epsilon\in\{\pm 1\}$ and $t_0=p/q 2\pi$ with $p\in\{0,1,\dots,q-1\}$ and $p=0$ if $\epsilon=-1$. We may furthermore conjugate $\varphi_0$ and $\sigma$ simultanously by an arbitrary $\psi_0\in\Aut\g$ because this corresponds to a quasiconjugation of $\varphi$ by $\psi:L_\alg(\g,\sigma)\to L_\alg(\g,\psi_0\sigma\psi_0^{-1})$ with $\psi u(t)=\psi_0(u(t))$, and the same remark applies to $\tilde\varphi_0$ and $\tilde\sigma$. In particular, in case $\F=\C$, we may assume by Lemma~\ref{8.18} that $\varphi_0,\tilde\varphi_0,\sigma$ and $\tilde\sigma$ leave a compact real form of $\g$ invariant and that thus $\varphi$ and $\tilde\varphi$ are complexifications of real automorphisms. Since these have also equal invariants (cf. the proof of~\ref{7.5} (i)) it is enough to consider the case $\g$ compact as we shall do in the following.

We study the two cases $\epsilon=\pm 1$ separately.

If $\epsilon=1$ the invariant of $\varphi$ is $(p,\alpha\varphi_0^{q'}\sigma^{p'}\alpha^{-1}, [\alpha\varphi_0^{-l}\sigma^m\alpha^{-1}])$ where $\alpha\in\Aut\g$ and $p',q',l,m$ are the integers with $p'/q'=p/q$, $lp'+ mq'=1$ and $0\le l<q'$ (cf.~\ref{4.3}). By eventually conjugating $\varphi_0$ and $\sigma$ simultanously by $\alpha$ we may assume $\alpha=\id$. Applying the same to $\tilde\varphi$, equality of the invariants yields $(p,\varphi_0^{q'}\sigma^{p'},[\varphi_0^{-l}\sigma^m])=(\tilde p,\tilde\varphi_0^{\tilde q'}\tilde\sigma^{\tilde p'}, [\tilde\varphi_0^{-\tilde l}\tilde\sigma^{\tilde m}])$ which means $\tilde p=p$, $\tilde q'=q'$, $\tilde p'=p'$, $\tilde l=l$, $\tilde m=m$,
\[
\tilde\varphi_0^{q'}\tilde\sigma^{p'}=\varphi_0^{q'}\sigma^{p'}=:\varrho\ ,
\]
and
\[
\tilde\varphi_0^{-l}\tilde\sigma^m=\gamma\delta\varphi_0^{-l}\sigma^m\delta^{-1}\ ,
\]
where $\gamma\in((\Aut\g)^\varrho)_0$ and $\delta\in(\Aut\g)^\varrho$. By eventually conjugating $\varphi_0$ and $\sigma$ further by $\delta$ we may assume $\delta=\id$. Let $\beta:=\varphi_0^{-l}\sigma^m$, $\tilde\beta:=\tilde\varphi_0^{-l}\tilde\sigma^m$, and $\tilde G:=(\Aut\g)^\varrho$. Since $\varphi_0$ and $\sigma$ as well as $\tilde\varphi_0$ and $\tilde\sigma$ commute (due to $\epsilon=1)$, $\varphi_0,\tilde\varphi_0,\sigma,\tilde\sigma,\beta,\tilde\beta$ are contained in $\tilde G$. Moreover $\gamma=\tilde\beta\beta^{-1}\in\tilde G_0$. We now try to find an isomorphism $\psi: L_\alg(\g,\sigma)\to L_\alg(\g,\tilde\sigma)$ of the form $\psi u(t)=\psi_t u(t)$ with $\psi\varphi\psi^{-1}=\tilde\varphi$, that is an algebraic curve $\psi_t$ in $\Aut\g$ with $\psi_{t+2\pi}=\tilde\sigma\psi_t\sigma^{-1}$ and $\psi_t\varphi_0\psi^{-1}_{t+{p\over q}2\pi}=\tilde\varphi_0$. Since ${1\over q'}=l{p\over q}+m$ these two equations imply
\begin{equation}\label{*}
\psi_{t+{2\pi\over q'}}=\tilde\beta\psi_t\beta^{-1}
\end{equation}
and in fact are equivalent to (\ref{*}) if we choose $\psi_t\in\tilde G$. We therefore make the ansatz $\psi_t:=\psi_0e^{\ad tX}$ with $X\in\g^\varrho$ and $\psi_0\in\tilde G$. Then (\ref{*}) is equivalent to
\[
e^{\ad{2\pi\over q'} X}=\psi_0^{-1}\gamma\beta\psi_0\beta^{-1}\mbox{ and } \beta X=X\ .
\]
This equation has a solution $(\psi_0,X)$ by Lemma~\ref{8.19} (i), and by~\ref{8.16} $\psi_t$ is algebraic.

If $\epsilon=-1$ the invariants of $\varphi$ and $\tilde\varphi$ are $[\varphi_0,\varphi_0\sigma^{-1}]$ and $[\tilde\varphi_0,\tilde\varphi_0\tilde\sigma^{-1}]$, respectively. Note that quasiconjugation of $\varphi$ by $\psi:L_\alg(\g,\sigma)\to L_\alg(\g,\sigma^{-1})$ with $\psi u(t):=u(-t+\pi)$ reverses the order of $\varphi_0$ and $\varphi_0\sigma^{-1}$ as $\psi\varphi\psi^{-1}u(t)=\varphi_0\sigma^{-1}(u(-t))$. Thus  we may assume that equality of the invariants implies $\tilde\varphi_0=\alpha\varphi_0\alpha^{-1}$ and $\tilde\varphi_0\tilde\sigma^{-1}=\beta\varphi_0\sigma^{-1}\beta^{-1}$ for some $\alpha,\beta\in\Aut\g$ with $\alpha^{-1}\beta\in((\Aut\g)^{\varphi_0^2})_0$. Furthermore we may assume $\alpha=\id$ by quasiconjugating $\varphi$ by $L_\alg(\g,\sigma)\to L_\alg(\g,\alpha\sigma\alpha^{-1})$, $u(t)\mapsto\alpha(u(t))$, which maps $(\varphi_0,\sigma)$ to $(\alpha\varphi_0\alpha^{-1},\alpha\sigma\alpha^{-1})$. Let $\tilde G:=(\Aut\g)^{\varphi_0^2}$. Since $\varphi$ and $\tilde\varphi$ are of the second kind, $\varphi_0=\sigma\varphi_0\sigma$, i.e. $(\varphi_0\sigma^{-1})^2=\varphi_0^2$ and similarly $(\tilde\varphi_0\tilde\sigma^{-1})^2=\tilde\varphi_0^2=\varphi_0^2$. Thus $\varphi_0,\varphi_0\sigma^{-1},\tilde\varphi_0$, and $\tilde\varphi_0\tilde\sigma^{-1}\in\tilde G$. Moreover $\tilde\varphi_0\tilde\sigma^{-1}=\beta\varphi_0\sigma^{-1}\beta^{-1}$ for some $\beta\in\tilde G_0$ and thus $\tilde\sigma=\beta\sigma\varphi_0^{-1}\beta^{-1}\varphi_0$.

Let $\varphi_+:=\varphi_0$ and $\varphi_-:=\varphi_0\sigma^{-1}$. Conjugation with $\varphi_\pm$ defines automorphisms $\varrho_\pm:\tilde G\to\tilde G$ with $\varrho_\pm^2=\id$. Hence there exist by Lemma~\ref{8.19} (ii) $\psi_0\in\tilde G^{\varrho_+}=\tilde G^{\varphi_0}$ and $X\in\g$ with $\varphi_\pm X=-X$ such that $e^{\ad 2\pi X}=(\psi_0^{-1}\beta)\varrho_-(\psi_0^{-1}\beta)^{-1}$, that is with
\[
\psi_0e^{\ad 2\pi X}=\tilde\sigma\psi_0\sigma^{-1}\ .
\]
Let $\psi_t:=\psi_0e^{\ad tX}$. Then $\psi_{t+2\pi}=\tilde\sigma\psi_t\sigma^{-1}$ and $\psi_t\varphi_0\psi_{-t}^{-1}=\varphi_0=\tilde\varphi_0$. Therefore $\psi:L_\alg(\g,\sigma)\to L_\alg(\g,\tilde\sigma)$ with $\psi u(t))=\psi_t u(t)$ is an isomorphism of the first kind that conjugates $\varphi$ into~$\tilde\varphi$.
\end{proof}

Combining the results of Theorem~\ref{8.14}, Proposition~\ref{8.20} and Corollary~\ref{8.10} we obtain:

\begin{theorem}\label{8.21}
The mappings
\[
\Aut_i^q\hat L_\alg(\g,\sigma)/\Aut_1\hat L_\alg(\g,\sigma)\to \J_i^q(\g,\sigma)
\]
and
\[
\Aut^q_i L_\alg(\g,\sigma)/\Aut_1 L_\alg(\g,\sigma)\to \J_i^q(\g,\sigma)
\]
induced by associating to each automorphism its invariant, are bijections. Moreover in case $i=2$ or $q=2$, $\Aut_1\hat L_\alg(\g,\sigma)$ and $\Aut_1L(\g,\sigma)$ can be replaced by $\Aut\hat L_\alg(\g,\sigma)$, resp., $\Aut L_\alg(\g,\sigma)$.\qed
\end{theorem}

The theorem together with Corollaries~\ref{4.12} and \ref{5.9} imply that the classification of automorphisms of finite order is the same in the smooth and the algebraic category. This applies in particular to involutions which are hence classified by Tables 2 and 3 also in the algebraic case.

\subsection{Real forms and Cartan decompositions}
The discussion of conjugate linear automorphisms and real forms in Chapter~7 carries over to the algebraic setting almost word by word.

In particular compact real forms of $\hat L_\alg(\g,\sigma)$ and $L_\alg(\g,\sigma)$, where $\g$ is complex, are defined and unique up to conjugation. The conjugacy classes of real forms of type $i\in\{1,2\}$ are e.g.~on $\hat L_\alg(\g,\sigma)$ in bijection with $\overline{\Aut}^2_i\hat L_\alg(\g,\sigma)/\Aut\hat L_\alg(\g,\sigma)$ and if $\u$ is a $\sigma$-invariant compact real form of $\g$, complex conjugate extensions of automorphisms of $\hat L_\alg(\u,\sigma)$ induce bijections $\Aut^1_1\hat L_\alg(\u,\sigma)\cup\Aut^2_1\hat L_\alg(\u,\sigma)/$ $\Aut\hat L_\alg(\u,\sigma)\to \overline{\Aut}^2_1\hat L_\alg(\g,\sigma)/$ $\Aut\hat L_\alg(\g,\sigma)$ and $\Aut^2_2\hat L_\alg(\u,\sigma)/\Aut\hat L_\alg(\u,\sigma)\to$ $\overline{\Aut}^2_2\hat L_\alg(\g,\sigma)/$ $\Aut\hat L_\alg(\g,\sigma)$, respectively. Moreover surjectivity of these mappings yields the existence of Cartan decompositions of real forms while injectivity yields their uniqueness up to conjugation.

%%%%%%%%%%%%%%%%%%%%%%%%%%%%%%%%%%%%%%%%%%%%%%%%%%%%%%%%%
\newpage

\begin{appendix}
\begin{center}
\bf \LARGE Appendix
\end{center}
\section{$\pi_0 ((\Aut\g)^\varrho)$ and representatives of its conjugacy classes}

Let $\g$ be a compact (real) simple Lie algebra and $\varrho\in\Aut\g$ an involution. The group $\pi_0((\Aut\g)^\varrho)$ of connected components of $(\Aut\g)^\varrho$ has been determined by Cartan \cite{Car} and Takeuchi \cite{Tak} (cf.~also \cite{Mur}, \cite{MatH}). A simplified, but still quite involved computation of these groups is contained in \cite{Loo} where one can also find a table of them.

Our purpose here is to describe representatives of their conjugacy classes and thereby verifying Table~1 of Section~6. But actually we will determine also the groups themselves since it turns out that this does not need much extra work.  Moreover some of the extra work (Lemmas~\ref{A.2} and \ref{A.8}) is needed also for other purposes, namely  the classification of involutions of the second kind.

In the following we fix $\g$ and $\varrho$ and let $\g=\k+\p$ be the splitting of $\g$ into the $\pm 1$ eigenspaces of $\varrho$. We let $\tau_p$, $J$ and $\mu$ as in Section~6.

If $\g$ is \textbf{classical} we only use the following well known facts (Lemmas~\ref{A.1} -- \ref{A.4}) to determine $\pi_0((\Aut\g)^\varrho)$ and representatives of its conjugacy classes.

\begin{lemma}\label{A.1}
$\pi_0((\Int\g)^\varrho)$ is a normal subgroup of $\pi_0((\Aut\g)^\varrho)$. The quotient $F$ is isomorphic to a subgroup of $\Aut\g/\Int\g$.
\end{lemma}

\begin{proof}
The first statement is clear since $(\Int\g)^\varrho$ is normal in $(\Aut\g)^\varrho$. The second statement follows from $\pi_0((\Aut\g)^\varrho)/\pi_0((\Int\g)^\varrho)$ being isomorphic to $(\Aut\g)^\varrho/(\Int\g)^\varrho$ which is embedded naturally in $\Aut\g/\Int\g$.
\end{proof}

In particular $\pi_0((\Aut\g)^\varrho)/\pi_0((\Int\g)^\varrho)$ is $1$ or $\Z_2$ unless $\g\cong\so(8)$ (in which case it is a subgroup of $S_3)$.

\begin{lemma}\label{A.2}
Let $\Int\k=\{e^{\ad X}:\k\to\k\mid X\in\k\}$. Then the mapping
\[
\pi_0((\Aut\g)^\varrho)\to \Aut\k/\Int\k
\]
induced by restriction is injective if $\varrho$ is inner and has kernel $\{1,\varrho\}$ if $\varrho$ is outer. Here we denote the image of $\varrho\in(\Aut\g)^\varrho$ in $\pi_0((\Aut\g)^\varrho)$ also by $\varrho$.
\end{lemma}

\begin{remark}\label{A.3}
{\rm Note that $\k$ is not necessary semisimple. This happens precisely if the corresponding symmetric space is Hermitian, in which case $\k\cong\R+\k^*$ with $\k^*$ semisimple and $\Aut \k/\Int \k\cong\Aut\R\times \Aut \k^*/\Int \k^*$. Since $\pi_0((\Aut\g)^\varrho)$ is finite its image lies in that case in $\{\pm \id\}\times \Aut \k^*/\Int \k^*$.}
\end{remark}

{\sl Proof of Lemma~\ref{A.2}.}
An element from the kernel can be represented by a $\varphi\in(\Aut\g)^\varrho$ with $\varphi_{\vert_\k}\in\Int\k$ and thus even with $\varphi_{\vert_\k}=\id$ since $(\Int\k)\hookrightarrow((\Aut\g)^\varrho)_0$ naturally. The action of this $\varphi$ on $\p$ commutes with the restriction of the action of $K:=\{e^{\ad  X}:\g\to\g\mid X\in\k\}$ on $\p$ which is irreducible and either (i) a real or (ii) a complex representation. The latter occurs precisely in the Hermitian case. Therefore either (i) $\varphi_{\vert_\p}=\pm \id$ (hence $\varphi=\id$ or $\varphi=\varrho$ and the claim follows) or (ii) $\varphi$ lies in the circle $\{e^{\ad tJ_0}\mid t\in\R\}$ (and thus $\varphi$ represents a trivial element) where $J_0$ spans the center of $\k$. \qed

Lemma~\ref{A.2} already implies $\pi_0((\Aut\g)^\varrho)=\{1,\varrho\}$ if $\varrho$ is outer unless the corresponding symmetric pair is $(\so(2n),\so(n)+\so(n))$ with $n$ odd or $(\su(2n),\so(2n))$. For in all other cases $\Aut\k/\Int\k$ is trivial.

\begin{lemma}\label{A.4}
The triality automorphism $\vartheta$ of $\so(8)$ commutes with $\Ad \tau_4$.
\end{lemma}

\begin{proof}
Let
$
X_1:={\tiny
\left(
\begin{array}{ccccc}
0&1&&&\cr
-1&0&&&\cr
&&0&&\cr
&&&\ddots&\cr
&&&&0
\end{array}
\right)}
,
\dots,X_4:={\tiny
\left(
\begin{array}{ccccc}
0&&&&\\
&\ddots&&&\\
&&0&&\\
&&&0&1\\
&&&-1&0
\end{array}
\right)}
$
be a basis of the standard torus of $\so(8)$ and $l_1,\dots,l_4$ the dual linear forms. Then $\pm l_i\pm l_j\ (1\le i<j\le 4)$ are the roots and $\alpha_1:=l_1-l_2$, $\alpha_2:=l_2-l_3$, $\alpha_3:= l_3-l_4$, $\alpha_4:=l_3+l_4$ is a basis of the root system with diagram

%%%%%%%%%%%%%%%%%%%%%%%%%%%%%%%%%%%%
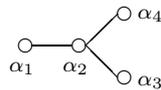
\begin{center}\begin{picture}(110,30)

           \put(36,12){\circle{5}}
           \put(56,12){\circle{5}}
           \put(73,24){\circle{5}}
           \put(73,0){\circle{5}}

           \put(38.5,12){\line(1,0){15}}
           \put(58.5,12){\line(1,1){12}}
           \put(58.5,12){\line(1,-1){12}}

                    \put(30,2){$\scriptstyle \alpha_1$}
           \put(50,2){$\scriptstyle \alpha_2$}
           \put(78,22){$\scriptstyle \alpha_4$}
           \put(78,-3){$\scriptstyle \alpha_3$}

\end{picture}\end{center}

%%%%%%%%%%%%%%%%%%%%%%%%%%%%%%%%%%%%

The triality automorphism $\vartheta$ corresponds to the diagram automorphism which cyclically permutes $\alpha_1,\alpha_3,\alpha_4$ and fixes $\alpha_2$. It therefore fixes $(X_1-X_2)+(X_3-X_4)+(X_3+X_4)+2(X_2-X_3)=X_1+X_2$ and hence commutes with $e^{\ad t(X_1+X_2)}$ for all $t\in\R$. But $e^{\ad \pi(X_1+X_2)}=\Ad \tau_4$.
\end{proof}

It is now completely elementary to determine $\pi_0((\Aut\g)^\varrho)$ and representatives of its conjugacy classes  for classical $\g$. We do this case by case, using Cartan's classification (cf.~\cite{Hel}).

\vskip 1 cm

\textbf{A I\hskip 0,5 cm
 $\g=\su(n)\ (n\ge 3),\ \varrho=\mu,\ \k=\so(n)$}

$\varrho$ is outer since $n>2$. If $n$ is \textbf{odd} then $\pi_0((\Aut\g)^\varrho)=\{1,\varrho\}$ by Lemma~\ref{A.2}.

If $n$ is \textbf{even} then $(\Int\g)^\varrho=\{\Ad X \mid X\in U(n),\ \bar X=\lambda X$ for some $\lambda\in\C$ with $\vert\lambda\vert=1\}=\{\Ad X\mid X=aY\in U(n),\ \bar Y=Y,\ a\in\C\}=\{\Ad X\mid X\in O(n)\}$. It has two connected components as the restriction of $\Ad \tau_1\in(\Int\g)^\varrho$ to $\k$ is not inner. Thus $\id$, $\Ad \tau_1$ (resp.~$\id$, $\Ad \tau_1$, $\varrho$, $\varrho \Ad \tau_1)$ are representatives of the conjugacy classes of $\pi_0((\Int\g)^\varrho)\cong\Z_2$ (resp.~$\pi_0(\Aut\g)^\varrho\cong\Z_2\times\Z_2)$.

\vskip 1 cm

\textbf{
A II\hskip 0,5 cm
 $\g=\su(2n)\ (n\ge 2),\ \varrho=\mu \Ad J,\ \k=\sp(n)$}

Since $\k$ has no outer automorphisms $\pi_0((\Aut\g)^\varrho)=\{1,\varrho\}$ by Lemma~\ref{A.2}.

\vskip 1 cm

\textbf{
A III \hskip 0,5 cm
$ \g=\su(n)\ (n\ge 2),\ \varrho=\Ad \tau_p,\ \k={\mathfrak s} (\u(p)+\u(q))\ (p+q=n)$}

$\varrho$ is inner and commutes with $\mu$ which is outer if $n\ge 3$. Now, $(\Int\g)^\varrho=\{\Ad X\mid X=
\left({\scriptsize
\begin{array}[pos]{cc}
    A&B\\
    C&D
\end{array}}
\right)
\in U(n)$ with $A$ of size $p\times p$ and
$\left({\scriptsize
\begin{array}[pos]{cc}
    A&-B\\
    -C&D
\end{array}}
\right)
=\lambda
\left({\scriptsize
\begin{array}[pos]{cc}
    A&B\\
    C&D
\end{array}}
\right)
$ for some $\lambda\in\C\}$.
Necessarily $\lambda=1$ or $-1$ and $(\Int\g)^\varrho$ has one connected component if $p\ne q$
$(X=
\left({\scriptsize
\begin{array}[pos]{cc}
    O&B\\
    C&O
\end{array}}
\right)$
would be singular) and two if $p=q$ with non trivial representative $\Ad J$. Hence %$\pi_0((\Aut\g)^\varrho)\cong\{1,\mu\}\cong\Z_2$ if $p\ne q$ or $n=2$, and $\cong\{1, \Ad J,\mu, \mu \Ad &J\}\cong\Z_2\times\Z_2$ otherwise. 
\[
\pi_0((\Aut\g)^\varrho)\cong\left\{
\begin{array}{ll}
\{1,\mu\}\cong\Z_2&\mbox{if $p\ne q$ or $n=2$,} \\
\{1, \Ad J,\mu, \mu \Ad J\}\cong\Z_2\times\Z_2&\mbox{otherwise.}
\end{array}
\right.
\]
Note that $\Ad J=\mu$ if $n=2$.
\vskip 1 cm

{\bf
B I \hskip 0,5 cm
$\g=\so(2n+1)\ (n\ge 2),\ \varrho=\Ad \tau_p,\ \k=\so(p)+\so(q)\ (p+q=2n+1)$}

$\varrho$ is inner and from Lemma~\ref{A.2} and Remark~\ref{A.3} (if $p=2$ or $q=2)$ we have $\vert\pi_0((\Aut\g)^\varrho)\vert\le 2$. Since $\Ad \tau_1\tau_{p+1}$ induces  an outer automorphism on $\k,\ \pi_0((\Aut\g)^\varrho)\cong\{\id, \Ad \tau_1\tau_{p+1}\}\cong\Z_2$.

\vskip 1 cm

{\bf
C I \hskip 0,5 cm
$\g=\sp(n)\ (n\ge 3), \ \varrho=\Ad (iE),\ \k=\u(n)$}
(where $\sp(n)$ and $Sp(n)$ below are viewed as sets of quaternionic matrices)

Let $A_0=iE$. Then $(\Aut\g)^\varrho=(\Int\g)^\varrho=\{\Ad X\mid X\in Sp(n), A_0XA_0^{-1}=\pm X\}=\{\Ad X\mid X\in U(n)$ or $X\in j U(n)\}$. Thus $\pi_0 ((\Aut\g)^\varrho)\cong\{\id, \Ad jE\}\cong\Z_2$.

\vskip 1 cm

{\bf
C II \hskip 0,5 cm
$ \g=\sp(n)\ (n\ge 3),\ \varrho=\Ad \tau_p,\ \k=\sp(p)\times\sp(q)\ (p+q=n)$}

$\varrho$ is inner and $\Aut\k/\Int\k=1$ if $p\ne q$ and $\cong\Z_2$ if $p=q$. In the latter case $\Ad J\in(\Aut\g)^\varrho$ restricts to an outer automorphism on $\k$. Hence $\pi_0((\Aut\g)^\varrho)=1$ if $p\ne q$ and $\pi_0((\Aut\g)^\varrho)\cong\{\id, \Ad J\}\cong\Z_2$ if $p=q$.

\vskip 1 cm

{\bf
D I \hskip 0,5 cm
$\g=\so(2n)\ (n\ge 3),\ \varrho=\Ad \tau_p,\ \k=\so(p)+\so(q)\ (p+q=2n)$}

For simplicity we only consider the case $p=q=n$ with $n$ even. The other cases are similar but easier. Then $(\Aut\g)^\varrho\supseteq\{\Ad X\mid X\in O(2n),\ \tau_nX\tau_n=\pm X\}=\{\Ad X\mid X=
\left({\scriptsize
\begin{array}[pos]{cc}
 A& \\
    &B
\end{array}}
\right)$
or
$
\left({\scriptsize
\begin{array}[pos]{cc}
    & A\\
    B&
\end{array}}
\right)$
with $A,B\in O(n)\}$ and equality holds if $n\ne 4$. Moreover $(\Int\g)^\varrho$ consists precisely of those elements with $\det A=\det B$. Hence $(\Int\g)^\varrho$ has $4$ connected components represented e.g.~by $\Ad X_i$ with $X_1=
\left({\scriptsize
\begin{array}[pos]{cc}
    E&\\
    &E
\end{array}}
\right)$,
$X_2=
\left({\scriptsize
\begin{array}[pos]{cc}
    \tau_1&\\
    &\tau_1
\end{array}}
\right)$,
$X_3=
\left({\scriptsize
\begin{array}[pos]{cc}
    &E\\
    -E&
\end{array}}
\right)$,
and
$X_4=
\left({\scriptsize
\begin{array}[pos]{cc}
    &\tau_1\\
    -\tau_1&
\end{array}}
\right)$,
and $\pi_0((\Int\g)^\varrho)\cong\Z_2\times\Z_2$.

Let $\bf n>4$. Then $(\Aut\g)^\varrho$ has 8 connected components, represented by $\Ad X_i$ with $X_1,\dots,X_4$ as above, $X_5=
\left({\scriptsize
\begin{array}[pos]{cc}
    \tau_1&\\
&E
\end{array}}
\right)$,
$X_6=
\left({\scriptsize
\begin{array}[pos]{cc}
    E&\\
    &-\tau_1
\end{array}}
\right)$,
$X_7=
\left({\scriptsize
\begin{array}[pos]{cc}
    &\tau_1\\
E&
\end{array}}
\right)$,
and $X_8=
\left({\scriptsize
\begin{array}[pos]{cc}
    &E\\
    -\tau_1&
\end{array}}
\right)$,
which form a group isomorphic to $\D_4$. Moreover the $\Ad X_i$ with $i\in\{1,2,3,5,7\}$ are representatives of the conjugacy classes of this group and hence of $\pi_0((\Aut\g)^\varrho)$.

Now let $\bf n=4$. Since $\varrho$ is inner $\pi_0((\Aut\g)^\varrho)$ is isomorphic to a subgroup of $\Aut\k/\Int\k$ by Lemma~\ref{A.2}, which in turn is isomorphic to the group of symmetries of the Dynkin diagram of $\so(4)+\so(4)$ and thus to $S_4$ (the symmetric group).
Moreover $\vartheta\in(\Aut\g)^\varrho$ by Lemma~\ref{A.4} and hence $\pi_0((\Aut\g)^\varrho)\cong S_4$. For
$\vert\pi_0((\Int\g)^\varrho)\vert=4$ while $(\Aut\g)^\varrho/(\Int\g)^\varrho$ is isomorphic to a subgroup of $\Aut\g/\Int\g\cong S_3$ by Lemma~\ref{A.1} and contains elements of order two and three and is thus isomorphic to $S_3$. The conjugacy classes of $S_4$ consist of the sets of cycles of order 1 to 4 and $Z:=\{(1,2)(3,4),(1,3)(2,4),(1,4)(2,3)\}$ of cardinality 1,6,8,6, and 3, respectively. Therefore $\pi_0((\Int\g)^\varrho)$, which is a normal subgroup of $\pi_0((\Aut\g)^\varrho)$ with four elements, corresponds to $\{(1)\}\cup Z$ and any two non trivial elements of $\pi_0((\Int\g)^\varrho)$ are conjugate in $\pi_0((\Aut\g)^\varrho)$. Thus $\vartheta$ and the $\Ad X_i$ above with $i\in\{1,3,5,7\}$ represent the conjugacy classes of $\pi_0((\Aut\g)^\varrho)$.

\vskip 1 cm
{\bf
D III \hskip 0,5 cm
$\g=\so(2n),\ \varrho=\Ad J,\ \k=\u(n)$}

In this case $(\Aut\g)^\varrho\supseteq\{\Ad X\mid X\in O(2n), JXJ^{-1}=\pm X\}=\{\Ad X\mid X\in U(n)\}\cup\{\Ad X\mid X\in\tau_n U(n)\}$ with equality if $n\ne 4$. But equality also holds if $n=4$ as $\vert\pi_0((\Aut\g)^\varrho)\vert\le 2$ by Lemma~\ref{A.2}. Thus $(\Aut\g)^\varrho$ has two connected components represented e.g.~by $\id$ and $\Ad \tau_n$ and $\Ad \tau_n$ is inner precisely if $n$ is even.

\vskip 1 cm

We now study the exceptional case. Our main tool here is the following result (cf.~\cite{Hel}, chapter VII, 7.2).

\begin{theorem}\label{A.5}
The fixed point set $G^\varrho$ of an involution $\varrho$ on a compact, connected and simply connected Lie group $G$ is connected. \qed
\end{theorem}

\begin{corollary}\label{A.6}
Let $G$ be the universal cover of $\Int\g$. Then $\pi_0((\Int\g)^\varrho)$ is isomorphic to a quotient of a subgroup of the center of $G$.
\end{corollary}

\begin{proof}
Let $p:G\to \Int\g$ be the universal covering and $\hat G:=p^{-1}((\Int\g)^\varrho)=\{g\in G\mid \varrho(g)g^{-1}\in Z(G)\}$ where $Z(G)$ is the center of $G$. The mapping $\hat G\to Z(G), g\mapsto\varrho(g)g^{-1}$, is a homomorphism and induces an exact sequence
\[
1\to G^\varrho\to\hat G\to Z(G)\ .
\]
Since $G^\varrho$ is connected by~\ref{A.5}, $\pi_0(\hat G)=\hat G/G^\varrho$ is isomorphic to a subgroup of $Z(G)$. Moreover $\pi_0((\Int\g)^\varrho)$ is isomorphic to a quotient of this subgroup as $\hat G\to(\Int\g)^\varrho$ is surjective.
\end{proof}

From this and Cartan's computations of the centers (cf.~\cite{Hel}, Table IV)
we immediately get $\pi_0((\Int\g)^\varrho)\cong 1$ or $\Z_3$ if $\g=\e_6,\ \pi_0((\Aut\g)^\varrho)=\pi_0 ((\Int\g)^\varrho)\cong1$ or $\Z_2$ if $\g=\e_7$ and $\pi_0((\Aut\g)^\varrho)=\pi_0((\Int\g)^\varrho)=1$ if $\g=\e_8,\f_4$ or $\g_2$. By the next lemma the $\Z_3$ in case of $\g=\e_6$ is excluded.

\begin{lemma}\label{A.7}
$\pi_0((\Int\g)^\varrho)\cong(\Z_2)^l$ for some $l\ge 0$.
\end{lemma}

\begin{proof}
Let $K:=\{e^{\ad X}:\g\to\g\mid X\in\k\}$. Then $K=((\Int\g)^\varrho)_0,\ \Int\g=K\cdot\{e^{\ad X}\mid X\in\p\}$ (as the exponential mapping of $\Int\g/K$ is surjective), and thus $(\Int\g)^\varrho=K\cdot\{e^{\ad X}\mid X\in\p\}^\varrho$. If $X\in\p$ then $e^{\ad X}\in(\Int\g)^\varrho$ if and only if $e^{-\ad X}=e^{\ad X}$, i.e. $(e^{\ad X})^2=\id$. Hence each non trivial element of $\pi_0((\Int\g)^\varrho)$ is of order $2$.
\end{proof}

\begin{lemma}\label{A.8}
The conjugacy classes of involutions on $\g=\e_6$ can be represented by commuting elements.
\end{lemma}

\begin{proof}
Up to conjugacy there are four involutions on $\e_6$, two inner and two outer. If $\rho$ is an outer involution and ${\mathfrak t}_0$ a maximal abelian subalgebra of $(\e_6)^\varrho$ then the second outer involution may be chosen to be of the form $\rho\cdot \Ad e^{X}$ with $X\in{\mathfrak t}_0$ (cf.~\cite{Loo}, Theorem 3.3, Chapter VII, \cite{Wol} Theorem 8.6.9) and hence to commute with $\rho$. Therefore it suffices to prove the existence of two non-conjugate involutions of the form $\Ad e^{X}$ with $X\in \mathfrak t_0$.

We choose $\rho$ to be the diagram automorphism of $\e_6$, after fixing a maximal torus $\mathfrak t$ of $\e_6$ and a basis, say $\alpha_1,\dots,\alpha_6$ of the root system. In general, if $\sum m_i\alpha_i$ is the maximal root and $X_1,\dots, X_n$ are the elements in $\mathfrak t$ with $\alpha_i(X_j)={1\over m_i}\vartheta_{ij}$ then $\Ad e^{X}$ with $X=X_i$ and $m_i=2$ or $X={1\over 2}(X_i+X_j)$ and $m_i=m_j=1,i\ne j$, are involutions. Moreover, involutions corresponding to the first case $(m_i=2)$ are not conjugate to those of the second $(m_i=m_j=1)$ as their fixed point algebras are semisimple and not semisimple, respectively (cf. \cite{Loo}, p. 121 - 123). Hence the desired result follows from the diagram

%%%%%%%%%%%%%%%%%%%%%%%%%%%%%%%%%%%%
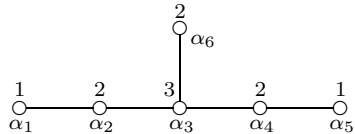
\begin{center}\begin{picture}(120,46)

          \multiput(4,8)(30,0){2}{\circle{5}}
          \multiput(94,8)(30,0){2}{\circle{5}}
           \put(64,8){\circle{5}}
           \put(64,38){\circle{5}}

           \put(6.5,8){\line(1,0){25}}
           \put(36.5,8){\line(1,0){25}}
           \put(66.5,8){\line(1,0){25}}
           \put(96.5,8){\line(1,0){25}}
           \put(64,10.5){\line(0,1){25}}

           \put(0,0){$\scriptstyle \alpha_1$}
           \put(30,0){$\scriptstyle \alpha_2$}
           \put(60,0){$\scriptstyle \alpha_3$}
           \put(90,0){$\scriptstyle \alpha_4$}
           \put(120,0){$\scriptstyle \alpha_5$}
           \put(68,32){$\scriptstyle \alpha_6$}

           \put(62,42){$\scriptstyle 2$}

           \put(2,12.5){$\scriptstyle 1$}
           \put(32,12.5){$\scriptstyle 2$}
           \put(58,12.5){$\scriptstyle 3$}
           \put(92,12.5){$\scriptstyle 2$}
           \put(122,12.5){$\scriptstyle 1$}

   \end{picture}\end{center}

%%%%%%%%%%%%%%%%%%%%%%%%%%%%%%%%%%%%

of $\e_6$, in which the superscripts denote the numbers $m_i$. Note that $X_6$ and $X_1+X_5$ are contained in $\mathfrak t^\rho$.
\end{proof}

From Lemmas~\ref{A.7} and \ref{A.8} we get $\pi_0((\Aut \e_6)^{\varrho_i})\cong\{1,\varrho_1\}\cong\Z_2$ for all $i\in\{1,\dots,4\}$ where $\varrho_1$ is an outer automorphism and $\varrho_1,\dots,\varrho_4$ are commuting representatives of the conjugacy classes of involutions on $\e_6$.

The most intricate case is $\g=\e_7$ which we consider now. Up to conjugacy $\g$ has three involutions $\varrho_1,\varrho_2,\varrho_3$ with fixed point algebras (i) $\su(8)$, (ii) $\so(12)+\su(2)$ and (iii) $\e_6+\R$, respectively. Since $\e_7$ has no outer automorphisms $\pi_0((\Aut\g)^\varrho)=\pi_0((\Int\g)^\varrho)$.

\begin{lemma}\label{A.9}
In case (ii), $\pi_0((\Aut\g)^\varrho)=1$.
\end{lemma}

\begin{proof}
Let $\varphi\in(\Aut\g)^\varrho\setminus ((\Aut\g)^\varrho)_0$. Thus $\varphi$ restricted to $\k=\so(12)+\su(2)$ is outer by Lemma~\ref{A.2} and after multiplying it by an appropriate element from $((\Aut\g)^\varrho)_0=\{e^{\ad X}\mid X\in\k\}$ we may assume that $\varphi$ is the standard involution on $\k=\so(12)+\su(2)$ with $\k^\varphi=\so(11)+\su(2)$ and $\dim\k^\varphi=58$. Then $\varphi$ is also an involution on $\g$. For $\varphi^2=\id$ on $\k$ implies that $\varphi^2_{\vert_\p}$ commutes with the $\k$-action on $\p$, which is irreducible. Therefore $\varphi^2_{\vert_\p}=\id$ (as claimed) or $\varphi^2(X)\ne X$ for all $X\in\p\setminus\{0\}$. But in the latter case $\k^\varphi=\g^\varphi$ which is in contradiction to $\rank\ \k^\varphi<\rank\ \k$. Let $\p_\pm$ be the $\pm 1$ eigenspaces of $\varphi$ on $\p$. Then the fixed point algebras $\k^\varphi+\p_\pm$ of $\varphi$ and $\varrho\varphi$ are symmetric and thus of dimension 63, 69 or 79. Hence $\dim\p_\pm\in\{5,11,21\}$ in contradiction to $\dim\p_++\dim\p_-=\dim\p=64$.
\end{proof}

The other two cases of $\e_7$ are handled by the next two results.

\begin{lemma}\label{A.10}
Let  $\varrho$ be an inner involution on $\g$ and $X\in\p\setminus\{0\}$ with $(\ad X)^3=-\pi^2 \ad X$. Then $e^{\ad X}$ represents a non trivial element of $\pi_0((\Aut\g)^\varrho)$.
\end{lemma}

\begin{proof}
The eigenvalues of $\ad X$ are $0$ and $\pm\ i\pi$. Therefore $\varphi:=e^{\ad X}$ has eigenvalues $\pm 1$ and $\ker(\ad X)$ is the eigenspace of $1$. Thus $\varphi$ is an involution and hence commutes with $\varrho$. Since $[X,\k^\varphi]=0$ but $X\not\in\k^\varphi,\ \rank\ \k^\varphi<\rank\ \g=\rank\ \k$ and the restriction of $\varphi$ to $\k$ is not inner. Thus $\varphi$ represents by Lemma~\ref{A.2} a non trivial element of $\pi_0((\Aut\g)^\varrho).$
\end{proof}

The $X\in\p$ with $(\ad X)^3=-\pi^2 \ad X$ are strongly related to extrinsic symmetric spaces in the sense of Ferus (\cite{Fer}, cf. also \cite{EH}) and it is known that their existence can be read off from the highest root of the symmetric space. More precisely, let $\Sigma$ be the restricted root system of $(\g,\k)$ with respect to a maximal abelian subspace $\a$ of $\p$, ${\mathfrak C}\subset \a$ a Weyl chamber, and  $\alpha_1,\dots,\alpha_r$ the corresponding basis of $\Sigma$. Up to conjugation we may assume $X\in{\bar{\mathfrak C}}$. On each root space $\g_\alpha=\{Y\in\g_{\C}\mid [H,Y]=\pi i\alpha(H)Y\ \forall\ H\in\a\},\ \ad X$ has eigenvalues $\pi i\alpha(X)$. Hence $(\ad X)^3=-\pi^2 \ad X$ is equivalent to $\alpha(X)\in\{0,\pm 1\}$. Let $X_1,\dots,X_r\in\a$ be dual to $\alpha_1,\dots,\alpha_r$ and $X=\sum x_iX_i$. Then $x_i\ge 0$ (because $X\in\bar{\mathfrak C})$ and $(\ad X)^3=-\pi^2 \ad X$ is equivalent to $x_i\in\{0,1\}$ (because $\alpha_i(X)\in\{0,\pm 1\})$ and $\sum m_ix_i\in\{0,1\}$ where $\delta=\sum m_i\alpha_i$ is the maximal root. Hence we have the following result.

\begin{lemma}\label{A.11}
A non trivial $X\in\p$ with $(\ad X)^3=-\pi^2 \ad X$ exists if and only if at least one of the coefficients $m_i$ of the highest root $\delta$ of $\Sigma$ is equal to one. \qed
\end{lemma}

Now in cases (i) and (iii) of $\g=\e_7,\ \Sigma$ is of type $\e_7$ and $\c_3$, respectively, and in both cases the highest roots have a coefficient $m_i=1$ (see e.g.~\cite{Hel}, ch.~X). Hence $\pi_0((\Aut\g)^\varrho)=\pi_0((\Int\g)^\varrho)\cong\Z_2$ by Lemmas~\ref{A.7} and \ref{A.11}. Moreover any $e^{\ad X}$ with $X\in\p,\ X\ne 0$ and $(\ad X)^3=-\pi^2\ \ad X$ represents the non trivial element.

\newpage

\section{Conjugate linear automorphisms of $\g$}

We recall and outline here some results about finite dimensional complex simple Lie algebras, in particular the relation between their complex linear and conjugate linear automorphisms.

Thus let $\g$ be a complex simple Lie algebra, $\u$ a compact real form of $\g$ and $\omega$ the conjugation with respect to $\u$. Let $\g^\R$ be the realification of $\g$. Then $\u$ is a maximal compact subalgebra of $\g^\R$ and $\u+i\u$ is a Cartan decomposition of $\g^\R$.  Moreover $\Aut\g^\R=\Aut\g\cup{\overline{\Aut}}\g$, where ${\overline\Aut}\g:=\omega\Aut\g$ denotes the set of conjugate linear automorphisms of $\g$ (since $(\g^\R)_\C$ is the sum of two simple ideals $\g_\pm$ which are either left invariant or interchanged by $\varphi_\C,\ \forall\ \varphi\in\Aut\g^\R)$.

Let $\Aut_\u\g^\R:=\{\varphi\in\Aut\g^\R\mid \varphi\u=\u\}$. Then, by classical results of Cartan (cf. e.g. \cite{GOV} Ch. 4, 3.2 and 3.3), the following holds
\begin{enumerate}
\item[(A)] The mapping $\Aut_\u\g^\R\times i\u\to\Aut\g^\R,\ (\varphi,X)\mapsto\varphi e^{\ad X}$, is a diffeomorphism.
\item[(B)] Each compact subgroup of $\Aut\g^\R$ is conjugate to a subgroup of $\Aut_\u\g^\R$ (and the conjugation can be achieved by an element from $\Int\g$ by $(A))$. 
\end{enumerate}
In fact, (A) follows from the Hadamard-Cartan Theorem applied to the symmetric space of nonpositive curvature $\Aut\g^\R/\Aut_\u\g^\R$ while (B) follows from the Cartan fixed point Theorem. As a consequence of (A) we have
\begin{enumerate}
\item[(C)] If $\varphi,\tilde\varphi\in\Aut_\u\g^\R$ are conjugate (resp. conjugate by an inner automorphism) in $\Aut\g^\R$ then they are already conjugate (by an inner automorphism) in $\Aut_\u\g^\R$.
\end{enumerate}

Indeed, if $\tilde\varphi=\alpha\varphi\alpha^{-1}$ for some $\alpha=\alpha_0e^{\ad X}$ with $\alpha_0\in\Aut_\u\g^\R$ and $X\in i\u$ then $\alpha_0^{-1}\tilde\varphi\alpha_0e^{\ad X}=\varphi e^{\ad \varphi^{-1}X}$ and thus $\alpha_0^{-1}\tilde\varphi\alpha_0=\varphi$ by (A).

We denote the sets of automorphisms of order $q$ by $\Aut^q$. The mapping $\Aut^q\u\to\Aut^q\g$ that maps $\varphi$ to its complex linear extension $\varphi_\C$ induces mappings
\[ \begin{array}{lcl}
\Aut^q\u/\Aut\u&\to&\Aut^q\g/\Aut\g \quad \mbox { and}\\
\Aut^q\u/\Int\u&\to&\Aut^q\g/\Int\g
\end{array}
\]
between the conjugacy classes. It follows from (B) that these mappings are surjective and from (C) that they are injective.

If instead of $\varphi_\C$ the conjugate linear extension $\varphi_\C\omega$ of $\varphi\in\Aut\u$ is used one gets mappings
\[
\Aut^q\u/\Aut\u (\Int\u) \to
\left\{
\begin{array}{lrl}
\overline{\Aut}\,^q\g/\Aut\g(\Int\g)&q& \mbox{ even}\\
\overline{\Aut}\,^{2q}\g/\Aut\g(\Int\g) &q& \mbox{ odd}
\end{array}
\right.
\]
which are by the same reasoning injective. By (B) they are also surjective if $q$ is even and induce a surjective and hence bijective mapping
\[
(\Aut^q\u\cup\Aut^{2q}\u)/\Aut\u(\Int\u)\to \overline{\Aut}^{2q}\g/\Aut\g(\Int\u)
\]
if $q$ is odd. Note that the order of $\varphi\in\Aut\u$ is equal to the order of $\varphi_\C\omega$ or to half of it depending on whether the latter is divisible by $4$ or not.

\bigskip
Summarizing, we have.

\begin{theorem}
The following mappings which are induced by complex linear, resp. conjugate linear extensions are bijective:
\begin{enumerate}
\item[(i)] $\Aut^q\u/\Aut\u(\u)\to\Aut^q\g/\Aut\g(\g)\ ,\ q\in\N$
\item[(ii)] $\Aut^{2q}\u/\Aut\u(\u)\to\overline{\Aut}^{2q}\g/\Aut\g(\g)$, $q$ even
\item[(iii)] $(\Aut^q\u\cup\Aut^{2q}\u)/\Aut\u(\u)\to\overline{\Aut}^q(\g)/\Aut\g(\g)$, $q$ odd. 
\end{enumerate}
$\Aut\u$ and $\Aut\g$ in the denominators may be replaced by $\Int\u$ and $\Int\g$, respectively.
\qed
\end{theorem}

In particular, $q=2$ gives the well known bijections
\[
\Aut^2\g/\Aut\g\leftrightarrow\Aut^2\u/\Aut\u\leftrightarrow\{\mbox{non compact real forms of }\g\}/\Aut\g\ .
\]
Note that $\Aut^1\u=\{id\}$ corresponds to the compact real forms. If $\varrho\in\Aut^2\u$ is an involution then the corresponding real form is the fixed point set of $\varrho_\C\omega$ and thus $\k+i\p$, where $\k$ and $\p$ are the $\pm 1$ eigenspaces of $\varrho$.

\begin{proposition}
Let $\varrho\in\Aut\,\u$. Complex linear resp.~conjugate linear extensions induce bijections between the conjugacy classes of the following (subsets of) groups
\begin{enumerate}
\item[(i)] $\pi_0((\Aut\,\u)^\varrho),\ \pi_0((\Aut\g)^{\varrho_\C})$, and $\pi_0((\Aut\g)^{\varrho_\C\omega})$
\item[(ii)] $\pi_0((\Aut\,\u)^\varrho),\ \pi_0((\overline{\Aut}\g)^{\varrho_\C})$, and $\pi_0((\overline{\Aut}\g)^{\varrho_\C\omega})$.
\end{enumerate}
\end{proposition}

\begin{proof}
(A) implies $(\Aut\g^\R)^{\varrho_\C}=(\Aut_\u\g^\R)^{\varrho_\C}\cdot\{e^{\ad X}\mid X\in i\k\}$ and $(\Aut\g^\R)^{\varrho_\C\omega}=(\Aut_\u\g^\R)^{\varrho_\C}\cdot\{e^{\ad X}\mid X\in i\p\}$. Thus the proposition follows from the isomorphism $\{id,\omega\}\times(\Aut\u)^\varrho\to(\Aut_\u\g^\R)^{\varrho_\C)}\ ,\ (\alpha,\beta)\mapsto\alpha\circ\beta_\C$\ .
\end{proof}

Here $\pi_0((\overline{\Aut}\g)^{\varrho_\C})$, for example, is considered as a subset of $\pi_0((\Aut\g^\R)^{\varrho_\C})$, the set of connected components of $(\Aut\g^\R)^{\varrho_\C}$.

\end{appendix}

%%%%%%%%%%%%%%%%%%%%%%%

\vskip0,5cm
Christian Gro{\ss}\\
Ernst Heintze\\
Institut für Mathematik, Universität Augsburg\\
Universitätsstrasse 14\\
D - 86159 Augsburg, Germany\\
e-mail: {\it gross@math.uni-augsburg.de}\\
\hspace*{1,3cm} {\it heintze@math.uni-augsburg.de}\\

\end{document}